\pdfoutput=1
\RequirePackage{ifpdf}
\ifpdf 
\documentclass[pdftex]{sigma}
\else
\documentclass{sigma}
\fi

\usepackage{mathtools}
\usepackage[arrow,matrix,curve,color]{xy}
\usepackage{pb-diagram,pb-xy}
\usepackage{enumitem}
\usepackage[normalem]{ulem}
\usepackage{tikz}
\usepackage{xspace}
\definecolor{light-blue}{rgb}{0.8,0.85,1}
\definecolor{light-red}{rgb}{1,.4,.4}
\definecolor{purp}{rgb}{.7,.3,1}
\definecolor{yel}{rgb}{1,1,.5}
\definecolor{cy}{rgb}{0,1,1}
\definecolor{m}{rgb}{1,0.1,1}

\usepackage{colortbl}
\usepackage{multirow}
\usepackage{array}

\numberwithin{equation}{section}

\newtheorem{Theorem}{Theorem}[section]
\newtheorem{Corollary}[Theorem]{Corollary}
\newtheorem{Lemma}[Theorem]{Lemma}
\newtheorem{Proposition}[Theorem]{Proposition}
 { \theoremstyle{definition}
\newtheorem{Definition}[Theorem]{Definition}
\newtheorem{Example}[Theorem]{Example}
\newtheorem{Remark}[Theorem]{Remark}}

\newcommand{\p}{\partial}
\newcommand{\fb}{\mbox{-}\mathrm{fb}}

\newcommand{\spin}{\mathsf{spin}}

\newcommand{\scal}{{\rm scal}}
\newcommand{\cyl}{{\rm cyl}}

\newcommand{\psc}{positive scalar curvature}
\newcommand{\co}{\colon}

\newcommand{\bR}{\mathbb R}

\newcommand{\bH}{\mathbb H}

\newcommand{\bZ}{\mathbb Z}

\newcommand{\bP}{\mathbb P}

\newcommand{\bQ}{\mathbb Q}

\newcommand{\cR}{\mathcal R}

\newcommand{\SU}{\mathop{\rm SU}}

\newcommand{\Sp}{\mathop{\rm Sp}}
\newcommand{\Spin}{\mathop{\rm Spin}}

\newcommand{\lp}{\textup{(}}
\newcommand{\rp}{\textup{)}}

\newcommand{\res}{\operatorname{res}}
\newcommand{\resS}{\res_\Sigma}

\newcommand{\Bott}{\mathrm{Bott}}

\newcommand{\bord}{\rightsquigarrow}
\begin{document}

\newcommand{\arXivNumber}{2005.02744}

\renewcommand{\thefootnote}{}

\renewcommand{\PaperNumber}{062}

\FirstPageHeading

\ShortArticleName{Positive Scalar Curvature on Spin Pseudomanifolds}

\ArticleName{Positive Scalar Curvature on Spin Pseudomanifolds:\\ the Fundamental Group and Secondary Invariants\footnote{This paper is a~contribution to the Special Issue on Scalar and Ricci Curvature in honor of Misha Gromov on his 75th Birthday. The~full collection is available at \href{https://www.emis.de/journals/SIGMA/Gromov.html}{https://www.emis.de/journals/SIGMA/Gromov.html}}}

\Author{Boris BOTVINNIK~$^{\rm a}$, Paolo PIAZZA~$^{\rm b}$ and Jonathan ROSENBERG~$^{\rm c}$}

\AuthorNameForHeading{B.~Botvinnik, P.~Piazza and J.~Rosenberg}

\Address{$^{\rm a)}$~Department of Mathematics, University of Oregon, Eugene OR 97403-1222, USA}
\EmailD{\href{mailto:botvinn@uoregon.edu}{botvinn@uoregon.edu}}
\URLaddressD{\url{http://pages.uoregon.edu/botvinn/}}

\Address{$^{\rm b)}$~Dipartimento di Matematica ``Guido Castelnuovo'', Sapienza Universit\`a di Roma,\\
\hphantom{$^{\rm b)}$}~Piazzale Aldo Moro 5, 00185 Roma, Italy}
\EmailD{\href{mailto:piazza@mat.uniroma1.it}{piazza@mat.uniroma1.it}}
\URLaddressD{\url{http://www1.mat.uniroma1.it/people/piazza/}}

\Address{$^{\rm c)}$~Department of Mathematics, University of Maryland, College Park, MD 20742-4015, USA}
\EmailD{\href{mailto:jmr@umd.edu}{jmr@umd.edu}}
\URLaddressD{\url{http://www2.math.umd.edu/~jmr/}}

\ArticleDates{Received May 26, 2020, in final form June 08, 2021; Published online June 24, 2021}

\Abstract{In this paper we continue the study of positive scalar curvature (psc) metrics on~a~depth-1 Thom--Mather stratified space $M_\Sigma$ with singular stratum $\beta M$ (a closed manifold of~positive codimension) and associated link equal to $L$, a smooth compact mani\-fold. We~briefly call such spaces manifolds with $L$-fibered singularities. Under suitable spin assumptions we give necessary index-theoretic conditions for the existence of wedge metrics of positive scalar curvature. Assuming in addition that $L$ is a simply connected homogeneous space of positive scalar curvature, $L=G/H$, with the semisimple compact Lie group~$G$ acting transitively on~$L$ by isometries, we investigate when these necessary conditions are also sufficient. Our main result is that our conditions are indeed sufficient for large classes of~examples, even when $M_\Sigma$ and $\beta M$ are not simply connected. We~also investigate the space of such psc~metrics and show that it often splits into many cobordism classes.}


\Classification{53C21; 58J22; 53C27; 19L41; 55N22; 58J28}

\renewcommand{\thefootnote}{\arabic{footnote}}
\setcounter{footnote}{0}

\section{Introduction}
This paper continues a program begun in~\cite{MR1857524,BR},
and our previous paper~\cite{BB-PP-JR} (hereafter called \emph{Part~I}), to understand obstructions to positive scalar curvature on
manifolds with fibered singularities, for metrics that are well
adapted to the singularity structure.

As in~\cite{BB-PP-JR}, the stratified spaces (or singular manifolds)
$M_\Sigma$ that we study are Thom--Mather pseudomanifolds of depth one,
where we take the two strata to be spin. Topologically, $M_\Sigma$~is
homeomorphic to a quotient space of a compact smooth manifold $M$ with
boundary $\partial M$. The~manifold $M$ is called
the \emph{resolution} of $M_\Sigma$, and the quotient map $M\to
M_\Sigma$ is the identity on~the interior $\mathring M$ of $M$, and on
$\partial M$, collapses the fibers of a fiber bundle
$\phi\co\partial M\to \beta M$, with fibers all diffeomorphic to a
fixed manifold $L$, called the \emph{link} of the singularity, and
with base~$\beta M$ sometimes called the Bockstein of
$M$ (by analogy with other cases in topology). While we do treat
the general case for some of our results, we shall eventually adopt
the main geometric assumptions from~\cite{BB-PP-JR}, namely that the
bundle $\phi$ comes from a principal $G$-bundle $p\co P\to \beta
M$, for some connected semisimple compact Lie group $G$ that acts
transitively on~$L$ by isometries for some fixed metric $g_L$, and
thus $\partial M = P\times_G L$. The transitivity of the action of~$G$
on~$L$ means that $L=G/H$ is a homogeneous space
which comes with a metric~$g_L$ with constant
positive scalar curvature, which we normalize to be the same as that
of a sphere of the dimension $\ell = \dim L$.

Given a Lie group $G$ and a link $L=G/H$ as above, we
say that a fiber bundle $E\to B$ with a~fiber~$L$ is a
\emph{geometric $(L,G)$-bundle} if its structure group is $G$, where
$G$ acts on $L$ by~iso\-metries of the metric $g_L$.

An \emph{adapted} metric $g$ on $M_\Sigma$ will be defined by the
following data: a Riemannian metric~$g_M$ on $M$, which is a product
metric ${\rm d}t^2 + g_{\partial M}$ in a small collar neighborhood
$\partial M\times [0, \varepsilon)$ of the boundary $\partial M$, a
 connection $\nabla^p$ on the principal $G$-bundle $p$, and a
 Riemannian metric $g_{\beta M}$ on $\beta M$. We~require the
 metrics $g_L$, $g_{\p M}$ and $g_{\beta M}$ to be compatible in the
 sense that the bundle projection $\phi\co (\partial M,g_{\p M}) \to
 (\beta M, g_{\beta M})$ is a Riemannian submersion with fibers $(L,
 g_L)$. Furthermore, since the structure group $G$ of the bundle
 $\phi\co \partial M\to \beta M$ acts by isometries of the metric
 $g_L$, we can make the special metric $g_L$ on the fibers orthogonal
 to the metric~$g_{\beta M}$ lifted up to the horizontal spaces for
 the connection $\nabla^p$. Then $M_\Sigma$ is the result of gluing
 together the Riemannian manifold $M$ and a tubular neighborhood $N$
 of $\beta M$ along their common boundary~$\partial M$. The
 complement of $\beta M$ in $N$ will look like a fiber bundle over
 $\beta M$ whose fibers are open cones $(0, R)\times L$, where $R$ is
 the radius of the cones. We~require these fibers to be actual
 metric cones with the metric ${\rm d}r^2 + r^2g_L$ (which is actually a
 warped product metric on the product of $L$ with the interval
 $(0,R)$), transitioning smoothly near $r=R$ to a product metric
 ${\rm d}r^2 + Cg_L$ for a suitable positive constant $C$. In this paper,
 we allow for $M$ and $\beta M$ to have non-trivial fundamental
 groups; the simply-connected case was considered in Part~I~\cite{BB-PP-JR}. Throughout this paper we assume that the link $L =
 G/H$ is a simply connected homogeneous space, where $G$ is a compact
 semisimple Lie group acting on $L$ by isometries of the metric
 $g_L$, and where
\begin{gather}\label{normalization}
\scal_{g_L}=\scal_{S^{\ell}}=\ell(\ell-1) .
\end{gather}
This normalization
makes the cone on $L$ scalar-flat with respect to the metric ${\rm d}r^2 + r^2g_L$.
This will be important for our eventual existence results since
this will guarantee that if $g_{\beta M}$ has {\psc}, then we can
make the scalar curvature uniformly positive on the tubular neighborhood
of $\beta M$.
In this setting, it is easy to see (since
$L$ is simply connected) that
$\pi_1(\p M) = \pi_1(\beta M)$, and then Van Kampen
implies that $\pi_1(M_\Sigma)=\pi_1(M)$.
\begin{Theorem}[obstruction theorem]\label{thm:psobstruction}
Let $M_\Sigma$ be an $n$-dimensional compact pseudomanifold with
resolution $M$, a spin manifold with boundary $\p M$. Assume the following:
\begin{enumerate}\itemsep=0pt
\item[{\rm (1)}] $M$ is a spin manifold with boundary $\p M$
 fibered over a connected spin manifold $\beta M$,
\item[{\rm (2)}] the fiber bundle $\phi \co\p M\to \beta M$ is a
 geometric $(L,G)$-bundle.
\end{enumerate}
Let $\Gamma = \pi_1(M)$, $\Gamma_\beta = \pi_1(\beta M)$. Assume $g$ is
an adapted Riemannian metric on $M_{\Sigma}$ and let~$g_M$ be the restriction
of $g$ to $M$. Then there are two
``alpha invariants'', generalized indices of~Clifford algebra-linear
Dirac operators,
\begin{gather*}
\alpha^\Gamma_{\rm cyl}(M,g_M)\in KO_n\big(C^*_{r,\bR} (\Gamma)\big) \qquad\text{and}\qquad \alpha^{\Gamma_\beta}(\beta M) \in KO_n\big(C^*_{r,\bR} (\Gamma_\beta)\big).
\end{gather*}
The indices $\alpha^\Gamma_{\rm cyl}(M,g_M)$ and
$\alpha^{\Gamma_\beta}(\beta M)$ do not depend on a choice of adapted
metric $g$, and they both vanish if there exists an adapted psc~metric
on $M_{\Sigma}$.
\end{Theorem}

{\sloppy\begin{Remark}
 In our case we fix the metric $g_L$ on $L$ and the
 connection on the bundle \mbox{$\phi \co\p M \to \beta M$} coming
 from a connection $\nabla^p$ on the principal bundle
 $p\co P\to \beta M$. (This is harmless since the space of such
 connections is contractible.) Then an adapted
 wedge metric $g$ on $M_{\Sigma}$ is
 determined (up to contractible choices) by a metric $g_{\beta M}$
 (which then determines the metric $g_{\p M}$) and by an extension
 $g$ of $g_{\p M}$ to $M$. Thus the space of adapted
 wedge metrics on $M_{\Sigma}$ is
 contractible, and from the analytic properties of Dirac
 operators we then obtain that the cylindrical $\alpha$-class
 $\alpha_{\text{cyl}}(M,g_M)$ is
 independent of $g$ (once the metric $g_{\beta M}$
 has been fixed). Hereafter, we will omit the metric $g$ from the
 notation. Notice that there is also a \emph{wedge} $\alpha$-class
 $\alpha^\Gamma_w (M_\Sigma,g)$ defined by considering the Dirac
 operator on the regular part of the pseudomanifold~$M_\Sigma$. This
 is also independent of $g$ for $(L,G)$-fibered pseudomanifolds. The
 two classes $\alpha^\Gamma_{\rm cyl}(M)$ and
 $\alpha^\Gamma_w (M_\Sigma)$ are equal if the metric $g_{\beta M}$ is psc.
 However, in more general situations, the
 wedge $\alpha$-class $\alpha^\Gamma_w (M_\Sigma,g)$
 and the cylindrical $\alpha$-class
 $\alpha^\Gamma_{\rm cyl}(M,g_M)$ both depend on the
 choice of metric $g$ and they are in general different.
 We shall make all this very precise in the next section.
\end{Remark}

}

Now we are ready for the existence result.
\begin{Theorem}[existence theorem]
Let $M_\Sigma$ be an $n$-dimensional compact pseudomanifold with
resolution $M$, a spin manifold with boundary $\p M$. Assume the following:
\begin{enumerate}\itemsep=0pt
\item[{\rm (1)}] $M$ is a spin manifold with boundary $\p M$
 fibered over a connected spin manifold $\beta M$,
\item[{\rm (2)}] the fiber bundle $\phi\co \p M\to \beta M$ is a
 geometric $(L,G)$-bundle.
\end{enumerate}
Let $\Gamma = \pi_1(M)$, $\Gamma_\beta = \pi_1(\beta M)$. Furthermore,
assume $n\geq \ell+6$, where $\ell=\dim L$, and that one of the
following condition holds:
\begin{enumerate}\itemsep=0pt
\item[$(i)$] either $L$ is the boundary of a
 spin $G$-manifold
 $\bar L$ equipped with a $G$-invariant psc~metric~$g_{\bar L}$,
 which is a product near the boundary and satisfies $g_{\bar L}|_L=g_L$,
\item[$(ii)$] or the embedding $\p M \to M$ induces an
 isomorphism on $\pi_1$, and moreover, $\p M=\beta M\times L$,
 where $L$ is an even quaternionic projective space, and
 $\Omega^\spin_*(B\Gamma)$ is free as an $\Omega^\spin_*$-module.
\end{enumerate}
Then, provided the Gromov--Lawson--Rosenberg conjecture holds for the groups
$\Gamma$ and $\Gamma_{\beta}$, the vanishing of
the invariants
\begin{gather*}
\alpha^\Gamma_{\rm cyl}(M_\Sigma)\in KO_n\big(C^*_{r,\bR} (\Gamma)\big) \qquad\text{and}\qquad \alpha^{\Gamma_\beta}(\beta M) \in KO_n\big(C^*_{r,\bR} (\Gamma_\beta)\big)
\end{gather*}
implies that $M_\Sigma$ admits an adapted psc~metric.
\end{Theorem}

\begin{Remark}
We notice that the condition $(i)$ holds when $L$ is a sphere, an odd
complex projective space, or when $L=G$.
\end{Remark}
In the last part of this paper, Section~\ref{sec:secondary}, we begin
to analyze the homotopy type of the space~$\mathcal R^+_w(M_\Sigma)$
of adapted metrics of {\psc} on $M_\Sigma$, in the case where this
space is non-empty. One of the key results is the following.

\begin{Theorem}
\label{thm:section}
Let $M_\Sigma$ be an $(L,G)$-fibered spin pseudomanifold, and assume
that $M_\Sigma$ admits an adapted psc~metric. Fix a connection
$\nabla^p$ on the associated principal $G$-bundle over $\beta M$. Let~$\resS\co \mathcal{R}_w^+(M_\Sigma)\to \mathcal{R}^+(\beta M)$ be the
forgetful map sending a psc metric $g$ on $M_\Sigma$, interpreted as a~pair $(g_M, g_{\beta M})$, to the metric $g_{\beta M}$ on $\beta M$.
Then $\resS$ is surjective onto $\mathcal{R}^+(\beta M)$, and it has a~{\lp}non-canonical{\rp} section. In particular,
there is a split injection of the homotopy groups
of~$\mathcal{R}^+(\beta M)$ into those of $\mathcal{R}_w^+(M_\Sigma)$.
\end{Theorem}

We use this result to detect non-trivial homotopy
groups of the space $\mathcal{R}_w^+(M_\Sigma)$. Namely, we let
$M_{\Sigma}= M\cup -N(\beta M)$, as before. Once we fix a base
point, a metric $g_0\in \cR^+_w(M_\Sigma)$, it~determines
corresponding metrics $g_{\beta M,0}\in \cR^+(\beta M)$ and
$g_{M,0}\in \cR^+(M)_{g_{\p M,0}}$ (where $g_{\p M,0}$ is given by
the metric $g_{\beta M,0}$ and $g_L$). Then we have
index-difference homomorphisms:
\begin{gather*}
 \mathrm{inddiff}_{g_{\beta M,0} }\colon\ \pi_q \big( \cR^+ (\beta M)\big)\to KO_{q+n-\ell}
 \end{gather*}
 and
\begin{gather*}
 \mathrm{inddiff}_{g_{\p M,0} }\colon\ \pi_q \big( \cR^+ (M)_{g_{\p M,0}}\big)\to KO_{q+n+1} ,
\end{gather*}
where $KO_k$ is the $k$-th coefficient group for real $K$-theory
(equal to $\bZ$ for $k$ a multiple of $4$ and~$\bZ/2$ for $k\equiv1,2\mod 8$);
see Section~\ref{sec:secondary}
and~\cite{BER-W,ER-W2}. According to Theorem~\ref{thm:section}, we can choose a splitting
\begin{gather*}
\pi_q \big(\cR^+_w(M_\Sigma)\big)=\pi_q \big(\cR^+(M)_{g_{\p M}}\big)\oplus \pi_q \big(\cR^+ (\beta M)\big).
\end{gather*}
(For $q=1$ one might get a semidirect product
instead of a direct sum,
as there is no obvious reason why the fundamental group should be
abelian.) In particular, we prove the
following result (see Corollary~\ref{cor:homot-groups2} for more details):
\begin{Theorem}
Let $M_\Sigma$ be a spin $(L,G)$-fibered compact pseudomanifold with
$L$ a simply connected homogeneous space of a compact semisimple Lie
group, and $n-\ell-1\geq 5$, where $\dim M=n$, $\dim L =\ell$.
Assume that $M_\Sigma$ admits an adapted psc~metric. Then the
composition
\begin{equation*}
  \xymatrix{\pi_q \big(\cR^+_w(M_\Sigma)\big)\ar[r]^(.3){\cong}\ar[rd]_{{\rm inddiff}_{g_0}}&
    \pi_q \big(\cR^+ (M)_{g_{\p M}} \big)\oplus \pi_q \big(\cR^+ (\beta M)\big)
    \ar[d]^{{\rm inddiff}_{g_{\p M,0} } \oplus
      {\rm inddiff}_{g_{\beta M,0}}}
    \\ & KO_{q+n+ 1}\oplus KO_{q+n- \ell}}
\end{equation*}
is surjective rationally and onto the $2$-torsion.
\end{Theorem}

We address the more general case when $M_{\Sigma}$ has
non-trivial fundamental group in Corolla\-ries~\ref{cor:homot-groups3} and~\ref{cor:stable-homot-groups}.

\section[KO-obstructions on L-fibered pseudomanifolds]
{$\boldsymbol{KO}$-obstructions on $\boldsymbol{L}$-fibered pseudomanifolds}\label{sec:obstructions-general}

\def\Di{\mathfrak{D}\kern-6.5pt/}
\def\Spi{\mathfrak{S}\kern-6.5pt/}
\def\tDi{\widetilde{\Di}}

\subsection[KOn-classes on L-fibered pseudomanifolds]
{$\boldsymbol{KO_n}$-classes on $\boldsymbol{L}$-fibered pseudomanifolds}

Let $(M_\Sigma,g)$ be a Thom--Mather space of depth one, endowed with
an adapted wedge metric $g$. We~denote as usual by $\beta M$ the
singular locus of $M_\Sigma$ and by $L$ the link. The resolved
manifold, a manifold with fibered boundary, will be denoted by
$M$.
\begin{Remark}
 We emphasize that for now the fibration $\phi\co \p M \to \beta M$ is
 assumed to be just a smooth fiber bundle with a fiber $L$, without
 any restriction on the structure group of that fibration. Then we
 say that $(M_\Sigma,g)$ is an $L$-fibered pseudomanifold, as opposed
 to an $(L,G)$-fibered pseudomanifold, which
 is the special case when the structure group $G$ of
 the fibration $\phi\co \p M \to \beta M$ acts on $L$ by isometries of
 a certain metric $g_L$. However, the relevant analytical constructions
 and results concerning the Dirac operators we need are well-studied
 and understood even for the general case of $L$-fibered pseudomanifolds.
\end{Remark}
Recall that $M_\Sigma=M\cup _{\p M} (-N(\beta M))$,
where $N(\beta M)$ is the tubular neighborhood of
the singular locus, which is also the total space of a
fiber bundle
\begin{gather}\label{eq1-b-new}
c(L) \to N(\beta M)\longrightarrow \beta M
\end{gather}
{\sloppy
with fiber equal to the cone over the link $L$. This
also determines a smooth fiber bundle $\phi\co \p M \to \beta M$ which is
a restriction of the bundle (\ref{eq1-b-new}) to the link $L$. We~denote by $T(\p M/\beta M)\to \beta M$ the corresponding vertical
tangent bundle. As in~\cite[Section 2.3]{BB-PP-JR}
(to which we refer for further
details) we fix a connection on the fiber bundle
$L\to \partial M\xrightarrow{\phi} \beta M$.
The metric $g$ on $M_{\Sigma}$ restricted to the regular part of
$N (\beta M)$ is assumed to have the following structure:
\begin{gather*}
g= {\rm d}r^2 + r^2 g_{\p M/\beta M} + \phi^* g_{\beta M} + O(r)
\end{gather*}}\noindent
with $r$ denoting the radial variable along the cone and where
the isomorphism between the horizontal bundle $H$ of the chosen connection
and $\phi^* T(\beta M)$ has been used.
(In this general setting we employ the notation $g_{X/B}$ for a metric
on the vertical tangent bundle of a smooth
fiber bundle $X\to B$.)

The resolved manifold $M$ inherits two metrics: the restriction of $g$
to $M$, a Riemannian metric of product type near the boundary, denoted
$g_M$, and the extension of the metric $g$ on~$\mathring{M}\equiv
M_\Sigma^{\rm reg}$ to the \emph{wedge metric} on the \emph{wedge tangent bundle}
${}^w TM\to M$. (This was defined in~\cite[Section~4]{almp-package} and
in~\cite{albin-gellredman-sigma}
under the alternate name \emph{incomplete edge tangent bundle};
see again Part I~\cite[Section~2.3]{BB-PP-JR}
for a quick introduction to the wedge tangent bundle). We~assume that~$M$, or equivalently
$M_\Sigma^{\rm reg}$, is given a spin structure. This fixes a spin
structure on $(\p M,g_{\p M})$ also. Then we assume that $\beta M$
is also spin and fix a spin structure for $(\beta M,g_{\beta M})$.
This also fixes a spin structure for the vertical
tangent bundle $T(\p M/\beta M) \to \beta M$ endowed with
the vertical metric $g_{\p M/\beta M}$.

Let us recall some basic facts in spin geometry. We~refer to~\cite[Chapter~II, Section~7]{lawson89:_spin} for further details.
Let $C\ell_n$ denote the Clifford algebra, and
$\ell\co {\rm Spin}_n\to {\rm Hom} (C\ell_n,C\ell_n)$
the representation given by
left multiplication. Then we denote by $\Spi_g (M)$ the bundle given
by~$P_{{\rm spin}}\times_{\ell} C\ell_n$ (here $P_{{\rm spin}}$
is the principal $\Spin(n)$-bundle defined by the spin structure).
There is a fiberwise action of
$C\ell_n$ on $\Spi_g (M)$ on the right which makes $\Spi_g (M)$ a
bundle of rank one $C\ell_n$-modules. The bundle
$\Spi_g (M)$ inherits a Levi-Civita connection $\nabla$.
Let $\Di_g$ be the associated
$C\ell_n$-linear Atiyah--Singer operator; thus, by definition,
$\Di_g = {\rm cl}\circ \nabla$.
The operator $\Di_g$ has the usual local expression
\begin{gather*}
\Di_g =\sum_j {\rm cl}\big(e^j\big)\nabla_{e_j}
\end{gather*}
with $\{e_j\}$ a local orthonormal frame of vector fields
and $\{e^j\}$ the dual basis defined by the metric.
$\Di_g$ is a $\bZ/2$-graded odd formally self-adjoint operator of
Dirac type commuting with the right action of $C\ell_n$.
For this operator the Schr\"odinger--Lichnerowicz formula holds:
\begin{gather*}
\Di^2_g=\nabla^*\nabla+ \frac{1}{4}\kappa_g,
\end{gather*}
with $\kappa_g$ denoting the scalar curvature of $g$.
See again~\cite[Chapter~II, Section~7]{lawson89:_spin} and also~\cite{Rosenberg-3} for more details on this crucial point.

\medskip
\noindent
{\bf Notation.} Unless absolutely necessary we shall omit the
reference to the wedge metric $g$ in the bundle
and the operator, thus denoting the
$C\ell_n$-linear Atiyah--Singer operator simply by
\begin{gather*}
 \Di\colon\ C^\infty\big(M,\Spi(M)\big)\to C^\infty\big(M,\Spi(M)\big).
\end{gather*}
Moreover, we
shall often use the shorter notation $\Spi$ instead of $\Spi (M)$.

\medskip
As in Part I, we can regard $\Di$ as a wedge differential operator of
first order on $L^2 (M,\Spi)$, initially with domain
equal to $C^\infty_c \big(M_\Sigma^{{\rm reg}} ,\Spi\big)\subset L^2
(M,\Spi)$ (more on this below). We~are looking for self-adjoint $C\ell_n$-linear
extensions of this differential operator in $L^2 (M,\Spi)$.

We are assuming that $\partial M$ is a fiber bundle of spin manifolds,
$L\to \p M\xrightarrow{\varphi} \beta M$,
 and we fix a connection on this bundle so as to have a well-defined
 notion of horizontal and vertical subspaces.
Consequently
\begin{gather*}
\mathfrak{S} (\p M)\simeq
\mathfrak{S}( \p M/\beta M)\hat{\otimes}\varphi^* \mathfrak{S}(\beta M),
\end{gather*}
where $\hat{\otimes}$ denotes the graded tensor product.
Extending work of Bismut and Cheeger on fib\-rations of spin manifolds, see
\cite[Section~4]{BC-adiabatic}, Albin and Gell-Redman
made a careful study of the Levi-Civita connection near
the singular stratum of a depth-one spin pseudomanifold~-- see Sections~2.2
and~3.1 in~\cite{albin-gellredman-sigma}; this study implies the
following structure of $\Di$ near the singular stratum:
\begin{gather}\label{AS-operator-near-stratum}
\Di= {\rm cl}({\rm d}r) \partial_r + {\rm cl}({\rm d}r) \frac{\ell}{2r} + \frac{1}{r}\Di_{\partial M/\beta M}\,\hat{\otimes}\, {\rm Id} +
 {\rm Id}\,\hat{\otimes}\,\widetilde{\Di}_{\beta M} + \mathfrak{B}
\end{gather}
with $\ell=\dim L$,
$\Di_{\p M/\beta M}$ the vertical family of Atiyah--Singer operators on
the fibration
\begin{gather*}
L\to \p M\xrightarrow{\varphi} \beta M,
\end{gather*}
$\widetilde{\Di}_{\beta M}$ an explicit horizontal operator, and
$\mathfrak{B}$ a bundle endomorphism
which is $O(r)$. Formula~\eqref{AS-operator-near-stratum} can be rewritten as
\begin{gather*}
\Di= r^{-1}\bigg( {\rm cl}({\rm d}r) r\partial_r + {\rm cl}({\rm d}r) \frac{\ell}{2} + \Di_{\partial M/\beta M}\hat{\otimes} {\rm Id} +
r {\rm Id}\hat{\otimes}\widetilde{\Di}_{\beta M} + r \mathfrak{B} \bigg)
\end{gather*}
and exhibits $\Di$ as a wedge differential operator of order 1: $\Di = r^{-1} \Di^e$,
with $\Di^e$ an edge differential operator. See Part~I
and of course~\cite{albin-gellredman-sigma}
for an introduction to edge and wedge operators.
With a small abuse of notation, widely used in family index theory,
we denote by $\Di_{L}$ the generic operator of the vertical family
$\Di_{\p M/\beta M}$ and by ${\rm spec}_{L^2}(\Di_{L})$ its spectrum.

The following result has been discussed in Part I:

\begin{Theorem}
Assume that
\begin{gather}\label{geometric-witt-part2}
 {\rm spec}_{L^2} (\Di_{L})\cap (-1/2,1/2)=\varnothing \qquad
 \text{for each fiber } L.
\end{gather} Then the following holds:
\begin{enumerate}\itemsep=0pt
\item[$(i)$] The operator $\Di$ with domain
 $C^\infty_c \big(M_\Sigma^{{\rm reg}},\Spi\big)\subset L^2 (M,\Spi)$ is
 essentially self-adjoint.
\item[$(ii)$]
 Its unique self-adjoint extension, still
 denoted by $\Di$, defines a $C\ell_n$-linear Fredholm operator and
 thus a class $\alpha_w (M_\Sigma,g)$ in $KO_n$, with $n=\dim M_\Sigma$.
\end{enumerate}
\end{Theorem}

As explained in Part I, $(i)$ and $(ii)$ are direct consequences of
the analysis
developed in~\cite{albin-gellredman-sigma} (which builds in turn on the slightly
more complicated case of the signature operator on Witt spaces, see~\cite{almp-package}).
The main step is the construction
of a $Cl_n$-linear parametrix for an operator~$\Di$ satisfying~\eqref{geometric-witt-part2}; this
is based crucially on the construction of a parametrix for
the edge operator~$\Di^e$ associated to $\Di$, using
Mazzeo's edge pseudodifferential calculus~\cite{MR1133743}.

The following result also follows directly from~\cite{albin-gellredman-sigma}
and from the Schr\"odinger--Lichnerowicz formula, which is valid for $\Di$:

\begin{Theorem}\label{geometric-witt-part2-addendum}
 If the tubular neighborhood of the singular stratum
 $\big(N(\beta M),g|_{N(\beta M)}\big)$ has non-negative
 scalar curvature, then~\eqref{geometric-witt-part2} holds. See
 \textup{\cite[Theorem~1.3]{albin-gellredman-sigma}}.

 If $\big(M_\Sigma^{{\rm reg}},g\big)$ has psc everywhere then~\eqref{geometric-witt-part2} holds {\lp}this is clear from above, given that
 $\big(N(\beta M)$, $g|_{N(\beta M)}\big)$ has psc{\rp} and the unique
 self-adjoint extension $\Di$ is $L^2$-invertible; in particular
\begin{gather*}
\alpha_w (M_\Sigma,g)=0\qquad\text{in}\quad KO_n.
\end{gather*}
\end{Theorem}

\begin{Definition}
We shall say that the stratified spin pseudomanifold $(M_\Sigma,g)$ is
{\it geometric-Witt} if the metric $g$ is
such that $ {\rm spec}_{L^2} (\Di_{L})\cap (-1/2,1/2)=\varnothing$
 for each fiber $L$.
\end{Definition}
This notion is taken from~\cite{albin-gellredman-sigma}, where it is
applied to any generalized wedge Dirac operator.

If $(M_\Sigma,g)$ does not satisfy \textup{\eqref{geometric-witt-part2}}
but is such that the vertical metric $g_{\p M/\beta M}$ is of psc along
the fibers, then we can still define
a wedge-alpha class $\alpha_w (M_\Sigma,g)\in KO_n$. Indeed,
by the Schr\"odinger--Lichnerowicz formula applied to the vertical family
we know that there exists $\epsilon > 0$ such that
\begin{gather}\label{eq:epsgap}
 {\rm spec}_{L^2} (\Di_{L})\cap (-\epsilon,\epsilon)
 =\varnothing \qquad \text{for each fiber } L.
\end{gather}
From~\eqref{eq:epsgap}
we can achieve condition~\eqref{geometric-witt-part2} above by
rescaling the vertical metric $g_{\p M/\beta M}$.
(This requires some adjustments in order to glue the
rescaled metric around the singular locus $\beta M$ to the original
metric on the complement of $N(\beta M)$, along the lines of~\cite[Appendix 6]{BW}.)
We~obtain in this way a wedge-alpha class
for $\Di$ by considering
 the wedge-alpha class of the operator associated to the rescaled metric.

There is a different (and in fact preferable) realization of this
class that makes use of a~natural self-adjoint domain defined for any
wedge Dirac operator with invertible vertical family along the fibers
of the boundary fibration. This is the so called
 vertical APS domain, see~\cite[Definition~2.3]{a-gr-dirac}. By
 applying the definition to our case we thus obtain a closed
 self-adjoint extension of $\Di$, deno\-ted~$\mathcal{D}_{{\rm VAPS}}
 (\Di)$. Following~\cite{a-gr-dirac} one proves that on this domain
 $\Di$ admits a parametrix, that is, an~inverse modulo compacts,
 which can be used in order to see that $\big(\Di, \mathcal{D}_{{\rm VAPS}} (\Di)\big)$ is a $C\ell_n$-linear Fredholm operator on its
 domain endowed with the graph norm. We~obtain in this way a~class
 $\alpha_w (M_\Sigma,g)\in KO_n$; one can prove that the class
 defined through the vertical APS domain and the class defined by the
 operator associated to the rescaled metric are in fact equal. See~\cite[Remark 4.10]{a-gr-dirac} for a sketch of the argument.

\begin{Definition}
We shall say that the stratified spin pseudomanifold $(M_\Sigma,g)$ is
{\it psc-Witt} if the metric $g$ is of psc along the links, i.e., if
the vertical metric $g_{\p M/\beta V}$ induces on each fiber $L$ a~metric of psc.
\end{Definition}

Given a {\it psc-Witt} stratified spin pseudomanifold $(M_\Sigma,g)$
we define its wedge alpha class in~$KO_n$ by considering the
$C\ell_n$-linear Fredholm operator $\big(\Di, \mathcal{D}_{{\rm VAPS}}(\Di)\big)$.

\begin{Remark}\label{(G,L)-remark}
In this article, which concentrates on $(L,G)$-pseudomanifolds, with
$L = G/H$ a simply connected homogeneous space and $G$ a compact
semisimple Lie group acting on $L$ by isometries,
we do not need to consider the vertical APS domain
 $\mathcal{D}_{{\rm VAPS}} (\Di)$ or, equivalently, the rescaled
 metric.
Indeed, the proof of~\cite[Theorem~1.3]{albin-gellredman-sigma} shows
that \textup{\eqref{geometric-witt-part2}} is automatic when
$\scal_{g_L}=\scal_{S^{\ell}}=\ell(\ell-1)$, which is the
normalization we have adopted, see~\eqref{normalization}.
Put it differently, for the purposes of this
 article we can and we shall exclusively treat the {\it
 geometric-Witt} case.
\end{Remark}

\begin{Theorem}\label{theo:geometric-witt-part2-bis}
Let $(M_\Sigma,g)$ be a 
geometric-Witt spin pseudomanifold. Then there is a
well-defined fundamental class $[\Di_g]\in KO_n (M_\Sigma)$.
\end{Theorem}

\begin{proof}
Following~\cite[Section 11]{debord-lescure-rochon},~\cite[Section 7.1]{PZ}
consider $\mathcal{A}$, the algebra of smooth functions on $M$ that are
constant along the fibers of $\partial M\xrightarrow{\varphi} \beta M$;
there is a dense inclusion $\mathcal{A}\hookrightarrow C(X_\Sigma)$.
Once a parametrix for $\Di$ is constructed, one obtains the existence of
the fundamental
class in $KO_n (M_\Sigma)$ by proceeding exactly as for the signature
operator on Witt spaces~--
see~\cite[Section 6.2]{almp-package}, but taking $\mathcal{A}$ as a
dense subalgebra of $C(M_\Sigma)$
instead of the subalgebra of Lipschitz functions employed in
\cite[Section 6.2]{almp-package}.
The details are very similar and thus we omit them.
\end{proof}

\begin{Remark}\label{push-forward}
As usual, we have $\alpha_w(M,g)=\pi_* [\Di_g]$,
with $\pi$ the mapping of the compact pseudomanifold $M_\Sigma$ to a point.
\end{Remark}

\begin{Remark}\label{dependance-metric}
We remark that the wedge $\alpha$-class and the fundamental class
$[\Di_g]$ are unchanged if $g(t)$, $t\in [0,1]$, is a 1-parameter
family of adapted wedge metrics that are geometric-Witt for any $t\in
[0,1]$. This is in fact a special case of Theorem~\ref{bordism-invariance} below.
\end{Remark}
\begin{Remark}
Theorem~\ref{theo:geometric-witt-part2-bis}, Remarks~\ref{push-forward}
and~\ref{dependance-metric} also hold when $(M_\Sigma,g)$ is
psc-Witt, with the condition in Remark~\ref{dependance-metric} being
that $g(t)$ is psc-Witt for each $t\in [0,1]$. As we shall not use
this case, we do not discuss the details.
\end{Remark}

\subsection[C*-algebras and KO classes associated to D]
{$\boldsymbol{C^*}$-algebras and $\boldsymbol{KO}$ classes associated to $\boldsymbol{\mathfrak{D}\kern-7.9pt/}$}

In this second part of our work we would like to bring in the
fundamental group of the pseudomanifold $M_\Sigma$. In fact, we will
begin, more generally, with Galois $\Gamma$-covering spaces since this
approach is more general and simplifies the techniques.

Let $\Gamma$ be a discrete group. A $\Gamma$-cover of $M_\Sigma$ will
be denoted by $M_\Sigma^\Gamma\xrightarrow{\pi} M_\Sigma$. We~shall
keep the notation $\widetilde M_\Sigma\to M_\Sigma$ for the universal
cover of $M_\Sigma$. We~notice that the
$\Gamma$-cover $M_\Sigma^\Gamma$ has a~natural structure of depth-1
pseudomanifold, with strata equal to the inverse image of the strata
of $M_\Sigma$ through the projection map $\pi$. See~\cite{PZ} and
\cite{Piazza-Vertman} for more on this background material. Hence,
the singular stratum $\beta M^\Gamma$ of
$M_\Sigma^\Gamma$\footnote{As we have denoted the singular locus of
$M_\Sigma$ by $\beta M$, we denote the singular locus of
$M_\Sigma^\Gamma$ by $\beta M^\Gamma$.} is a $\Gamma$-cover of $\beta M$,
\begin{gather*}
\Gamma\to \beta M^\Gamma \to \beta M
\end{gather*}
with a link diffeomorphic again to $L$; similarly, the regular part
$\big(M_\Sigma^{\Gamma}\big)^{\rm reg}$ of $M_\Sigma^\Gamma$ is a $\Gamma$-cover of~the
regular part $M_\Sigma^{{\rm reg}}$ of $M_\Sigma$.
Even if our Galois cover is the universal cover of $M_\Sigma$, then~$\beta M^\Gamma$ need not be, in general, the
universal covering $\widetilde{\beta M}$ of $\beta M$.
A similar remark applies to the two regular parts.

Finally, if $M^\Gamma$ is the resolution of $M_\Sigma^\Gamma$ then it
is easy to see that $M^\Gamma$ is in a~natural way a~Galois covering
of $M$. The interior of $M^\Gamma$ is clearly a Galois $\Gamma$-cover
of the interior of $M$. If
\begin{gather*}
L\to \p M\xrightarrow{p} \beta M
\end{gather*}
is the boundary fibration of $M$ and
$L\to \p M^\Gamma\xrightarrow{p_\Gamma} \beta M^\Gamma$ is the boundary
fibration of $M^\Gamma$ and if~$\pi$ denotes the quotient map induced
by the action of $\Gamma$, then there is a commutative diagram
\begin{gather*}
\xymatrix{\p M^{\Gamma}\ar[r]^(.5){\pi}\ar[d]^{p_\Gamma}&\p
 M\ar[d]^{p}\\ \beta M^\Gamma\ar[r]^(.5){\pi}& \beta M }
\end{gather*}
with $\p M^\Gamma\xrightarrow{\pi} \p M$ inducing a diffeomorphism
between the fiber of $p_\Gamma$ over $x\in \beta M^\Gamma$ and the
fiber of $p$ over $\pi (x)$. In the sequel we shall always endow
$M_\Sigma^\Gamma$, or rather its regular part, with the lift of an
adapted wedge metric on $M$; this metric extends as a
$\Gamma$-equivariant wedge metric on all of the resolution $M^\Gamma$;
we denote this $\Gamma$-invariant metric by $g_\Gamma$. Let us assume
once again that~$M$ and $\beta M$ are spin; then also $M^\Gamma$ and
$\beta M^\Gamma$ will be spin. Consequently, $M_\Sigma^{{\rm reg}}$
and $\big(M^\Gamma_\Sigma\big)^{{\rm reg}}$ are spin.

Let $\Di_\Gamma$ be the associated $\Gamma$-equivariant $Cl_n$-linear
Atiyah--Singer operator. Notice that since the link of
$M^\Gamma_\Sigma$ is again $L$ the analysis to be developed for
understanding the properties of $\Di_\Gamma$ is obtained by combining
the usual analysis of $\Gamma$-equivariant operators on
$\Gamma$-covers of smooth compact manifolds and the wedge analysis on
compact stratified spaces. This principle is explained in detail in
\cite{Piazza-Vertman}.

For the next theorem, however, we shall rather use the
Mishchenko-Fomenko operator asso\-ci\-a\-ted to $\Di$. Define the
Mishchenko bundle $\mathcal{V}:=M^\Gamma_\Sigma\times_\Gamma
C^*_{r,\bR}(\Gamma)$, a bundle of finitely gene\-ra\-ted projective
$C^*_{r,\bR}(\Gamma)$-modules of rank 1 {\it over the whole
 $M_\Sigma$.} By definition $\Di_{{\rm MF}}$ is equal to $\Di$
twisted by $\mathcal{V}$ restricted to the regular part of
$M_\Sigma$. We~shall also refer to this operator as~the
Atiyah--Singer--Mishchenko operator and if necessary we shall denote it
more precisely by~$\Di_{{\rm MF},g}$. We~shall also consider the
associated edge-operator, $\Di^e_{{\rm MF}}$, obtained by twisting
$\Di^e$ by~the Mishchenko bundle $\mathcal{V}$.

\begin{Theorem}\label{theo:geometric-witt-K}
Let $M$ and $\beta M$ be spin and let $g$ be a wedge adapted
metric. Consider a Galois $\Gamma$-cover $M_\Sigma^\Gamma$ of
$M_\Sigma$ and endow the regular part with the associated
$\Gamma$-equivariant metric $g_\Gamma$ as above. Let us assume that
$(M_\Sigma, g)$ is geometric-Witt. Consider the Atiyah--Singer--Mishchenko operator
$\Di_{{\rm MF},g}$. Then 
there is a unique
self-adjoint extension of $\Di_{{\rm MF},g}$, with domain denoted
$\mathcal{D} (\Di_{{\rm MF},g})$ such that the following hold:
\begin{enumerate}\itemsep=0pt
\item[{\rm(1)}] the pair $\big(\Di_{{\rm
 MF},g},\mathcal{D} (\Di_{{\rm MF},g})\big) $ defines an unbounded
 Kasparov $ (\mathbb{R},C^*_{r,\mathbb{R}} \Gamma \otimes
 C\ell_n)$-bimodule and thus an index class ${\rm Ind}_w \big(\Di_{{\rm
 MF},g},M_\Sigma^\Gamma\big)\in KKO_n \big(\mathbb{R},C^*_{r,\mathbb{R}}
 \Gamma\big)$;
 as usual we identify the group
 $KKO_n\big(\mathbb{R},  C^*_{r,\mathbb{R}} \Gamma\big)$
 with the isomorphic group $KO_n \big(C^*_{r,\mathbb{R}} \Gamma\big),$

\item[{\rm(2)}] if $g$ has psc everywhere then ${\rm Ind}_w \big(\Di_{{\rm
 MF},g},M_\Sigma^\Gamma\big)=0$ in $KO_n \big(C^*_{r,\mathbb{R}} \Gamma\big),$

\item[{\rm(3)}] if $\mathrm{ass}
 \co KO_n(B\Gamma)\to KO_n \big(C^*_{r,\mathbb{R}}
 \Gamma\big)$ is the assembly map and if $f\co M_\Sigma\to B\Gamma$ is the
 classifying map of the $\Gamma$-cover, then
\begin{gather*}
 {\mathrm{ass}} \big(f_* [\Di_g]\big)=
 {\rm Ind}_w \big(\Di_{{\rm MF},g},M_\Sigma^\Gamma\big)\qquad\text{in}\quad KO_n \big(C^*_{r,\mathbb{R}} \Gamma\big).
\end{gather*}
\end{enumerate}
\end{Theorem}

\begin{proof}[Sketch of the proof]
The proof is an easy adaptation of the corresponding result for the
signature operator on Witt spaces, see~\cite[Proposition 6.4, Theorem~6.6]{almp-package}, and for this reason we shall be brief. The
operator $\Di_{{\rm MF}}$ acts on the sections of $\Spi\otimes
\mathcal{V}$. Since the Mishchenko bundle on the tubular neighborhood
of the singular locus $\beta M$ is the pull-back of a bundle on $\beta
M$, we see that the analysis to be developed in order to understand
$\Di_{{\rm MF}}$ is no more difficult than the one already developed
for $\Di$. More precisely, following the proof of~\cite[Proposition~6.4]{almp-package}, we~see that $N_q (\Di^e_{{\rm MF}})$, the normal
operator of $\Di^e_{{\rm MF}}$ at $q\in \beta M$ is equal, up to
conjugation by~a~bundle isomorphism, to $N_q (\Di)\otimes {\rm
 Id}_{C^*_{r,\bR}\Gamma}$ and so its invertibility properties, that
are crucial in the analysis developed in
\cite{albin-gellredman-sigma}, are a consequence of those already
established for $N_q (\Di)$. The~proof now proceeds parallel to the
one given in~\cite[Section 6.3]{almp-package} for the signature
operator on Witt spaces. Notice that in (2) we use again the fact that
a metric which is psc everywhere is psc in
$\big(N(\beta M),g|_{N(\beta M)}\big)$ and thus geometric-Witt.
\end{proof}

\noindent
 {\bf Notation.} We shall briefly denote the index class
${\rm Ind}_w \big(\Di_{{\rm MF},g},M_\Sigma^\Gamma\big)$ as
\begin{gather*}
\alpha^{\Gamma}_w (M_\Sigma,g)\in KO_* \big(C^*_{r,\mathbb{R}} \Gamma\big).
\end{gather*}
In case we want to be very precise about the Galois cover
$M^\Gamma_\Sigma$ involved in the definition of this class we shall
also write
\begin{gather*}
\alpha^{\Gamma}_w (M_\Sigma,g,f)\in KO_* \big(C^*_{r,\mathbb{R}} \Gamma\big), \qquad \text{where}\quad
f\colon\ M_\Sigma\to B\Gamma\quad\text{and}\quad M^\Gamma_\Sigma=f^* E\Gamma.
\end{gather*}

\begin{Remark}\label{dependance-metric-Gamma}
Similarly to Remark~\ref{dependance-metric}, the class
$\alpha^{\Gamma}_w (M_\Sigma,g)\in KO_* \big(C^*_{r,\mathbb{R}} \Gamma\big)$
is unchanged if~$g(t)$, $t\in [0,1]$, is a 1-parameter family of
adapted wedge metrics that are geometric-Witt for any $t\in
[0,1]$. This is in fact a special case of Theorem~\ref{bordism-invariance} below.
\end{Remark}
\begin{Remark}
The above results can also be established for psc-Witt pseudomanifold,
either by using the vertical APS domain or by rescaling in the fiber
direction on the boundary. We~shall not need this more general case.
\end{Remark}

\begin{Remark}
The existence of a fundamental class $[\Di]\in K_* (M_\Sigma)$
and of a class $\alpha^\Gamma (M_\Sigma,g)$ can be established
more generally for spin pseudomanifolds of arbitrary
depth, assuming of~course suitable Witt conditions along the links. This is based heavily on the general edge pseudodifferential calculus
developed by Albin and Gell-Redman in~\cite{a-gr-dirac}.
Details will appear in~\cite{AG-RP}.
\end{Remark}

\subsection[Cylindrical KO-theory classes and a gluing formula]
{Cylindrical $\boldsymbol{KO}$-theory classes and a gluing formula}

We decompose
\begin{gather*}
M_\Sigma= M\cup_{\p M} (-N(\beta M))\qquad\text{and}\qquad
M_\Sigma^{{\rm reg}}=M\cup_{\p M} (-N(\beta M)^{{\rm reg}}).
\end{gather*}
Let $g$ be an adapted wedge metric on $M_\Sigma$. Denote by $g_M$ the
Riemannian metric equal to the restriction of $g$ to $M$ and by
$g_{N(\beta M)}$ the metric equal to the restriction of $g$ to
$N(\beta M)$, the collar neighborhood of the singular stratum. By
assumption,
\begin{gather*}
g_{N(\beta M)}={\rm d}r^2 + r^2 g_{\p M/\beta M} + p^* g_{\beta M} + O(r).
\end{gather*}
Recall that an adapted wedge metric $g$ is such that $g_M$ and
$g_{N(\beta M)}$ are of product type in a~collar neighborhood of $\p
M$. We~make the hypothesis that $g$ is geometric-Witt and that the
whole metric $g_{\p M}$ is of psc. Now attach an infinite cylinder to
$M$ along the boundary $\p M$ and extend the metric to be constant on
the cylinder; similarly, attach an infinite cylinder to $N(\beta M)$
and extend the metric. Because of the hypothesis that $g_{\p M}$ is
of psc we have well defined classes
\begin{gather*}
 {\rm Ind}_{{\rm cyl}} \big(\Di_{{\rm MF},g_M},M^\Gamma\big),\qquad
 {\rm Ind}_{{\rm cyl},w} \big(\Di_{{\rm MF},g_{N(\beta M)}},N(\beta M)^\Gamma\big)\qquad\text{in}\quad
 KO_* \big(C^*_{r,\mathbb{R}} \Gamma\big).
\end{gather*}
For the existence of the cylindrical class on $M^\Gamma$ we refer the reader to~\cite{Bunke-Callias,Le-Pi-psc,Wahl-APS} and refe\-ren\-ces therein. The existence of the index class
$ {\rm Ind}_{{\rm cyl},w} \big(\Di_{{\rm MF}},N(\beta M)^\Gamma\big)$ follows from the above theorem and these references.

\medskip\noindent
\textbf{Notation.}
 We shall briefly denote the above index classes as
\begin{gather*}
\alpha^{\Gamma}_{{\rm cyl}} (M,g_M)\qquad\text{and}\qquad
\alpha^{\Gamma}_{{\rm cyl},w} \big(N(\beta M),g_{N(\beta M)}\big)
\in KO_* \big(C^*_{r,\mathbb{R}} \Gamma\big).
 \end{gather*}

\begin{Proposition}\label{prop:gluing}
Under the same assumptions as above, namely that $g$ is geometric-Witt and that $g|_{\p M}$ is of psc, the following gluing formula holds:
\begin{gather*}
\alpha^{\Gamma}_w (M_\Sigma,g)=\alpha^{\Gamma}_{{\rm cyl}} (M,g_M)+
\alpha^{\Gamma}_{{\rm cyl},w} \big(N(\beta M),g_{N(\beta M)}\big)
\qquad\text{in}\quad KO_* \big(C^*_{r,\mathbb{R}} \Gamma\big).
\end{gather*}
Consequently, if $\big(N(\beta M),g_{N(\beta M)}\big)$ if of psc then
\begin{gather}\label{w=cyl}
\alpha^{\Gamma}_w \big(M_\Sigma^\Gamma,g\big)=\alpha^{\Gamma}_{{\rm cyl}}
\big(M^\Gamma,g_M\big)\qquad\text{in}\quad KO_* \big(C^*_{r,\mathbb{R}} \Gamma\big).
\end{gather}
\end{Proposition}

\begin{proof}
The existence of these classes has been discussed above. The gluing
formula has been discussed for the signature operator on Witt and
Cheeger spaces by Albin--Piazza~\cite{ap-surgery}, based on a~well
known technique due to Bunke~-- see~\cite{Bunke-Callias}; the same
arguments, with minor modifications, apply here.
\end{proof}

\subsection{Spin bordism of geometric-Witt pseudomanifolds}
\label{subsect:psc-witt-bordism}
It is well-known that in the smooth context the
Fredholm index of the spin-Dirac operator and the index class of the
spin-Dirac operator
twisted by the Mishchenko bundle define group homomorphisms
\begin{gather}\label{spin-bordism-complex-smooth}
{\rm ind}\co\ \Omega^{{\rm spin}}_n \to \bZ,\qquad
{\rm Ind}_\Gamma\co\ \Omega^{{\rm spin}}_n (B\Gamma)\to
K_* (C^*_{r} \Gamma),
\end{gather}
where for the numeric index we limit ourselves to $n=4k$.
More generally,
the $\alpha$ invariant of the Atiyah--Singer operator $\Di$
defines a group homomorphism
\begin{gather*}
\alpha\co\ \Omega^{{\rm spin}}_n \to KO_n
\end{gather*}
while the $\alpha^\Gamma$-invariant of the Atiyah--Singer--Mishchenko
operator $\Di_{{\rm MF}}$ defines a group homomorphism
\begin{gather*}
\alpha^\Gamma\co\ \Omega^{{\rm spin}}_n (B\Gamma)\to
KO_* \big(C^*_{r,\mathbb{R}} \Gamma\big).
\end{gather*}
The results for the spin-Dirac operator, see~\eqref{spin-bordism-complex-smooth},
is well known and treated in several references. See for example
the detailed treatment in~\cite{baum-douglas-taylor} and its extension
to the $C^*_r\Gamma$-index class in~\cite{Leichtnam-Piazza-GAFA}.
For the passage to the Atiyah--Singer operator, see~\cite[Chapter~IV, Section~4, equation~(4.3)]{lawson89:_spin} and~\cite{Rosenberg-3} or, for a purely analytic argument,
the appendix to this paper,
Section~\ref{sec:appendix}.

\begin{Definition}\label{bordism-of-psc-Witt}
Let $(M_\Sigma,g)$ and $(M_\Sigma ', g')$ be two spin-stratified
pseudomanifolds of dimension $n$ and let us assume that they are both
geometric-Witt. Let $(W_\Sigma, {\bar g})$ be a spin-stratified
pseudomanifold of dimension $(n+1)$ with \emph{collared boundary} and
adapted wedge metric, also of~product-type near the boundary, such
that $\partial W_\Sigma=M_\Sigma\cup (-M_\Sigma ')$ and ${\bar
 g}|_{\partial W_\Sigma}= g\cup g'$. We~shall say that
$(W_\Sigma,{\bar g})$ provides a geometric-Witt
bordism between $(M_\Sigma,g)$ and $(M_\Sigma ',g')$ if the adapted
wedge metric ${\bar g}$ is \emph{globally} psc-Witt.
\end{Definition}
If $\Gamma$ is a discrete group as above and $f\co M_\Sigma\to B\Gamma$
and $f'\co M_\Sigma '\to B\Gamma$ are classifying maps, then we shall
say that $(M_\Sigma,g,f)$ is geometric-Witt bordant to $(M_\Sigma ', g',f')$
if there exists $(W_\Sigma,{\bar g})$ as above and a continuous
map $F\co W_\Sigma\to B\Gamma$ restricting to $f$ and $f'$ on the
boundary. We~use the notation $(W_\Sigma,{\bar g},F) \co
(M_\Sigma,g,f)\bord (M_\Sigma ', g', f')$ for such a geometric-Witt bordism.
Notice that $f$, $f'$ and $F$ define in a natural way Galois coverings
$M_\Sigma^\Gamma$, $(M^\prime_\Sigma)^\Gamma$, $W_\Sigma^\Gamma$ and
thus Mishchenko bundles on $M_\Sigma$, $M_\Sigma '$ and $W_\Sigma$. We~denote briefly by $\Di_{{\rm MF}}$ and $\Di\,^\prime_{{\rm MF}}$ the
corresponding Atiyah--Singer--Mishchenko operators on $M_\Sigma$ and
$M_\Sigma^\prime$.

\begin{Theorem}\label{bordism-invariance}
If there exists a geometric-Witt bordism
$(W_\Sigma,\bar g)\co (M_\Sigma,g)\bord (M_\Sigma ', g')$, then
\begin{gather*}
\alpha_w (M_\Sigma,g)=\alpha_w (M_\Sigma ', g')\qquad\text{in}\quad KO_n.
\end{gather*}
If there exists a geometric-Witt bordism $(W_\Sigma,\bar g,F)\co
(M_\Sigma,g,f)\bord (M_\Sigma ', g', f')$, then
\begin{gather*}
\alpha^\Gamma_w ( M_\Sigma,g,f)=\alpha^\Gamma_w (M^\prime_\Sigma,g',f')
\qquad\text{in}\quad
KO_n \big(C^*_{r,\mathbb{R}} \Gamma\big),
\end{gather*}
where we recall that the above classes are defined in terms of
$M_\Sigma^\Gamma=f^* E\Gamma$ and $(M^\prime_\Sigma)^\Gamma=(f')^* E\Gamma$.
\end{Theorem}

\begin{proof}[Sketch of proof]
Once an analytic proof of the spin-bordism
invariance is given in the smooth closed case, the argument
for geometric-Witt pseudomanifolds is the same. For this reason,
we just give a sketch of the argument. Indeed, by general principles,
the proof is reduced to a collar neighborhood of the boundary where
the two directions,
the one of the boundary and the one normal to the boundary, are
completely decoupled.
Now, for smooth spin manifolds we have provided
a purely analytic proof of the spin-bordism invariance
of the $\alpha$-class and of the $\alpha^\Gamma$-class
in the appendix;\footnote{This is well known to the experts,
 but we could not pin down a quotable reference.}
it is very easy to extend
this proof to the geometric-Witt case, exactly as it is done for the
signature operator on Witt pseudomanifolds in order to establish the
Witt-bordism invariance
of the signature index and of the $C^*_r \Gamma$ signature-index class.
See~the proof of~\cite[Theorem~7.1]{almp-package}.
\end{proof}

\begin{Remark}
One can extend Definition~\ref{bordism-of-psc-Witt}
to the psc-Witt case and the analogue of Theo\-rem~\ref{bordism-invariance} holds.
\end{Remark}

\section[KO-obstructions on (L,G)-fibered pseudomanifolds]
{$\boldsymbol{KO}$-obstructions on $\boldsymbol{(L,G)}$-fibered pseudomanifolds}

In this section we shall finally meet the obstructions to the
existence of a wedge metric of psc on an $(L,G)$-fibered
pseudomanifold. We~treat first the general case of $L$-fibered
pseudomanifolds and then specialize to the case of $(L,G)$-fibered
pseudomanifolds; the latter are the pseudomanifolds for which we shall
prove an existence theorem.

\subsection[KO-obstructions for general L-fibered pseudomanifolds]
{$\boldsymbol{KO}$-obstructions for general $\boldsymbol{L}$-fibered pseudomanifolds}
Let $\Gamma$ be a finitely generated
discrete group, and assume that $M_\Sigma^\Gamma$ is a $\Gamma$-cover
of an $L$-fibered pseudomanifold $M_\Sigma=M\cup
_{\p M} (-N(\beta M))$. If $M_\Sigma$ admits an adapted
wedge metric of psc then we know that $M_\Sigma$ is geometric-Witt and
 $\alpha^\Gamma_w (M_\Sigma,g)=0$ in $KO_* (C^*_{r,\mathbb{R}}
 (\Gamma))$. See Theorems~\ref{geometric-witt-part2-addendum} and~\ref{theo:geometric-witt-K}(2). Since
$\big(N(\beta{M}),g_{N(\beta M)}\big)$ is of psc we also have that
$\alpha^{\Gamma}_w (M_\Sigma,g)=\alpha^{\Gamma}_{{\rm cyl}} (M,g_M)$
in~$KO_* \big(C^*_{r,\mathbb{R}} \Gamma\big)$. See Proposition~\ref{prop:gluing}, and in particular~\eqref{w=cyl}.

Summarizing, if $g$ is an adapted
 wedge metric of psc on an $L$-fibered pseudomanifold of~dimen\-sion
 $n$, then the following necessary condition is fulfilled:
\begin{gather*}
\alpha^{\Gamma}_w (M_\Sigma,g)=\alpha^{\Gamma}_{{\rm cyl}} (M,g_M)=0 \qquad\text{in}\quad
KO_n \big(C^*_{r,\mathbb{R}} (\Gamma)\big).
\end{gather*}
Consider now the singular locus of $M^\Gamma_\Sigma$, denoted $\beta
M^\Gamma$. We~already observed that the cove\-ring map
$M_\Sigma^\Gamma\to M_\Sigma$ induces a $\Gamma$-cover $\beta
M^\Gamma\to \beta M$. Consider the Atiyah--Singer operator on~$(\beta
M,g_{\beta M})$ and let us denote it by $\Di^\beta_{g_{\beta M}}$. We~also have an Atiyah--Singer--Mishchenko operator $\Di_{{\rm MF},g_{\beta
 M}}^\beta$, obtained by twisting $\Di^\beta_{g_{\beta M}}$ with
the Mishchenko bundle associated to the $\Gamma$-cover $\Gamma \to
\beta M^\Gamma \to \beta M$. We~can consider ${\rm Ind}
\big(\Di^\beta_{{\rm MF},g_{\beta M}}, \beta M^\Gamma\big)$ in $KO_{n-\ell-1}
\big(C^*_{r,\mathbb{R}} (\Gamma)\big)$ that we denote, as usual, by
$\alpha^\Gamma (\beta M,g_{\beta M})$. Since $\beta M$ is a closed
smooth spin manifold we can in fact adopt the notation $\alpha^\Gamma
(\beta M)\in KO_{n-\ell-1} \big(C^*_{r,\mathbb{R}} (\Gamma)\big)$, given that
this class does not depend on the particular choice of the metric
$g_{\beta M}$. See~\cite{Rosenberg-3}. Unless further assumptions
are made on the fibration $L\to \partial M\to \beta M$ we cannot infer
that also $ \alpha^\Gamma (\beta M)=0$ in $ KO_{n-\ell-1}
\big(C^*_{r,\mathbb{R}} (\Gamma)\big)$, where $\ell=\dim L$.

In the next section, on the other hand, we shall
specialize this discussion to the case of~$(L,G)$-fibered
pseudomanifolds and get consequently more
precise information.

\subsection[KO-classes on (L,G)-fibered pseudomanifolds]
{$\boldsymbol{KO}$-classes on $\boldsymbol{(L,G)}$-fibered pseudomanifolds}
We follow the notation of the previous section but we now assume that
$M_\Sigma$ is an $(L,G)$-fibered pseudomanifold endowed with an
adapted wedge metric $g$. We~know that
 $(M_\Sigma,g)$ is geometric-Witt~-- see Remark~\ref{(G,L)-remark}.
 By Theorem~\ref{theo:geometric-witt-K} there exists a well defined
 wedge-alpha class $\alpha_w^\Gamma (M_\Sigma,g)$ in $KO_n
 \big(C^*_{r,\mathbb{R}} (\Gamma)\big)$. Moreover, from Theorem~3.5 of Part 1
 (item (1)) we know that $g|_{\partial M}$ is of psc and so we also
 have a cylindrical class $\alpha_{\cyl}^\Gamma (M,g_M)\in KO_n
 \big(C^*_{r,\mathbb{R}} (\Gamma)\big)$. As~already remarked in the
 Introduction, the space of adapted wedge metrics on an
 $(L,G)$-fibered pseudomanifold is contractible. Thus
 $\alpha_w^\Gamma (M_\Sigma,g)\in KO_n \big(C^*_{r,\mathbb{R}} (\Gamma)\big)$
 does not depend on the choice of the adapted wedge metric $g$. See
 Remark~\ref{dependance-metric-Gamma}. Similarly, the cylindrical
 $\alpha$ class $\alpha_{\cyl}^\Gamma (M,g_M)$ does not depend on the
 choice of $g$; indeed, if $g$ and $g'$ are two adapted wedge metrics
 joined by a 1-parameter family of adapted wedge metric, then $g_M$
 and $g'_M$ are joined by a path of metrics that are uniformly of psc
 on the boundary $\partial M$. The result then follows from well-known
 properties of index classes on manifolds with cylindrical ends. See
 for example~\cite{MR720933} or~\cite{MR1979016}. Thus we can adopt
 the notation
 \begin{gather*}
 \alpha_w^\Gamma (M_\Sigma)\qquad \text{and}\qquad \alpha_{\cyl}^\Gamma (M)
 \end{gather*}
for the two classes.

Finally, we remark that \emph{under
 the assumption that $(\beta M,g_{\beta M})$ is of psc}, the equality
\begin{gather*}
\alpha^{\Gamma}_w (M_\Sigma,g)=\alpha^{\Gamma}_{{\rm cyl}} (M,g_M)
\end{gather*}
holds generally; indeed, in this case the
wedge metric~$g$ restricted to $N(\beta M)$ has
psc~-- see again Theorem~3.5 of Part I (item (2))~--
and thus the result follows from Proposition~\ref{prop:gluing}.

Assume now that our $(L,G)$-fibered pseudomanifolds is
endowed with an adapted wedge metric of psc. We~then certainly have that
\begin{gather*}
\alpha^{\Gamma}_w (M_\Sigma,g)=\alpha^{\Gamma}_{{\rm cyl}} (M,g_M)=0 \qquad\text{in}\quad
KO_n \big(C^*_{r,\mathbb{R}} (\Gamma)\big)
\end{gather*}
since this is true even in the general $L$-fibered case.
However, in this particular case, because of Theorem~3.5 of
Part I (item (3)), we have additionally that
\begin{gather*}
\alpha^\Gamma (\beta M)=0\in
KO_{n-\ell-1} \big(C^*_{r,\mathbb{R}} (\Gamma)\big).
\end{gather*}
This discussion proves the obstruction theorem, Theorem~\ref{thm:psobstruction} in the introduction.

Our task in the next
two sections will be to show that under suitable additional
assumptions these necessary conditions for the existence of an adapted
wedge metric of psc are also sufficient.

\medskip
Recall now that for these special pseudomanifolds we have defined a
bordism theory\linebreak $\Omega^{\spin,(L,G)\fb}_*(-)$; see Part I, Section 4. We~end this section by observing
that in this special case of $(L,G)$-fibered pseudomanifolds we can
frame the alpha classes of this section in the following
elegant way.
\begin{Proposition}\label{prop:bordism-map-index}
Let $\Omega^{\spin,(L,G)\fb}_*(-)$ be the singular spin bordism group.
Then we have well-defined homomorphisms:
\begin{gather}
\alpha^\Gamma_w\co\ \Omega^{\spin,(L,G)\fb}_* (B\Gamma) \to KO_* \big(C^*_{r,\mathbb{R}} (\Gamma)\big),\nonumber
\\
\alpha^\Gamma_{{\rm cyl}}\co\ \Omega^{\spin,(L,G)\fb}_* (B\Gamma) \to KO_* \big(C^*_{r,\mathbb{R}} \nonumber (\Gamma)\big),
\\\label{bordism-hom-3}
 \alpha^\Gamma_{\beta}\co\ \Omega^{\spin,(L,G)\fb}_* (B\Gamma)
 \to KO_{*-\ell-1} \big(C^*_{r,\mathbb{R}} (\Gamma)\big).
\end{gather}
\end{Proposition}

\begin{proof}
Let $[M_\Sigma,f\colon M_\Sigma\to B\Gamma]\in \Omega^{\spin,(L,G)\fb}_* (B\Gamma)$. We~know that $(L,G)$-fibered
pseudomanifolds are geometric-Witt. Fix an adapted wedge metric $g$ on
$M_\Sigma$ and set
\begin{gather*}
 \alpha^\Gamma_w \big([M_\Sigma,f\co M_\Sigma\to B\Gamma]\big):= \alpha^\Gamma_w
 (M_\Sigma,g,f)\in KO_* \big(C^*_{r,\mathbb{R}} (\Gamma)\big).
\end{gather*}
By definition,
 bordant elements in
the group $\Omega^{\spin,(L,G)\fb}_*(B\Gamma)$ are in fact geometric-Witt
bordant. Thus the first result follows immediately from Theorem~\ref{bordism-invariance}. Regarding the cylindrical class, exactly
the same arguments given for the $\alpha_{{\rm cyl}}$ invariant in
Part I, apply to $\alpha^\Gamma_{{\rm cyl}}$ , thus showing that
\begin{gather*}
\alpha^\Gamma_{{\rm cyl}}\co\ \Omega^{\spin,(L,G)\fb}_* (B\Gamma) \to
KO_* \big(C^*_{r,\mathbb{R}} (\Gamma)\big)
\end{gather*}
is well defined.

Regarding~\eqref{bordism-hom-3}: we have a homomorphism
$\beta_{B\Gamma}\co \Omega^{\spin,(L,G)\fb}_*(B\Gamma)\to \Omega_{*-\ell-1}(B\Gamma)$,
associating to $[M_\Sigma,f\co M_\Sigma\to B\Gamma]$
the class $[\beta M,f|_{\beta M}\co \beta M \to B\Gamma]$.
The homomorphism~\eqref{bordism-hom-3} is obtained by composing this
homomorphism $\beta_{B\Gamma}$ with the well-known alpha
homomorphism for closed spin manifolds.
Thus $\alpha^\Gamma_\beta [M_\Sigma,f\co M_\Sigma\to B\Gamma]$
is equal to the $C^*_{r,\bR} (\Gamma)$-index
of the Atiyah--Singer--Mishchenko operator
$\Di_{{\rm MF},g_{\beta M}}^\beta$ obtained by twisting
$\Di^\beta_{g_{\beta M}}$
with the Mishchenko bundle associated to the $\Gamma$-cover
$\Gamma \to \beta M^\Gamma \to \beta M$ defined by
$f|_{\beta M}\co \beta M \to B\Gamma$; it is therefore well defined.
\end{proof}

\subsection[Fundamental groups and KO-obstructions for (L,G)-fibered pseudomanifolds]
 {Fundamental groups and $\boldsymbol{KO}$-obstructions for $\boldsymbol{(L,G)}$-fibered \\pseudomanifolds}

 We assume now that $\Gamma=\pi_1
(M_\Sigma)$. We~also consider $\pi_1 (\beta M)$. In this section
we want to be particularly precise about the discrete groups involved
in the various Galois coverings.
\begin{itemize}\itemsep=0pt
\item \emph{In all of this section we assume that $M_\Sigma$ has
 an adapted wedge metric $g$ of psc.}
\end{itemize}
Under this assumption we observe that $(M,g_{M})$ is of psc and we
certainly have that
\begin{gather*}
\alpha_w^{\pi_1 (M_\Sigma)} (M_\Sigma,g)=0=
\alpha^{\pi_1 (M_\Sigma)}_{{\rm cyl}} (M,g_M)
\qquad\text{in}\quad KO_* \big(C^*_{r,\mathbb{R}} (\pi_1 (M_\Sigma))\big),
\end{gather*}
where the first class is relative to $\widetilde{M_\Sigma}$, the
universal cover of $M_\Sigma$, and the second class is relative to
$M^{\pi_1 (M_\Sigma)}$, the restriction of $\widetilde{M_\Sigma}$ to
the inverse image of
$M$. Also, since $\beta M$ admits a metric of psc, we also have, as
already remarked,
\begin{gather*}
 \alpha^{\pi_1 (M_\Sigma)}(\beta M)=0\qquad\text{in}\quad
 KO_{n-\ell-1} \big(C^*_{r,\mathbb{R}} (\pi_1 (M_\Sigma))\big),
\end{gather*}
where this class is relative to the $\pi_1 (M_\Sigma)$-cover of
$\beta M$ obtained by taking the restriction of $\widetilde{M_\Sigma}$ to
the inverse image of
$\beta M$ (this is in fact the singular locus of
$\widetilde{M_\Sigma}$). In addition we also have
\begin{gather*}
\alpha^{\pi_1 (\beta M)}(\beta M)=0\qquad\text{in}\quad
 KO_{n-\ell-1} \big(C^*_{r,\mathbb{R}} (\pi_1 (\beta M))\big),
\end{gather*}
where this class is relative to $\widetilde{\beta M}$,
the universal cover of $\beta M$,
$\pi_1 (\beta M)\to \widetilde{\beta M}\to \beta M$,
and it is defined in terms of
$\Di_{{\rm MF}, \pi_1 (\beta M)}$, the operator $\Di^\beta_{g_{\beta M}}$
twisted by the Mishchenko bundle
\begin{gather*}
\widetilde{\beta M}\times_{ \pi_1 (\beta M)} C^*_{r,\bR} (\pi_1 (\beta M)).
\end{gather*}
We remark that if $L$ is simply connected
(for example a sphere, complex or quaternionic projective space, or
complex Grassmannian), then the fibration $L\to \p
M\xrightarrow{\varphi} \beta M$ gives that $\varphi_*\co \pi_1(\p
M)\to \pi_1(\beta M)$ is an isomorphism; moreover, the tubular
neighborhood $N$ of $\beta M$ has a deformation retraction down to
$\beta M$. An easy application of Van Kampen's theorem then proves that
\begin{gather*}
\pi_1(M_\Sigma)\simeq \pi_1 (M).
\end{gather*}
Thus in this case there are exactly two fundamental groups to keep track of,
$\pi_1(M_\Sigma)\equiv \pi_1 (M)$ and $\pi_1(\beta M)$.

We summarize the index obstructions when $L$ is 1-connected and we
take into account the fundamental groups, in the following
 proposition.
\begin{Proposition}
Let $M_\Sigma$ be a $(L,G)$-fibered pseudomanifold.
Assume that $L$ is 1-connected so that $\pi_1(M_\Sigma)\equiv \pi_1 (M)$.
Let $g$ be an adapted wedge metric of psc on $M_\Sigma$. We~endow $M$ with the metric $g_M:=g|_{M}$.
Under the above hypothesis we have the following vanishing results:
\begin{gather*}
 \alpha_w^{\pi_1 (M_\Sigma)}(M_\Sigma,g)=
\alpha^{\pi_1 (M_\Sigma)}_{{\rm cyl}} ( M,g_M)=0
\qquad\text{in}\quad KO_n\big(C^*_{\bR,r}(\pi_1(M_\Sigma))\big),
\\
\alpha^{\pi_1 (M_\Sigma)} (\beta M)=0 \qquad\text{in}
\quad KO_{n-\ell-1} \big(C^*_{\mathbb{R},r} (\pi_1 (M_\Sigma))\big),
\\
\alpha^{\pi_1 (\beta M)}(\beta M)=0 \qquad\text{in}\quad KO_{n-\ell-1}\big(C^*_{\bR,r}(\pi_1(\beta M))\big),
\end{gather*}
where the first two classes are relative to
$\widetilde{M_\Sigma}$, the universal cover of $M_\Sigma$, and
to its restriction to the inverse image of
$M$ respectively; the third class is relative
to the restriction of $\widetilde{M_\Sigma}$ to
the inverse image of $\beta M$,
and the fourth class is defined in terms of
$\widetilde{\beta M}$, the universal cover of~$\beta M$.
\end{Proposition}

\section{Relevant surgery and bordism theorems}
\label{sec:main}
\subsection{The case of a single fundamental group}
We need two simple generalizations of our results from Part I
\cite[Theorems 4.6 and 4.7]{BB-PP-JR}. Given a simply connected link
manifold $L$ and a simply connected Lie group $G$, acting on~$L$, we
consider $(L,G)$-singular spin pseudomanifolds
$M_{\Sigma}$, where $M_{\Sigma}=M\cup_{\p M} N(\beta
 M)$, and the manifolds $M$ and $\beta M$ have non-trivial
fundamental groups. The first case to analyze is when the
pseudomanifold $M_{\Sigma}$ and its strata $M$ and $\beta M$ have the
same fundamental group $\Gamma$. While this case is also covered by
the more complicated case to be treated below in~Section~\ref{sec:twogroups},
we~have separated this case out because the argument is simpler and
clearer. In particular, such an~$M_{\Sigma}$ determines an element in
the bordism group $\Omega^{\spin,(L,G)\fb}_{*}(B\Gamma)$. In more
detail, this means that there are maps $\xi\co M_{\Sigma}\to B\Gamma$
and $\xi_{\beta}\co\beta M\to B\Gamma$ such that the following diagram
of isomorphisms is commutative:
\begin{gather}\label{eq-new1}
\begin{diagram}
 \node{\pi_1(M_{\Sigma})}
 \arrow{e,t}{\xi_*}
 \node{\Gamma}
 \\
 \node{\pi_1(M)}
 \arrow{n,t}{i_*}
 \arrow{ne,b}{(\xi|_M)_*}
 \node{\pi_1(\p M)}
 \arrow{w,t}{j_*}
 \arrow{e,t}{p_*}
 \node{\pi_1(\beta M).}
 \arrow{nw,t}{(\xi_{\beta})_*}
\end{diagram}
\end{gather}
Here $i\co M\hookrightarrow M_{\Sigma}$, $j\co \p M \hookrightarrow M$ are
canonical embeddings.

\begin{Theorem}[a refined surgery theorem]
\label{thm:surg1a}
Assume $G$ is a simply connected Lie group, $L$ is simply
connected and spin, and $\Gamma$ is a discrete
group. Let the element $x\in \Omega^{\spin,(L,G)\fb}_{*}(B\Gamma)$ be~represented by two $(L,G)$-singular spin pseudomanifolds $\xi\co
M_{\Sigma}\to B\Gamma$ and $\xi'\co M_{\Sigma}'\to B\Gamma$ of~dimension $n\ge 6 + \ell$, where $M_{\Sigma}=M\cup_{\p M} -N(\beta M)$,
$M'_{\Sigma}=M'\cup_{\p M} -N(\beta M')$ are given together with maps
$\xi_{\beta}\co \beta M\to B\Gamma$ and $\xi_{\beta}'\co \beta M'\to
B\Gamma$ and with structure maps $f \co\beta M \to BG$ and $f'
\co\beta M \to BG$. Assume the pseudomanifold $M_{\Sigma}$ and its
strata $M$ and $\beta M$ have the same fundamental group $\Gamma$,
i.e., \textup{(\ref{eq-new1})} is a commutative diagram of
isomorphisms.

Then there exists an $(L,G)$-bordism $W_{\Sigma}\co
M_{\Sigma}\rightsquigarrow M_{\Sigma}'$, $W_{\Sigma}= W\cup_{\p W}
N(\beta W)$ over $B\Gamma$, with a structure map
$\bar f \co \beta W\to BG$ restricting to the structure
maps on $\beta M$ and $\beta M'$, such that the pairs
$(\beta W, \beta M)$ and $(W, M)$ are $2$-connected.
\end{Theorem}
\begin{Remark}
We emphasize that the maps $\xi'\co M_{\Sigma}'\to
B\Gamma$ and $\xi_{\beta}'\co \beta M'\to B\Gamma$ are not assumed
to be isomorphisms on fundamental groups; the only requirement here
is that $(M_{\Sigma}',\xi',\xi_{\beta}')$ represents the same
bordism class $x\in \Omega^{\spin,(L,G)\fb}_{*}(B\Gamma)$ as
$(M_{\Sigma},\xi,\xi_{\beta})$.
\end{Remark}
\begin{proof}
We start with some $(L,G)$-bordism $W_{\Sigma}\co
M_{\Sigma}\rightsquigarrow M_{\Sigma}'$ over $B\Gamma$. In particular,
we have maps $\bar \xi \co W_{\Sigma}\to B\Gamma$ and $\bar \xi_{\beta}
\co \beta W \to B\Gamma$ such that $\bar \xi|_{M_{\Sigma}} = \xi$ and
$\bar \xi_{\beta}|_{\beta M_{\Sigma}} =
\xi_{\beta}$.
\begin{Remark}
For future use, we recall that
$W_{\Sigma}=W\cup_{\p^{(1)} W} N(\beta W)$, where $W$ is a manifold
with corners: $\p W = \p^{(0)} W \cup \p^{(1)} W$, where
\begin{gather*}
 \p^{(0)} W = M \sqcup - M', \qquad
 \p\big(\p^{(0)}W\big) = - \p \big(\p^{(1)} W\big) = \p M \sqcup - \p M',
\end{gather*}
and the $(L,G)$-fiber bundle $\p^{(1)} W$ is given by the map $\bar f
\co \beta W \to BG$, where $\bar f|_{\beta M}=f$ and~$\bar f|_{\beta
 M'}=f'$. Then we denote by
 \begin{gather*}
 \delta W_{\Sigma} = \p^{(0)} W
\cup_{\p(\p^{(0)}W)} N(\p (\beta W)) = M_{\Sigma}\sqcup M_{\Sigma}'.
\end{gather*}
\end{Remark}
By assumption, the maps $\xi$ and $\xi_{\beta}$ induce isomorphisms of
fundamental groups. As already remarked, the restrictions
$\bar\xi|_{M_{\Sigma}'} = \xi'$ and
$\bar \xi_{\beta}|_{\beta M_{\Sigma}'} = \xi_{\beta}'$
may not be isomorphisms on $\pi_1$. However,
the commutativity of the diagrams
\begin{gather*}
\begin{diagram}
\setlength{\dgARROWLENGTH}{1.95em}
 \node{\pi_1(W)}
 \arrow{e,t}{\bar \xi_*}
 \node{\Gamma}\\
 \node{\pi_1(M)}
 \arrow{n}
 \arrow{ne,r}{\xi_*}
\end{diagram}
\qquad\text{and}\qquad
\begin{diagram}
\setlength{\dgARROWLENGTH}{1.95em}
 \node{\pi_1(\beta W)}
 \arrow{e,t}{(\bar \xi_{\beta})_*}
 \node{\Gamma}\\
 \node{\pi_1(\beta M)}
 \arrow{n}
 \arrow{ne,r}{(\xi_{\beta})_*}
\end{diagram}
\end{gather*}
{\sloppy (with the vertical arrows induced by the inclusions) implies that the
homomorphisms $\bar \xi_* $: $\pi_1(W)\to \Gamma$ and $(\bar
\xi_{\beta})_*\co \pi_1(\beta W)\to \Gamma$ are surjective.

}

Then the first step is to do surgery on $\beta W$ to make the pair
$(\beta W, \beta M)$ $2$-connected. We~begin by killing the kernel of
$\pi_1(\beta W)\twoheadrightarrow \Gamma$. This proceeds as in the
proof of~\cite[Theorem~4.7]{BB-PP-JR}, the point being that since this
kernel maps trivially to $\pi_1(B\Gamma)$, the surgery can be done
over $B\Gamma$ without any extra effort. After taking care of the
fundamental group, the situation is the same as in~\cite[Theorem~4.7]{BB-PP-JR}. The process for making the pair $(W, M)$
$2$-connected is completely analogous.
\end{proof}
Once we have established Theorem~\ref{thm:surg1a}, we get a
corresponding bordism theorem; see Theorem~4.6 from Part I.
\begin{Theorem}[a refined bordism theorem]
\label{thm:bordism1}
Assume $G$ is a simply connected Lie group, $L$ is simply connected
and spin, and $\Gamma$ is a discrete group. Let the
element $x\in \Omega^{\spin,(L,G)\fb}_{*}(B\Gamma)$ be represented by
two $(L,G)$-singular spin pseudomanifolds $\xi\co M_{\Sigma}\to B\Gamma$
and $\xi'\co M_{\Sigma}'\to B\Gamma$ of dimension $n\ge 6 + \ell$,
where $M_{\Sigma}=M\cup_{\p M} -N(\beta M)$, $M'_{\Sigma}=M'\cup_{\p M}
-N(\beta M')$ are given together with maps $\xi_{\beta}\co \beta M\to
B\Gamma$ and $\xi_{\beta}'\co \beta M'\to B\Gamma$ and with structure
maps $f \co\beta M \to BG$ and $f' \co\beta M \to BG$. Assume the
pseudomanifold $M_{\Sigma}$ and its strata $M$ and $\beta M$ have the
same fundamental group $\Gamma$, i.e., {\rm (\ref{eq-new1})} is a
commutative diagram of isomorphisms. Furthermore, assume $M'_{\Sigma}$
has an adapted psc~metric $g'$.

Then there exists an $(L,G)$-bordism
$W_{\Sigma}\co M_{\Sigma}\rightsquigarrow M_{\Sigma}'$
over $B\Gamma$ together with
an adapted psc~metric $\bar g$ which is a product metric near the
boundary $\delta W_{\Sigma}= M_{\Sigma}\sqcup - M_{\Sigma}'$ such that
$\bar g|_{M_{\Sigma}'}=g'$. In particular, $M_{\Sigma}$ admits an
adapted psc~metric $g$.
\end{Theorem}
\subsection{The case of two fundamental groups} \label{sec:twogroups}
Next, we would like to address the natural situation when the
fundamental groups of $\beta M$ and~$M$ are different. For that, we
need to redo some of the bordism theory of~\cite[Section~4]{BB-PP-JR} and
fix \emph{two} discrete groups $\Gamma_{\beta}$ and $\Gamma$ and a
homomorphism $\theta\co \Gamma_{\beta}\to \Gamma$. This data will be
held fixed, and we consider pseudomanifolds
$M_\Sigma=M\cup_{\p M} -N(\beta M)$
with $(L,G)$-fibered singularities\footnote{Where as before $L$
 and $G$ are simply connected and $G$ acts transitively on $L$.} but
now $\Gamma_{\beta}$ is the fundamental group of $\beta M$ and $\p M$
and $\Gamma$ is the fundamental group of $M$, and $\theta \co
\Gamma_{\beta}\to \Gamma $ is the map of groups induced by the
inclusion. To describe a bordism group
$\Omega^{\spin,(L,G)\fb}_{n}\big(B\Gamma_{\beta}\xrightarrow{\theta_*}B\Gamma\big)$,
we need maps $\xi\co M\to B\Gamma$ and $\xi_{\beta}\co \beta M\to
B\Gamma_{\beta}$ such that the following diagram commutes:
\begin{gather*}
\begin{diagram}
 \node{\Gamma}
 \node{}
 \node{\Gamma_{\beta}}
 \arrow[2]{w,t}{\theta}
 \\
 \node{\pi_1(M)}
 \arrow{n,l}{\xi_*}
 \node{\pi_1(\p M)}
 \arrow{ne,t}{(\xi_{\p})_*}
 \arrow{w,t}{j_*}
 \arrow{e,t}{p_*}
 \node{\pi_1(\beta M).}
 \arrow{n,r}{(\xi_{\beta})_*}
\end{diagram}
\end{gather*}
(Here $\xi_{\p}=\xi_{\beta}\circ p$.)
In the case of interest, since $L$ is simply connected,
$p_*\co \pi_1 (\p M)\to \pi_1(\beta M)$,
$\xi_*\co \pi_1 M\to \Gamma$, and
$(\xi_{\beta})_*\co\pi_1 (\beta M )\to \Gamma_{\beta}$ are all isomorphisms.
Notice that in our case, as~already remarked,
by Van Kampen's theorem, the fundamental group of $M_{\Sigma}$ is
$\pi_1 (M)*_{\pi_1 (\beta M)} \pi_1 (\beta M) = \pi_1 (M) = \Gamma$.

\begin{Definition}
\label{def:relbordism-2}
Fix a group homomorphism $\theta\co\Gamma_\beta\to \Gamma$. We~define
$\Omega_n^{\spin,
 (L,G)\fb}\big(B\Gamma_\beta\xrightarrow{\theta_*}B\Gamma\big)$ to be the
bordism group of $n$-dimensional closed $(L,G)$-singular spin
pseudomanifolds $M_{\Sigma}$, where $M_{\Sigma}
 =M\cup_{\p M} - N(\beta M)$ is equipped with maps $\xi\co M\to
B\Gamma$ and
$\xi_{\beta}\co\beta M\to B\Gamma_{\beta}$ such that the diagram
\begin{gather*}
\begin{diagram}
 \node{B\Gamma}
 \node{}
 \node{B\Gamma_\beta}
 \arrow[2]{w,t}{B\theta}
 \\
 \node{M}
 \arrow{n,l}{\xi}
 \node{\p M}
 \arrow{ne,t}{\xi_{\p}}
 \arrow{w,t}{j}
 \arrow{e,t}{p}
 \node{\beta M}
 \arrow{n,r}{\xi_{\beta}}
\end{diagram}
\end{gather*}
commutes up to homotopy (here $\xi_{\partial}=\xi_\beta\circ p$).
A bordism
\begin{gather*}
 \big(W_{\Sigma},\bar \xi, \bar \xi_{\beta}\big)\co\
 (M_{\Sigma},\xi,\xi_{\beta}) \rightsquigarrow
 (M_{\Sigma}',\xi',\xi_{\beta}'),
\end{gather*}
between such objects is an $(L,G)$-singular spin pseudomanifold
$W_{\Sigma}= W\cup_{\p W} N(\beta W)$ with the boundary $\delta
W_{\Sigma}=M_{\Sigma}\sqcup -M_{\Sigma}'$ given together with maps
$\bar \xi\co W\to B\Gamma$ and $\bar \xi_{\beta }\co \beta W\to
B\Gamma_{\beta}$ such that
\begin{gather*}
 \bar \xi|_M=\xi, \qquad
 \bar \xi|_{M'}=\xi';\qquad
\bar \xi_{\beta }|_{\beta M}= \xi_{\beta},\qquad
\bar \xi_{\beta }|_{\beta M'}= \xi_{\beta}'.
\end{gather*}
\end{Definition}

\begin{Remark}
Note that, in general, $W$ here is a manifold with corners with
$\partial W = \partial^{(0)} W\cup \partial^{(1)} W$, where
$\partial^{(0)} W\cap \partial^{(1)} W = \partial M \sqcup -\partial
M'$ and $\partial^{(1)} W$ has an $(L,G)$-fibration structure over~$\beta W$, so that the map $\bar \xi_{\p }$ is just the bundle
projection followed by $\bar \xi_{\beta }$.
\end{Remark}
\begin{Remark}
Together with structure maps $f\co \beta M \to BG$,
$f'\co \beta M' \to BG$ and $\bar f \co\beta W \to BG$, we obtain a
spin bordism
$\big(\bar f \times \bar \xi_{\beta}\co\beta W \to BG\times B\Gamma_{\beta}\big)$
between
$(f \times \xi_{\beta}\co\beta M \to BG\times B\Gamma_{\beta})$ and
$(f' \times \xi_{\beta}'\co\beta M' \to BG\times B\Gamma_{\beta})$.
\end{Remark}
\begin{Remark}
Let $M_{\Sigma}=M\cup_{\p M} - N(\beta M)$ be as above,
i.e., $\pi_1 M= \Gamma$, $\pi_1\p M = \pi_1\beta M=\Gamma_{\beta}$,
and the homomorphism $\theta \co \Gamma_{\beta}\to \Gamma$ coincides
with the one given by the inclusion $\p M \xhookrightarrow{} M$. Then
there is a canonical map $u\co M_\Sigma\to B\Gamma$ such that
$u|_{M}=\xi$, and $u|_{N(\beta M)}$ is given by the composition
$N(\beta M)\xrightarrow{\mathrm{proj}}\beta
M\xrightarrow{\xi_{\beta}} B\Gamma_{\beta} \xrightarrow{\theta} B\Gamma$.
This construction defines a natural forgetful map
\begin{gather*}
\Omega^{\spin,(L,G)\fb}_{n}\big(B\Gamma_{\beta}\xrightarrow{B\theta}B\Gamma\big)
\to \Omega^{\spin,(L,G)\fb}_{n}(B\Gamma),
\end{gather*}
which sends $(M_{\Sigma},\xi,\xi_{\beta})$ to $(M_{\Sigma},u)$ as above.
\end{Remark}
To generalize the arguments in~\cite{BB-PP-JR} for proving existence of
well-adapted {\psc} metrics in some cases, we now need an
exact sequence into which the bordism group of Definition
\ref{def:relbordism-2} will fit. In general this is quite
complicated, so we only do an easy case as follows.

\begin{Theorem}
\label{prop:bordismtriangle}
Let $\Omega_n^{\spin,
 (L,G)\fb}\big(B\Gamma_\beta\xrightarrow{B\theta}B\Gamma\big)$ be a bordism
group as in Definition~$\ref{def:relbordism-2}$.
Assume that $L$ is a $G$-equivariant spin boundary, i.e., that there
is a spin $G$-manifold $\overline L$ with $\partial{\overline L} = L$
{\lp}as $G$-manifolds{\rp}. Then there is an exact sequence of bordism groups
\begin{gather}
\label{eq:bordismtriangle}
0 \to \Omega^{\spin}_{n}(B\Gamma) \xrightarrow{i}
\Omega_n^{\spin, (L,G)\fb}\big(B\Gamma_\beta\xrightarrow{B\theta}B\Gamma\big)
\xrightarrow{\beta} \Omega^{\spin}_{n-\ell-1}(B\Gamma_\beta\times BG)
\to 0.
\end{gather}
\end{Theorem}
\begin{proof}
First of all, note that the assumption on $L$ gives the surjectivity
of $\beta$. For if the map $\beta M\to B\Gamma_\beta\times BG$ represents an
element of $\Omega^{\spin}_{n-\ell-1}(B\Gamma_\beta\times BG)$ and we
form the associated $(L,G)$-bundle
\begin{gather*}
\xymatrix{L \ar[r] & X\ar[r]\ar[d]^p &B\Gamma_\beta\\
 & \beta M \ar[r]&B\Gamma_\beta \times BG,}
\end{gather*}
then by replacing $L$ by $\overline L$, we get a spin manifold
$\overline X$ of dimension $n$ with boundary $X$, fibering over
$\beta M$ with fiber $\overline L$ and fiber
bundle projection $\bar p$. Now we have a commutative
diagram
\begin{gather*}
\xymatrix{X \ar[r]^p \ar[d] &\beta M\ar[d]^= \ar[r]
 & B\Gamma_\beta\ar[d]^=\ar[dr]^{B\theta}\\
 \overline X \ar[r]^{\bar p} & \beta M\ar[r]
 &\ar[r] B\Gamma_\beta\ar[r]^{B\theta} & B\Gamma.}
\end{gather*}
The result is a
spin manifold $\overline X$ with boundary $X$, representing a class
in
\[
\Omega_n^{\spin, (L,G)\fb}\big(B\Gamma_\beta\xrightarrow{\theta}B\Gamma\big),
\]
mapping under $\beta$ to the class of $\beta M\to B\Gamma_\beta\times BG$.

The map $i\co \Omega^{\spin}_{n}(B\Gamma) \to
\Omega_n^{\spin, (L,G)\fb}\big(B\Gamma_\beta\xrightarrow{\theta}B\Gamma\big)$ just
comes from thinking of a closed manifold as a manifold with empty
fibered singularities. We~need to show this map is injective.
Suppose $M\xrightarrow{\xi} B\Gamma$
represents a class in $\Omega^{\spin}_{n}(B\Gamma)$. Suppose
it bounds as a manifold with $(L,G)$-fibered singularities
(in the bordism group in the center). Since $\p M$ and
$\beta M$ are empty, that means we
have an $(L,G)$-singular spin pseudomanifold with
boundary $W_{\Sigma}= W\cup_{\p W} N(\beta W)$ with the boundary
$\p W_{\Sigma}=M$ and maps $\bar \xi\co
W\to B\Gamma$, $\bar \xi_{\p }\co \p^{(1)} W\to B\Gamma_{\beta}$ and
$\bar \xi_{\beta }\co \beta W\to B\Gamma_{\beta}$ such that
$\bar \xi|_M=\xi$. Also, in this case, there are no corners,
i.e., $X= \partial^{(1)}W$ is disjoint from $\partial^{(0)}W = M$.
Now we want to modify $W$ to convert it to a manifold with
boundary exactly equal to $M$. Once again, we use the assumption that $L$
is a spin $G$-boundary to build an $\big(\overline L, G\big)$-fiber
bundle $\overline X$ over $\beta W$ bounding $X$, and
$W\cup_{\partial^{(1)}W} \overline X$ is our required null-bordism of $M$.
Thus the class of $M\xrightarrow{\xi} B\Gamma$
was already zero in $\Omega^{\spin}_{n}(B\Gamma)$.

Finally, we need to prove exactness in the middle. Suppose given
$M_\Sigma$ defining a class in
\begin{gather*}
\Omega_n^{\spin, (L,G)\fb}\big(B\Gamma_\beta\xrightarrow{B\theta}B\Gamma\big)
\end{gather*}
with $[\beta M\to B\Gamma_\beta\times BG]$ trivial in
$\Omega^{\spin}_{n-\ell-1}(B\Gamma_\beta\times BG)$. Choose a spin
manifold $Y$ bounding~$\beta M$, with a map to $B\Gamma_\beta\times BG$
extending the original map on $\beta M$. Form the associated
$(L,G)$ fiber bundle $X$ over $Y$. Now by gluing
$-X$ to $M$ along the common boundary $\p M$ we get a~closed
spin $n$-manifold $M_1$ with a map to $B\Gamma$ coming from $M$.
Indeed, this map extends over~$-X$ via the composite
$X\to B\Gamma_\beta \xrightarrow{B\theta}B\Gamma$.
Then $M_1$ represents the same class as
$M_\Sigma$, because
\begin{gather*}
 W = \big(M\times [0, 1]\big)\cup_{\p M\times \{1\}}X
 \end{gather*}
is a bordism between them in the sense of Definition
\ref{def:relbordism-2}.
\end{proof}

\begin{Theorem}[surgery theorem, general version]
\label{thm:surg2}
Assume $G$ is a simply connected Lie group,
$L^{\ell}$ is simply connected and spin, $\theta\co \Gamma_{\beta}\to
\Gamma$ is a group homomorphism, and $n\ge 6 + \ell$. Let $x\in
\Omega^{\spin,(L,G)\fb}_{n}\big(B\Gamma_{\beta}\xrightarrow{\theta_*}B\Gamma\big)$
be represented by two $(L,G)$-singular spin pseudomanifolds
$(M_{\Sigma},\xi,\xi_{\beta})$, $M_{\Sigma}=M\cup_{\p M}
-N(\beta M)$, and $(M_{\Sigma}',\xi',\xi_{\beta}')$,
$M'_{\Sigma}=M'\cup_{\p M} -N(\beta M')$ given together with maps
\begin{gather*}
\xi_{\beta}\co\ \beta M\to B\Gamma_{\beta}, \qquad \text{and} \qquad
\xi_{\beta}'\co\ \beta M'\to B\Gamma_{\beta},
\end{gather*}
such that $\xi_*\co \pi_1 M\to \Gamma$ and $(\xi_{\beta})_*\co
\pi_1\beta M \to \Gamma_{\beta}$ are isomorphisms.
Then there exists an $(L,G)$-bordism
 \begin{gather*}
 \big(W_{\Sigma},\bar \xi, \bar \xi_{\beta}\big)\co\
 (M_{\Sigma},\xi,\xi_{\beta}) \rightsquigarrow (M_{\Sigma}',\xi',\xi_{\beta}'),
\end{gather*}
as in Definition~$\ref{def:relbordism-2}$
with the additional property that the pairs
$(\beta W, \beta M)$ and $(W, M)$ are $2$-con\-nected.
\end{Theorem}

\begin{proof}
Start with an $(L,G)$-bordism
$W_{\Sigma}\co M_{\Sigma}\rightsquigarrow M_{\Sigma}'$
in the sense of Definition~\ref{def:relbordism-2}. We~have a commutative diagram
$$
\begin{diagram}
\setlength{\dgARROWLENGTH}{1.95em}
 \node{}\node{}
\node{\pi_1(\beta W)}\arrow{s} \arrow[2]{se}\\
\node{}\node{}\node{\pi_1(W)}\arrow{se}\\
\node{\pi_1(\beta M)}\arrow{e,l}{\theta}\arrow[2]{ne}
\node{\pi_1(M)}\arrow[2]{e,l}{\cong}\arrow{ne}
\node{}\node{\Gamma}
\node{\Gamma_{\beta},}\arrow{w,l}{\theta}
\end{diagram}
$$
where the composite $\pi_1(\beta M)\to \pi_1(\beta W)\to \Gamma_{\beta}$
(along the top triangle) is
an isomorphism, and thus the homomorphism
$\pi_1(\beta W)\to \Gamma_{\beta}$ is surjective. The first step is to do
surgery on~$\beta W$ to make the pair $(\beta W, \beta M)$
$2$-connected. We~begin by killing the kernel of the surjection
$\pi_1(\beta W)\twoheadrightarrow \Gamma_{\beta}$, just as in the proof of
Theorem~\ref{thm:surg1a}, so that the inclusion of $\beta
M\hookrightarrow \beta W$ induces an isomorphism on $\pi_1$. Once this
is done, we have an exact sequence
\begin{gather*}
\pi_2(\beta M) \to \pi_2(\beta W) \to \pi_2(\beta W, \beta M)
\xrightarrow{0}
\pi_1(\beta M) \xrightarrow{\cong} \pi_1(\beta W).
\end{gather*}
Then we need to do surgery on certain embedded $2$-spheres in $\beta W$,
generating $\pi_2(\beta W, \beta M)$. These embedded 2-spheres
have trivial normal bundle because of the spin assumption, and map
null-homotopically to the classifying space $BG$
(since the Lie group
$G$ is simply connected, hence automatically
$2$-connected~-- $\pi_2$ of any Lie group vanishes~-- and
thus $BG$ is $3$-connected). So the necessary
surgery is possible. The process for making the pair $(W, M)$
$2$-connected is completely analogous.
\end{proof}
The following result is a straightforward
generalization of known bordism results for psc metrics. Here Theorem~\ref{thm:surg2} provides all the necessary tools.
\begin{Theorem}[bordism theorem, general version]
\label{thm:bordism2}
Assume $G$ is a simply connected Lie group,
$L^{\ell}$ is simply connected and spin, $\theta\co\Gamma_{\beta}\to
\Gamma$ is a group homomorphism, and $n\ge 6 + \ell$. Let $x\in
\Omega^{\spin,(L,G)\fb}_{n}\big(B\Gamma_{\beta}\xrightarrow{\theta_*}B\Gamma\big)$
be represented by two $(L,G)$-singular spin pseudomanifolds
$(M_{\Sigma},\xi,\xi_{\beta})$, $M_{\Sigma}=M\cup_{\p M}
N(\beta M)$, and $(M_{\Sigma}',\xi',\xi_{\beta}')$,
$M'_{\Sigma}=M'\cup_{\p M} N(\beta M')$ given together with maps
\begin{gather*}
\xi \co\ M\to B\Gamma, \qquad
\xi_{\beta}\co\ \beta M\to B\Gamma_{\beta} \qquad\mbox{and} \qquad
\xi' \co\ M'\to B\Gamma, \qquad
\xi_{\beta}'\co\ \beta M'\to B\Gamma_{\beta}
\end{gather*}
such that $\xi_*\co \pi_1 M_{\Sigma}\to \Gamma$ and $(\xi_{\beta})_*\co
\pi_1\beta M \to \Gamma_{\beta}$ are isomorphisms. Furthermore, assume
$M'_{\Sigma}$ has an adapted psc~metric $g'$. Then there exists an
$(L,G)$-bordism
 \begin{gather*}
 \big(W_{\Sigma},\bar \xi, \bar \xi_{\beta}\big)\co\
 (M_{\Sigma},\xi,\xi_{\beta}) \rightsquigarrow (M_{\Sigma}',\xi',\xi_{\beta}'),
\end{gather*}
together with an adapted psc~metric $\bar g$ which is a product metric
near the boundary $\delta W_{\Sigma}= M_{\Sigma}\sqcup - M_{\Sigma}'$
such that $\bar g|_{M_{\Sigma}'}=g'$. In particular, $M_{\Sigma}$
admits an adapted psc~metric $g$.
\end{Theorem}

\begin{proof}
This proceeds like the proof Theorem~4.6 from Part I,
using Theorem~\ref{thm:surg2}.
\end{proof}

\section{The existence theorem in the non-simply connected case}
\label{sec:existence-GLR}
In this section we will use the results from Section~\ref{sec:main} to
obtain existence theorems for well-adapted {\psc} metrics on singular
pseudomanifolds in suitable cases. These results are parallel to
those in~\cite[Section~6]{BB-PP-JR}, but without the assumption that $M$ and
$\p M$ are simply connected. However, we will have to assume
something about the fundamental groups. We~denote by
$\Bott^8$ a simply connected spin $8$-manifold with
$\widehat A$-genus $1$, which one can take to be Ricci-flat~\cite{Joyce}.
This is often called a ``Bott manifold'', because taking the
product with it in~spin bordism corresponds, after applying the
$\alpha$-transformation to $KO_*$, to the Bott periodicity map
$KO_*\xrightarrow{\cong}KO_{*+8}$.
Let $\Gamma$ be a discrete group. We~need the following definitions.
 \begin{enumerate}\itemsep=0pt
\setlength{\leftskip}{0.55cm}
 \item[(GLR)] A group $\Gamma$ satisfies the
 \emph{Gromov--Lawson--Rosenberg conjecture} (GLR conjecture) if a~closed connected spin manifold $M$ with dimension $n\ge 5$ and
 fundamental group~$\Gamma$ admits a psc metric if and only if the
 generalized index $\alpha^\Gamma(M)\in KO_n(C^*_{r,\bR}(\Gamma))$ of
 the Dirac operator (i.e., the
 Atiyah--Singer--Mishchenko operator) on $M$ vanishes.
\item[($^s$GLR)] A group $\Gamma$ satisfies the \emph{stable
 Gromov--Lawson--Rosenberg conjecture} ($^s$GLR conjecture) if
 for any closed connected spin
 manifold $M$ with fundamental group $\Gamma$,
 the vanishing of the generalized index $\alpha^\Gamma(M)\in
 KO_n(C^*_{r,\bR}(\Gamma))$ implies that $M\times \big(\Bott^8\big)^k$
 admits a~psc metric for some sufficiently large $k$.
\end{enumerate}
We recall that Stolz (\cite{MR1403963} and~\cite[Section~3]{MR1937026})
sketched a proof that the $^s$GLR conjecture holds whenever the
Baum--Connes assembly map is injective, which is true for a very large
class of groups (conjecturally, all groups!). The GLR conjecture is
satisfied for a more restricted class of groups, but including free
groups, fundamental groups of oriented surfaces, and free abelian
groups (as long as one takes $n$ bigger than the rank of the group). It~also holds~\cite{MR1484887} for finite groups with periodic
cohomology, which includes finite cyclic groups and quaternion groups.
But it fails for some groups~\cite{MR2048721,MR1632971}.
Our first
main result is a generalization of~\cite[Theorem~6.3]{BB-PP-JR}. Recall the setting:

We fix a simply connected
homogeneous space $L = G/H$, where $G$ is a compact semisimple Lie
group acting on $L$ by isometries of the metric $g_L$, where
$\scal_{g_L}=\scal_{S^{\ell}}=\ell(\ell-1)$, where $\ell=\dim L$.
\begin{Theorem}
\label{thm:existenceGLR-Lbdy}
Let $M_\Sigma$ be an $n$-dimensional compact pseudomanifold with
resolution $M$, a spin manifold with boundary $\p M$. Assume the following:
\begin{enumerate}\itemsep=0pt
\item[{\rm (1)}] $M$ is a spin manifold with boundary $\p M$
 fibered over a connected spin manifold $\beta M$,
\item[{\rm (2)}] the fiber bundle $\phi\co \p M\to \beta M$ is a
 geometric $(L,G)$-bundle.
\end{enumerate}
Let $\Gamma = \pi_1(M)$, $\Gamma_\beta = \pi_1(\beta M)$. Furthermore,
assume $n\geq \ell+6$ and that the following
condition holds:
\begin{enumerate}\itemsep=0pt
\item[$\bullet$] the link $L$ is a spin $G$-boundary of a $G$-manifold
 $\bar L$ equipped with a psc~metric $g_{\bar L}$
 which is a product near the boundary and satisfies $g_{\bar L}|_L=g_L$.
\end{enumerate}
Then, provided both groups $\Gamma$ and $\Gamma_{\beta}$ satisfy
GLR conjecture, the vanishing of
the invariants
\begin{gather*}
\alpha^\Gamma_{\rm cyl}(M,g_M)\in
KO_n\big(C^*_{r,\bR} (\Gamma)\big) \qquad\text{and}\qquad
\alpha^{\Gamma_\beta}(\beta M) \in KO_n\big(C^*_{r,\bR} (\Gamma_\beta)\big)
\end{gather*}
implies that $M_\Sigma$ admits an adapted psc~metric.
\end{Theorem}
\begin{proof}
Before we start with the proof we
emphasize that under the assumptions above the
index $\alpha^\Gamma_{{\rm cyl}} (M)$ is well defined,
independent of the chosen adapted wedge metric $g$; see Remark~\ref{dependance-metric} and~Proposition~\ref{prop:bordism-map-index}. According to
Theorems~\ref{thm:bordism2} and~\ref{prop:bordismtriangle}, it
suffices to show that the bordism class\footnote{Here we use a
 short notation $M_\Sigma$ instead of $(M_{\Sigma},\xi,\xi_{\beta})$.}
\begin{gather*}
[M_\Sigma] \in
 \Omega_n^{\spin,
 (L,G)\fb}\big(B\Gamma_\beta\xrightarrow{B\theta}B\Gamma\big)
\end{gather*}
is a sum of classes of pseudomanifolds with well-adapted {\psc}
metrics. We~use the exact sequence~\eqref{eq:bordismtriangle}. Since
$\alpha^{\Gamma_\beta}(\beta M)=0$ and $\Gamma_\beta$ satisfies the
GLR conjecture, and since $\dim \beta M = n-\ell-1 \ge 5$,
the manifold $\beta M$ admits a psc
metric. Then, since the fibration $\phi \co \p M\to \beta M$
is a geometric $(L,G)$-bundle; we can identify $\p M$ with
$P\times_{G}L$, where $P\to \beta M$ is the corresponding principal
$G$-bundle. We~use a bundle connection on $P$ to construct a
well-adapted psc metric on the tubular neighborhood $N$ of $\beta M$
in $M_\Sigma$. In particular, we obtain a psc metric
on the boundary $\p N$, where the $(L,G)$-fiber bundle $\p N \to\beta M$
is also a~Riemannian submersion.

We use the $G$-manifold $\bar L$ which bounds $L$ to
construct an $\bar L$-bundle $M'=P\times_G \bar L$ over $\beta M$
associated to the above principal $G$-bundle $P\to \beta M$.
By assumptions, $\bar L$ is given a psc metric $g_{\bar L}$
which is a product metric near $\p \bar L=L$. Then $M'$ has a
bundle metric of {\psc}, and joining $M'$ to $N$, we get an
$(L,G)$-singular spin manifold
\begin{gather*}
M'_\Sigma=M'\cup_{\p M} -N
\end{gather*}
with a well-adapted psc metric.

Since the pseudomanifolds
$M'_\Sigma$ and $M_\Sigma$ coincide near $\beta M$, by
\eqref{eq:bordismtriangle} their bordism classes differ by a class
$[M'']$ in the image of $\Omega_n^\spin(B\Gamma)$.
Since, by assumption, we have $\alpha_{{\rm cyl}}^\Gamma(M)=0$,
and also $\alpha_{{\rm cyl}}^\Gamma(M')=0$
(since $M'$ has a psc metric), we have by additivity of the
$\alpha$-invariant (Proposition~\ref{prop:gluing}) that
$\alpha^\Gamma(M'')=0$, and so we can take
$M''$ to be a closed spin $n$-manifold
with fundamental group $\Gamma$ equipped with some psc metric.
Now $M_{\Sigma}$ is in the same bordism class in
\begin{gather*}
\Omega_n^{\spin, (L,G)\fb}\big(B\Gamma_\beta\xrightarrow{B\theta}B\Gamma\big)
\end{gather*}
as
$M'_{\Sigma}\sqcup M''$, and we
can apply Theorem~\ref{thm:bordism2} to get the conclusion.
\end{proof}
\begin{Theorem}
Let a link $L$ and pseudomanifold $M_\Sigma$ be as in \textup{Theorem~\ref{thm:existenceGLR-Lbdy}}, but this time only assume that the
groups $\Gamma$ and $\Gamma_\beta$ satisfy the $^s\!$GLR conjecture.
Then $M_\Sigma$ stably admits a~well-adapted metric of positive scalar
curvature, i.e., $M_\Sigma\times \big(\text{\textup{Bott}}^8\big)^k$ admits
such a metric for some~$k\ge 0$.
\end{Theorem}

\begin{proof}
The proof of this is exactly the same as for Theorem~\ref{thm:existenceGLR-Lbdy}, the only difference being that first
we need to cross with some copies of the Bott manifold to get
{\psc} on~$\beta M$, and then we might need to cross with additional
copies of the Bott manifold to get {\psc} on $M''$.
\end{proof}
\begin{Remark}
 It is possible without great effort to adapt the arguments above to
 the case where~$\beta M$ is disconnected. We~leave details to the
 reader.
\end{Remark}
Now we present a generalization of~\cite[Theorem~6.7]{BB-PP-JR}. It~was shown in~\cite[Lemma 6.6]{BB-PP-JR} that the class of
$\bH\bP^{2k}$, $k\ge 1$, is not a zero-divisor in the spin bordism
ring $\Omega_*^\spin$. It~then follows that if~$H_*(B\Gamma; \bZ)$ is
torsion free and the Atiyah-Hirzebruch spectral sequence converging to
$\Omega_*^\spin(B\Gamma)$ collapses, so that $\Omega_*^\spin(B\Gamma)$
is a free $\Omega_*^\spin$-module
$\Omega_*^\spin\otimes_{\bZ}H_*(B\Gamma; \bZ)$, then the class of
$\bH\bP^{2k}$ does not annihilate any non-zero class in
$\Omega_*^\spin(B\Gamma)$. This condition on $\Gamma$ is satisfied if~$B\Gamma$ is stably homotopy-equivalent to a wedge of spheres~-- for
example, if~$\Gamma$ is a free group or free abelian group.
\begin{Theorem}\label{thm:existenceGLR-HPeven}
Let $(M_{\Sigma},\xi,\xi_{\beta})$ be a closed
$(L,G)$-singular spin pseudomanifold, $\dim M_{\Sigma}=n$, with
$L=\bH\bP^{2k}$ and $G=\Sp(2k+1)$, $k\ge 1$. Assume
that
\begin{enumerate}\itemsep=0pt
\item[$(i)$] the manifolds $M$ and $\beta M$ are connected and the
 inclusion $\p M\hookrightarrow M$ induces an isomorphism on
 fundamental groups with $\pi_1 M \cong \pi_1 \p M =\Gamma$,
\item[$(ii)$] the $(L,G)$-bundle $\p M \to \beta M$ is trivial, i.e.,
 $\p M = \beta M \times L$,
\item[$(iii)$] the group $\Gamma$ satisfies the GLR conjecture,
\item[$(iv)$] $\Omega_*^\spin(B\Gamma)$ is free as an
 $\Omega_*^\spin$-module.
\end{enumerate}
Then if $n-8k\ge 6$, the pseudomanifold $M_\Sigma$ has an adapted psc
metric if and only if the $\alpha$-invariant $\alpha_{{\rm
 cyl}}^{\Gamma}(M) \in KO_n\big(C^*_{r,\bR}(\Gamma)\big)$ vanishes.
\end{Theorem}
\begin{Remark}
Note that in this case we have isomorphisms $\pi_1 M \cong \pi_1 \p M
\cong \pi_1 \beta M$. Also recall that if the bundle $\p M \to \beta M$
is trivial, then
the singularities are of Baas--Sullivan type.
\end{Remark}
\begin{proof}
Necessity was proved in Section~3, so we need to prove existence of an
adapted metric of {\psc} assuming the vanishing condition. Since $\p
M = \beta M \times L$ is a spin boundary over $B\Gamma$ (namely, it is
the boundary of $M$), its class is trivial in
$\Omega_*^\spin(B\Gamma)$. But the class of $L$ cannot annihilate any
non-zero class in $\Omega_*^\spin(B\Gamma)$, by~\cite[Lemma~ 6.6]{BB-PP-JR}, so in fact $\beta M$ must be a spin boundary over
$B\Gamma$. (In particular, $\alpha^{\Gamma}(\beta M) $ automatically
vanishes.) Now by~\cite[Proposition~4.4]{BB-PP-JR}, or equivalently, by the
proof of middle-exactness in the proof of~Theorem~\ref{prop:bordismtriangle} (this part of the proof doesn't require $L$
to be a spin boundary), $M_\Sigma$ is bordant over $B\Gamma$ to a
closed manifold~$M'$ with fundamental group $\Gamma$ and no
singularities. Then we must have $\alpha^\Gamma(M')=\alpha_{{\rm
 cyl}}^{\Gamma}(M) = 0$, and since $\Gamma$ satisfies the GLR
conjecture, $M'$ admits a metric of~{\psc}. The conclusion now
follows from Theorem~\ref{thm:bordism1}.
\end{proof}
\begin{Theorem}
Let $(M_{\Sigma},\xi,\xi_{\beta})$ be a closed
$(L,G)$-singular spin pseudomanifold, $\dim M_{\Sigma}=n$, with
$L=\bH\bP^{2k}$ and $G=\Sp(2k+1)$, $k\ge 1$. Assume
that
\begin{enumerate}\itemsep=0pt
\item[$(i)$] the manifolds $M$ and $\beta M$ are connected and the
 inclusion $\p M\hookrightarrow M$ induces an isomorphism on
 fundamental groups with $\pi_1 M \cong \pi_1 \p M =\Gamma$,
\item[$(ii)$] the $(L,G)$-bundle $\p M \to \beta M$ is trivial, i.e.,
 $\p M = \beta M \times L$,
\item[$(iii)$] the group $\Gamma$ satisfies the $^s\!$GLR conjecture,
\item[$(iv)$] $\Omega_*^\spin(B\Gamma)$ is free as an
 $\Omega_*^\spin$-module.
\end{enumerate}
Then the pseudomanifold $M_\Sigma\times \big(\text{\textup{Bott}}^8\big)^m$
has an adapted metric of {\psc} for some
$m\ge 0$ if and only if the $\alpha$-invariant $\alpha_{{\rm
 cyl}}^{\Gamma}(M) \in KO_n(C^*_{r,\bR}(\Gamma))$ vanishes.
\end{Theorem}
\begin{proof}
This is exactly the same as the proof of Theorem~\ref{thm:existenceGLR-HPeven}, except that we need to multiply~$M'$ by some product of copies of the Bott manifold to get {\psc}.
\end{proof}

\section[The space Rw+(MSigma) of adapted wedge metrics of positive scalar curvature]
{The space $\boldsymbol{\mathcal{R}_w^+ (M_\Sigma)}$
 of adapted wedge metrics \\of positive scalar curvature}
\label{sec:secondary}
Let $M_\Sigma$ be a $(L,G)$-fibered pseudomanifold which admits a
well-adapted wedge psc metric. In~this section we shall study the space
$\cR_w^+ (M_\Sigma)$ of well-adapted wedge psc metrics on $M_\Sigma$.
Our first goal is to establish a relationship between the
space $\cR_w^+ (M_\Sigma)$ and corresponding spaces
$\cR^+ (\beta M)$, $\cR^+ (\p M)$, $\cR^+(M)$ of psc
metrics. We~consider the map
$\resS\co \mathcal{R}_w^+ (M_\Sigma)\to \mathcal{R}^+ (\beta M)$
which takes an adapted wedge metric $g$ on
$M_{\Sigma}$ to its restriction $g_{\beta M}$ on $\beta M$,
and we will observe that this map is a Serre fiber bundle. We~show (this was stated before as Theorem~\ref{thm:section})
that there exists a non-canonical section
$s\co \mathcal{R}^+ (\beta M)\to \cR_w^+(M_\Sigma)$;
this fact, together with recent results on the homotopy
type of the space of psc metrics on spin manifolds,
\cite{BER-W,ER-W2}, implies that the homotopy groups of the space
$\mathcal{R}_w^+ (M_\Sigma)$ are rationally at least as complicated
as the groups $KO_*\big(C^*_r(\pi_1(\beta M))\big)\otimes \bQ$.
Next, we elaborate on a \emph{direct} definition
of the index-difference homomorphism which enables us to detect
some nontrivial classes in $\pi_q\big(\mathcal{R}_w^+ (M_\Sigma)\big)$.
Finally, if the fundamental group of $M_\Sigma$ has an element of
finite order, we show that the Cheeger--Gromov rho invariant on a
depth-1 pseudomanifold, introduced and studied in~\cite{Piazza-Vertman}, can be used to show that $\big\vert \pi_0\big(\cR_w^+ (M_\Sigma)\big) \big\vert=\infty$.

\subsection[Controlling the homotopy type of Rw+(MSigma)]{Controlling the homotopy type of $\boldsymbol{\mathcal{R}_w^+ (M_\Sigma)}$}
We fix an $(L,G)$-fibered compact pseudomanifold
$M_\Sigma= M\cup_{\p M} N(\beta M)$. Here $M$ is a compact manifold
with boundary $\p M$ which is a geometric $(L,G)$-bundle over $\beta M$.
Recall that the $G$ acts on $L$ by isometries of a given metric
$g_L$, and the $(L,G)$-bundle $\phi \co \p M\to \beta M$ is given by a
principal $G$-bundle $p\co P \to \beta M$ so that $\p M = P\times_G L$.
Then we identify the space $\mathcal{R}_w^+ (M_\Sigma)$ of
well-adapted wedge psc metrics on $M_\Sigma$ with the space of
associated triples $(g_M, g_{\beta M}, \nabla^p)$, where $g_M$ and
$g_{\beta M}$ are psc metrics on $M$ and $\beta M$ respectively, and
$\nabla^p$ is a connection on the principal bundle $P$. By
definition of well-adapted wedge metrics, the triples
$(g_M, g_{\beta M}, \nabla^p)$ are subject to the following conditions:
\begin{enumerate}\itemsep=0pt
\item[$(i)$] $g_M={\rm d}t^2 + g_{\p}$ near the boundary $\p M\subset M$,
 for some Riemannian metric $g_\p$ on $\partial M$;
\item[$(ii)$] the bundle map $\phi$ gives a Riemannian submersion
 $\phi\co (\p M,g_{\p}) \to (\beta M, g_{\beta M})$, and the
 connection $\nabla^p$ on $P$ defines a $G$-connection $\nabla^\phi$
 for the submersion.
\end{enumerate}
We denote by $\mathcal{R}^+(M,\p M)$ the space of Riemannian psc
metrics $g$ on $M$ such that $g={\rm d}t^2 + g_{\p}$ near the boundary $\p M$,
for some metric $g_{\p}$ on $\p M$. There is an obvious restriction
map
\begin{gather}\label{eq-res}
 \res\co \ \mathcal{R}^+(M,\p M) \to \mathcal{R}^+(\p M),
\end{gather}
which is known to be a Serre fiber bundle, see~\cite[Theorem~1.1]{Ebert-Frenck}. For a given metric $h\in \mathcal{R}^+(\p M)$,
we denote by $\mathcal{R}^+(M)_{h} = \res^{-1}(h)$ the
corresponding fiber of (\ref{eq-res}).

We notice that there is a map
$\iota_L\co \mathcal{R}^+(\beta M) \to \mathcal{R}^+(\p M)$
which is given by lifting a met\-ric~$g_{\beta M}$ on $\beta M$ to a
Riemannian submersion metric $g_{\p M}$ on the total space $\p M$
of the $(L,G)$-bundle $\phi\co \p M \to \beta M$, using the
connection $\nabla^\phi$ of (ii) above,
by putting the metric $g_L$ on each fiber~$L$. This map
is injective since the metrics $g_{\beta M}$ and $\iota_L(g_{\beta M})$
determine one another.
Then, by definition, a wedge psc metric $g$ on
$M_{\Sigma}$ determines a unique psc metric $g_{\beta M}\in
\mathcal{R}^+(\beta M)$. This gives a well-defined map
$\resS\co \mathcal{R}_w^+(M_\Sigma)\to
\mathcal{R}^+(\beta M)$. Thus we
see that the space $\mathcal{R}_w^+ (M_\Sigma)$ is a~pull-back in
the following diagram:
\begin{gather}\label{pull-back}
 \begin{diagram}
 \setlength{\dgARROWLENGTH}{1.95em}
\node{\mathcal{R}_w^+ (M_\Sigma)}
 \arrow{e,t}{\resS}
 \arrow{s}
\node{\mathcal{R}^+(\beta M)}
 \arrow{s,l,J}{\iota_L}
 \\
\node{\mathcal{R}^+(M,\p M)}
 \arrow{e,t}{\res}
\node{\,\mathcal{R}^+(\p M).}
\end{diagram}
\end{gather}
As we mentioned above, the restriction map
$\res\co\mathcal{R}^+(M,\p M)\to \mathcal{R}^+(\p M)$ is a Serre fiber
bundle, and it is easy to see that the map
$\resS\co \cR_w^+ (M_\Sigma)\to \cR^+(\beta M)$ from
(\ref{pull-back}) is also a Serre fiber bundle.\footnote{\ This
 requires a word of explanation, as follows: by~\cite[Theorem~3.5(2)]{BB-PP-JR}, given a metric $g_\beta \in \mathcal{R}^+(\beta
 M)$ such that $\iota_L(g_\beta)$ extends to a metric in
 $\mathcal{R}^+(M,\p M)$, then some rescaling of $g_\beta$ (by a
 constant conformal factor) lifts to $\mathcal{R}_w^+(M_\Sigma)$.}
Let $g_{\beta M}\in \mathcal{R}^+(\beta M)$ and $g_{\p M}\in \cR^+(\p M)$
be such that $\iota_L(g_{\beta M})=g_{\p M}$. We~have the following commutative diagram of fiber bundles:
\begin{gather}\label{pull-back2}
 \begin{diagram}
 \setlength{\dgARROWLENGTH}{1.95em}
 \node{\mathcal{R}_w^+ (M_\Sigma)_{g_{\beta M}}}
 \arrow{e}
 \arrow{s,l}{\simeq}
 \node{\mathcal{R}_w^+ (M_\Sigma)}
 \arrow{e,t}{\resS}
 \arrow{s}
\node{\mathcal{R}^+(\beta M)}
 \arrow{s,l}{\iota_L}
 \\
 \node{\mathcal{R}^+ (M)_{g_{\p M}}}
 \arrow{e}
 \node{\mathcal{R}^+(M,\p M)}
 \arrow{e,t}{\res}
\node{\,\mathcal{R}^+(\p M),}
\end{diagram}
\end{gather}
where $\cR_w^+ (M_\Sigma)_{g_{\beta M}}$ and $\cR^+(M)_{g_{\p M}}$ are
the corresponding fibers. Note that the downward map on the left is a
homotopy equivalence, in fact a homeomorphism, since any element of
$\cR^+(M)_{g_{\p M}}$ defines a unique psc wedge metric.

In particular, if we have an isotopy $(g_{\beta M})_t\in \cR^+(\beta
M)$, and $(g_{\beta M})_0 = \resS(g_0)$, then the isotopy $(g_{\beta
 M})_t$ lifts to an isotopy $g_t$ in $\cR_w^+(M_\Sigma)$ with
$\resS(g_t)=(g_{\beta M})_t$.

\begin{Theorem}
\label{thm:pi0ofbdyinjects}
Let $M_\Sigma$ be an $(L,G)$-fibered compact pseudomanifold with $L$ a
simply connected homogeneous space of a compact semisimple Lie group.
Also assume that $M_\Sigma$ admits an adapted wedge metric of {\psc}.
Then there exists a section
$s\co \mathcal{R}^+(\beta M)\to \mathcal{R}_w^+ (M_\Sigma)$
to $\resS$. In particular, there is an injection of homotopy groups
\begin{gather*}
 s_* \co\ \pi_q\big(\mathcal{R}^+ (\beta M)\big)\to \pi_q\big(\mathcal{R}_w^+ (M_\Sigma)\big), \qquad
 q= 0, 1,\ldots.
\end{gather*}
\end{Theorem}

\begin{proof}
We want to construct a (non-canonical) section to the ``forgetful map''
\begin{gather*}
\resS\co \cR_w^+(M_\Sigma)\to \cR^+(\beta M),
\end{gather*}
sending a psc metric $g$ (interpreted as a triple $(g_M, g_{\beta M},
\nabla^p)$, as above) to the metric $g_{\beta M}$. As~before, we
will fix the connection $\nabla^p$ once and for all, and we fix a
metric $g_{\beta M}\in \cR^+(\beta M)$. Given another metric
$g'_{\beta M}$ in $\cR^+(\beta M)$, we consider the linear
homotopy $g_{\beta M}(t)= (1-t)g_{\beta M}+ t g'_{\beta M}$ within the
space $\cR(\beta M)$ of all Riemannian metrics. \emph{This
homotopy may go out of the sub\-space~$\cR^+ (\beta M)$}, but the
scalar curvature function along the homotopy is bounded below by
some constant $-c$, $c > 0$. If we scale the metrics
$g_{\beta M}(t)$ by $\lambda^2$, then the scalar curvature scales by~$\scal_{\lambda^2g_{\beta M}(t)}= \lambda^{-2}\scal_{g_{\beta M}(t)}$. We~consider the family of metrics
$g_{\p M}(t)=\iota_L\big(\lambda^2g_{\beta M}(t)\big)$ on~$\p M$.
Recall that $\scal_L=\ell(\ell-1)>0$, so we obtain
\begin{gather*}
\scal_{g_{\p M}(t)}= \scal_{g_L} + \lambda^{-2}\scal_{g_{\beta M}(t)}
\geq \ell(\ell-1)- \lambda^{-2}c.
\end{gather*}
Clearly, there exists $\lambda>0$ such that $\scal_{g_{\p M}(t)}>0$
for all $t>0$. Then we lift the curve of met\-rics~$g_{\p M}(t)$ to a
curve of met\-rics~$g(t)\in \cR^+(M, \p M)$
using homotopy lifting in the fiber bundle
\begin{gather*}
\res\co\ \mathcal{R}^+(M,\p M) \to \mathcal{R}^+(\p M).
\end{gather*}
In particular, we obtain an adapted wedge metric $(g_M',g_{\beta
 M}',\nabla^p) \in \mathcal{R}_w^+(M_\Sigma)$. It~is important to
notice the following:
\begin{enumerate}\itemsep=0pt
\item[1.] The curve of the adapted wedge metrics $(g_M(t),g_{\beta
 M}(t),\nabla^p)$ is not, in general, a curve in the space
 $\cR_w^+(M_\Sigma)$. Indeed, the scalar curvature on the
 tubular neighborhood $N(\beta M)$ is dominated by the sum of
 zero scalar curvature on the cone over $L$ and by the scalar curvature
 $\lambda^{-2}\scal_{g_{\beta M}(t)}$.\footnote{An important comment:
 to apply the argument of~\cite[Theorem~3.5(2)]{BB-PP-JR} in order
 to rescale a metric $g'_{\beta M}$ in~$\mathcal{R}^+(\beta M)$ so
 that $\iota_L(g'_{\beta M})$ extends to a psc metric on the whole
 tubular neighborhood $N$ of $\beta M$, one needs to make use of a
 positive lower bound on the scalar curvature. So while $g'_{\beta
 M}$ lifts to $\mathcal{R}_w^+(M_\Sigma)$ (after rescaling by a~positive constant), the same is \emph{not} true for the other
 metrics in the path in $\mathcal{R}(\beta M)$ from $g_{\beta M}$
 to $g'_{\beta M}$.}
\item[2.] The correspondence
 $g'_{\beta M}\mapsto (g'_M,g_{\beta M}',\nabla^p)$
 can be made continuous (in the $C^2$ topology)
 as long as the metrics $g'_{\beta M}$ stay inside
 a compact set $K\subset \mathcal R^+(\beta M)$.
\end{enumerate}
This is enough to get a (non-canonical) section
$s_{g_{\beta M}}\co \mathcal{R}^+(\beta M)\to \mathcal{R}_w^+(M_\Sigma)$
of the fiber bundle $\resS$ in~\eqref{pull-back}.
\end{proof}

\begin{Corollary}\label{short-exact}
The fiber bundle~\eqref{pull-back2} gives a split short exact sequence:
\begin{gather}\label{eq:sgor-exact}
 0 \to \pi_q\big(\mathcal{R}_w^+ (M_\Sigma)_{g_{\beta
 M}}\big)\xrightarrow{i_*} \pi_q\big(\mathcal{R}_w^+
 (M_\Sigma)\big)\xrightarrow{(\mathrm{res}_{\Sigma})_*}\pi_q\big(\mathcal{R}^+(\beta
 M)\big)\to 0
 \\ \hphantom{0 \to }
 \mbox{for each} \quad q=0,1,\ldots.\nonumber
\end{gather}
\end{Corollary}

\begin{proof}
Indeed, a section $s_{g_{\beta M}}\co \mathcal{R}^+(\beta M)\to
\mathcal{R}_w^+(M_\Sigma)$ from Theorem~\ref{thm:pi0ofbdyinjects}
implies that the homo\-morphism $(\mathrm{res}_{\Sigma})_*$ is
split surjective for all $q=0,1,\ldots$.
\end{proof}
An immediate consequence of this corollary is that
$\pi_q(\cR_w^+(M_\Sigma))$ admits $\pi_q\big(\cR^+(\beta M)\big)$
as a direct summand. In particular, $\pi_0(\cR_w^+(M_\Sigma))$
contains $\pi_0\big(\cR^+(\beta M)\big)$ as a direct summand.

\begin{Example}
Here is an interesting example of obtaining information about
$\pi_0\big(\mathcal{R}^+(M_\Sigma)\big)$. Let $L=G=\SU(2)=S^3$. Since
$\pi_7(BS^3) = \pi_6\big(S^3\big) = \bZ/12$, there are $12$ distinct principal
$S^3$-bundles over $\beta M = S^7$, all of them of the rational
homotopy type of $S^7\times S^3$. For any one of~these $S^3$-bundles
over $S^7$, the total space is $\p M = \p N$,
where $M=N$ is the disk bundle of our~$S^3$-bundle over $S^7$.
Define $M_\Sigma=N\cup_{\p M}-N$, the double of $N$.
Note that $N$ is actually the unit disk bundle of a quaternionic
line bundle, so our $M_\Sigma$ in this case is in fact a smooth
manifold, though not every Riemannian metric on $M_\Sigma$ is
adapted for the $(L,G)$-fibered structure, and we are only interested
in this special subclass of metrics. We~know that $\pi_0\big(\cR_w^+(M_\Sigma)\big)$
contains $\pi_0\big(\cR^+(\beta M)\big)=\pi_0\big(\cR^+\big(S^7\big)\big)$ as a direct summand.
By~\cite[Theorem~4.47]{MR720933}, $\pi_0\big(\cR^+\big(S^7\big)\big)$ contains a copy of
$\bZ$, and the same argument shows $\pi_0\big(\cR^+\big(S^{11}\big)\big)$ contains a copy of
$\bZ$. By taking a~connected sum of $(M_\Sigma, g_\Sigma)$ with
$\big(S^{11},g_k\big)$ for suitable metrics $g_k$ on $S^{11}$ (the connected
sum takes place on the interior of $M$, so it doesn't change the
adapted wedge structure of the metric near $\beta M=S^7$), we see
that $\pi_0\big(\cR_w^+(M_\Sigma)\big)$ contains $\bZ\oplus\bZ$, with one summand
coming from~$\beta M$ and one coming from the interior of $M$.
\end{Example}

\begin{Remark}
We recall that for a closed manifold $X$, the homotopy
type of the space $\mathcal{R}^+(X)$ does not change under
\emph{admissible} surgeries; in particular, this implies that the
homotopy type of~$\mathcal{R}^+(X)$ depends only on the bordism
class of $X$ in a relevant bordism group
\cite[Theorem~1.5]{Ebert-Frenck}. As we have seen,
on an $(L,G)$-pseudomanifold $M_{\Sigma}$, we can do two
types of surgeries: surgeries on the interior of the resolution $M$,
and surgeries on $\beta M$. It~turns out that the homotopy type of the space
$\mathcal{R}^+_w(M_{\Sigma})$ is also invariant with respect to
corresponding admissible surgeries,~\cite{BW}.
\end{Remark}

There are several methods for detecting non-trivial
elements in $\pi_q \big(\cR^+(\beta M)\big)$:
the index-difference homomorphisms and (higher) rho-invariants are certainly
two very efficient tools in~this direction. In the Sections~\ref{subsect:index-diff}
we shall elaborate further on these two tools. Notice
that some of the results proven by
rho-invariants can also be recovered using the index-difference
homomorphisms~\eqref{ind-diff2b} explained below; see
\cite[Remark 1.1.2]{ER-W2}.

\subsection{The index-difference homomorphism}
\label{subsect:index-diff}
First we recall some results from~\cite{BER-W,ER-W2}. Let $\bar X$ be
a compact closed spin manifold, $\dim \bar X= k\geq 5$. Let $\bar h_0\in
\mathcal{R}^+\big(\bar X\big)$ be a base point. Assume
there is a map $f\co \bar X \to B\Gamma$, where $\Gamma$ is a~discrete
group such that $f_* \co\pi_1\bar X \to \Gamma $ is split surjective.
(If $\pi_1\bar X$ is trivial, we assume $\Gamma=0$.)

We also need the case when the boundary $\p \bar X = X$ is non-empty,
with a given metric $h\in \mathcal{R}^+(X)$. As before, we denote by
$\mathcal{R}^+\big(\bar X\big)_h$ the subspace of metrics $\bar h$ in
$\mathcal{R}^+(\bar X,X)$ which restrict to $h$ on the boundary
$X$. Let $\bar h_0\in \mathcal{R}^+\big(\bar X\big)_h$ be a base point. Again,
we assume that there is a map $f\co \bar X \to B\Gamma$, such that
$f_* \co \pi_1 \bar X \to \Gamma$ is split surjective.
In both these settings
(without using split surjectivity of $f_*$ yet)
we have the index-difference homomorphisms
\begin{gather}
 \mathrm{inddiff}_{\bar h_0}^{\Gamma}\co\ \pi_q\mathcal{R}^+\big(\bar X\big)\to
 KO_{q+k+1}\big(C^*_r(\Gamma)\big), \nonumber
 \\
 \mathrm{inddiff}_{\bar h_0}^{\Gamma}\co\ \pi_q\mathcal{R}^+\big(\bar X\big)_h\to
 KO_{q+k+1}\big(C^*_r(\Gamma)\big).\label{ind-diff2b}
\end{gather}
Since it will be useful later, we remind the reader about where
these index difference maps
come from. Say for simplicity that $\bar X$ is a closed
manifold. A class in $\pi_q\cR^+\big(\bar X\big)$ is represented by
a~family $g_t$, $t\in S^q$, of psc metrics on $\bar X$, with
$g_{t_0}$, $t_0$ the basepoint in $S^q$, equal to our basepoint psc
metric. Let us denote here by $g_0$ this base-point metric.
From this family we get a warped product metric $g$
on $S^q\times \bar X$, that restricts on the copy of $\bar X$ over
$t\in S^q$ to $g_t$, and that is the usual flat (if $q\le 1$) or
round (if $q\ge 2$) metric on the copy of $S^q$ over each point in~$\bar X$. By the usual argument that ``isotopy implies concordance''
\cite[Lemma 3]{MR577131}, we can, without changing the homotopy
class of our map $t\mapsto g_t$ in $\pi_q\cR^+\big(\bar X\big)$, assume that
$g$ has positive scalar curvature on $S^q\times \bar X$.
Our class in $\pi_q\cR^+\big(\bar X\big)$ is trivial if and
only if it can be extended over the disk~$D^{q+1}$. So~extend $g$
to a metric on $D^{q+1}\times \bar X$ which restricts to a
product metric ${\rm d}r^2 + g$ on a~neighborhood of the boundary.
Extend it to a complete metric on $\bR^{q+1}\times \bar X$ by allowing~$r$, the distance to the center of the disk, to go to $\infty$,
and taking the metric to be ${\rm d}r^2 + g$ for all $r\ge 1$.
Since the metric has psc outside a compact set, the Dirac
operator on $\bR^{q+1}\times \bar X$ for this metric is Fredholm,
and we get an index in $KO_{q+1+k}\big(C^*_r(\Gamma)\big)$ as usual.
One way to see this property is to use $b$-calculus techniques; this
point of view will be useful later on in this section.
This index is an obstruction to triviality of our class in
$\pi_q\cR^+\big(\bar X\big)$, since it would be $0$ if we could extend
our family of metrics to a psc family over the disk.
Summarizing, we have defined a homomorphism
\begin{gather*}
 \mathrm{inddiff}_{g_{0}}^{\Gamma}\co\ \pi_q\mathcal{R}^+\big(\bar X\big)\to
KO_{q+k+1}\big(C^*_r(\Gamma)\big).
\end{gather*}
An alternative way to construct the index difference is the
following. Once we choose a base point $g_0=g_{t_0}\in
\cR^+\big(\bar X\big)$, then for any metric $g\in \cR^+\big(\bar X\big)$, there is a
linear path \mbox{$g_t=(1-t)g_0+ tg$} of~metrics in the space $\cR\big(\bar X\big)$
of all Riemannian metrics. Then we have a curve of the corresponding
Dirac operators $D_{g_t}$ with the ends of the curve, $D_{g_0}$ and
$D_{g}$ being in the contractible space of $C\ell_k$-linear invertible
operators. Such a curve determines a loop in the space $\mathbf{KO}_k$
which classifies $C\ell_k$-linear Fredholm operators (where the
space of invertible operators is collapsed to a point). One can easily
replace this by the classifying space $\mathbf{KO}_k\big(C^*_r(\Gamma)\big)$
for the $K$-theory $KO_{*+k}\big(C^*_r(\Gamma)\big)$ in the non-simply connected case.
Thus in homotopy we get the index-difference homomorphism
\begin{gather*}
 \mathrm{inddiff}_{g_{0}}^{\Gamma}\co\ \pi_q\cR^+\big(\bar X\big)\to
 \pi_q\big(\Omega\mathbf{KO}_k\big(C^*_r(\Gamma)\big)\big)=KO_{q+1+k}\big(C^*_r(\Gamma)\big),
\end{gather*}
which depends on a choice of the base point $g_0\in
\cR^+\big(\bar X\big)$.
The equivalence of these two app\-ro\-a\-ches to the index difference (with
a discussion of how to trace them back to the work of~Gromov--Lawson
and Hitchin, respectively) may be found in~\cite{Ebert2,Ebert1} (see also~\cite{Bugg}).

Next we recall the following relevant results:
\begin{Theorem}[see~\cite{BER-W,ER-W2}]
\label{thm:new-index-diff}
\leavevmode
\begin{enumerate}\itemsep=0pt
\item[$(a)$] Let $\bar X$ be a simply connected spin manifold,
 with $\p \bar X$ possibly non-empty, $\dim \bar X= k\geq 6$,
 with $\bar h_0\in \mathcal{R}^+\big(\bar X\big)\neq \varnothing$ $\big($or
 $\bar h_0\in \cR^+\big(\bar X\big)_h\neq \varnothing \big)$. Then both
 index-difference homomorphisms $\mathrm{inddiff}_{\bar
 h_0}^{\Gamma}$ from~\eqref{ind-diff2b} are non-trivial
 whenever the target group $KO_{q+k+1}$ is.
\item[$(b)$] Let $\bar X$ be a non-simply connected spin manifold,
 with $\p \bar X$ possibly non-empty, $\dim X= k\geq 6$ with
 $\bar h_0\in \cR^+\big(\bar X\big)\neq \varnothing $ $\big($or
 $\bar h_0\in \cR^+\big(\bar X\big)_h\neq \varnothing\big)$.
 Assume there is a map $f\co \bar X \to B\Gamma$,
 such that $f_* \co \pi_1 \bar X \to \Gamma$ is
 split surjective. Furthermore, we assume that
\begin{itemize}\itemsep=0pt
\item $\Gamma$ satisfies the rational Baum--Connes conjecture,
\item $\Gamma$ is torsion free and has finite rational homological
 dimension $d$, and
\item $q > d-2k-1$.
\end{itemize}
Then the images of both index-difference homomorphisms
$\mathrm{inddiff}_{\bar h_0}^{\Gamma}$ from~\eqref{ind-diff2b}
generate the target group $KO_{q+k+1}\big(C^*_r(\Gamma)\big)\otimes \bQ$
as a $\bQ$-vector space.
\end{enumerate}
\end{Theorem}

Recall also the following \emph{stable result}. As
in the beginning of Section~\ref{sec:existence-GLR}, we denote by
$\Bott^8$ a Bott manifold, which is a simply connected spin $8$-manifold
with $\hat{A}(\Bott)=1$. By the work of Joyce~\cite{Joyce}, we may
choose a Ricci flat metric $g_{\flat}$ on $\Bott$. Then for a closed
manifold $\bar X$ there are induced maps
\begin{gather*}
 \cR^+\big(\bar X\big) \xrightarrow{\times(\Bott,g_{\flat})}\cR^+\big(\bar X\times \Bott\big)
 \xrightarrow{\times(\Bott,g_{\flat})} \cR^+\big(\bar X\times \Bott\times \Bott\big)
 \xrightarrow{\times(\Bott,g_{\flat})}\cdots ,
\end{gather*}
and write $ \cR^+\big(\bar X\big)\big[\Bott^{-1}\big]$ for the homotopy
colimit. Similarly, if $\p \bar X = X\neq \varnothing $, and
$h\in \cR^+(X)$, we define the space $ \cR^+\big(\bar X\big)_{h}\big[\Bott^{-1}\big]$.
Assume we have a reference map $f \co \bar X \to B\Gamma$
which is surjective on the fundamental groups. Then the
index-difference homomorphisms $\mathrm{inddiff}_{\bar h_0}^{\Gamma}$
extend to the corresponding index-difference homomorphisms
\begin{gather}
 \mathrm{inddiff}_{\bar h_0}^{\Gamma}\co\ \pi_q\mathcal{R}^+\big(\bar X\big)\big[\Bott^{-1}\big]\to
 KO_{q+k+1}\big(C^*_r(\Gamma)\big), \nonumber
 \\[1ex]
 \mathrm{inddiff}_{\bar h_0}^{\Gamma}\co\ \pi_q\mathcal{R}^+\big(\bar X\big)_h\big[\Bott^{-1}\big]\to
 KO_{q+k+1}\big(C^*_r(\Gamma)\big).\label{eq:stable-ind}
\end{gather}
Here is one more relevant result:
\begin{Theorem}[see {\cite[Theorem~B]{ER-W2}} and~\cite{Bugg}]
\label{thm:stable-homot}
If $\Gamma$ is a torsion-free group satisfying the Baum--Connes
conjecture, and $f\co \bar X \to B\Gamma$ is split surjective on the
fundamental groups, then the homomorphisms~\eqref{eq:stable-ind}
are surjective for all $q\geq 0$.
\end{Theorem}
Now we are ready to apply those results to the case of
$(L,G)$-pseudomanifolds. Let $M_\Sigma$ be an $(L,G)$-fibered compact
pseudomanifold. We~notice that if the space $\cR_w^+(M_\Sigma)$
is not empty and $g_0\in \cR_w^+ (M_\Sigma)$ is a
base point, then it determines base points in the corresponding spaces:
the metrics $g_{\beta M,0}\in \cR^+ (\beta M)$,
$g_{\p M,0}\in \cR^+ (\p M)$ and $g_{M,0}\in \cR^+ (M)_{g_{\p M,0}}$.

As we noticed above the spaces $\cR_w^+ (M_\Sigma)_{g_{\beta M,0}}$
and $\cR^+ (M)_{g_{\p M,0}}$ from~\eqref{pull-back2} are
homotopy equivalent and the groups
$\pi_q\big(\cR_w^+(M_\Sigma)_{g_{\beta M,0}}\big)$ in the exact sequence
\eqref{eq:sgor-exact} could be replaced by~$\pi_q\big(\cR^+(M)_{g_{\p M,0}}\big)$. Thus we obtain the homomorphism
\begin{gather*}
 \mathrm{inddiff}_{g_{M,0}}\co\ \pi_q\cR_w^+ (M_\Sigma)_{g_{\beta M,0}}\xrightarrow{\cong}
\pi_q\big(\cR^+(M)_{g_{\p M,0}}\big)\xrightarrow{\mathrm{inddiff}_{g_{\p M,0} }} KO_{q+n+1},
\end{gather*}
which, together with the index-difference homomorphism
\begin{gather*}
\mathrm{inddiff}_{g_{\beta M,0} }\co\ \pi_q ( \cR^+ (\beta M))\to KO_{q+n-\ell},
\end{gather*}
determines the homomorphism
\begin{gather*}
\mathrm{inddiff}_{g_0}\co\ \pi_q (\cR^+_w(M_\Sigma))\xrightarrow{\mathrm{inddiff}_{g_{\p M,0} }\oplus\mathrm{inddiff}_{g_{\beta M,0} }}
 KO_{q+n+1}\oplus KO_{q+n-\ell}.
\end{gather*}
\begin{Corollary}\label{cor:homot-groups2}
Let $M_\Sigma$ be a $(L,G)$-fibered compact pseudomanifold with $L$ a
simply connected homogeneous space of a compact semisimple Lie group,
and $n-\ell-1\geq 5$, where $\dim M=n$, $\dim L =\ell$. Let $g_0\in
\cR^+_w(M_{\Sigma})\neq \varnothing $ be a base point giving
corresponding base points, the met\-rics~$g_{\beta M,0}\in \cR^+(\beta M)$, $g_{\p M,0}\in \cR^+ (\p M)$
and $g_{M,0}\in \cR^+ (M)_{g_{\p M,0}}$.

If $M_\Sigma$ is spin and
simply connected, then we have the following commutative diagram:
\begin{gather*}\label{pull-back1}
 \begin{diagram}
 \setlength{\dgARROWLENGTH}{1.95em}
\node{0\to \pi_q\cR_w^+ (M_\Sigma)_{g_{\beta M,0}}}
 \arrow{e,t}{j_*}
 \arrow{s,l}{\mathrm{inddiff}_{g_{M,0}}}
\node{\pi_q\mathcal{R}_w^+ (M_\Sigma)}
 \arrow{e,t}{(\resS)_*}
 \arrow{s,l}{\mathrm{inddiff}_{g_0}}
\node{\pi_q\mathcal{R}^+(\beta M)\to 0}
 \arrow{s,l}{\mathrm{inddiff}_{g_{\beta M,0}}}
 \\
\node{\!\!\! 0\longrightarrow KO_{q+n+1}}
 \arrow{e}
\node{KO_{q+n+1}\oplus KO_{q+n-\ell} }
 \arrow{e}
\node{KO_{q+n-\ell}\to 0,}
\end{diagram}
\end{gather*}
where the homomorphisms $\mathrm{inddiff}_{g_{M,0}}$ and
$\mathrm{inddiff}_{g_{\beta M,0}}$ are both nontrivial whenever the
target groups are. In particular, the homomorphism
\begin{gather*}
 \mathrm{inddiff}_{g_0}\co\ \pi_q\mathcal{R}_w^+ (M_\Sigma) \to
 KO_{q+n+1}\oplus KO_{q+n-\ell}
\end{gather*}
is surjective rationally and surjective onto the torsion
of $KO_{q+n+1}\oplus KO_{q+n-\ell}$.
\end{Corollary}
\begin{proof}
We use Corollary~\ref{short-exact} to choose a splitting
\begin{gather*}
\pi_q\big(\cR_w^+ (M_\Sigma)\big)\cong \pi_q\big(\cR_w^+ (M_\Sigma)_{g_{\beta
 M}}\big)\oplus \pi_q\big(\cR^+(\beta M)\big)\cong
\pi_q\big(\cR^+(M)_{g_{\p M,0}}\big)\oplus \pi_q\big(\cR^+(\beta M)\big).
\end{gather*}
Then we can construct the homomorphism
$\mathrm{inddiff}_{g_0}$ as a direct sum
\begin{gather*}
 \mathrm{inddiff}_{g_0}=\mathrm{inddiff}_{g_{M,0}}\oplus \mathrm{inddiff}_{g_{\beta M,0}}
 \co\ \pi_q\mathcal{R}_w^+ (M_\Sigma) \to KO_{q+n+1}\oplus KO_{q+n-\ell}.
\end{gather*}
By Theorem~\ref{thm:new-index-diff}$(a)$, the individual index difference
maps in this decomposition are non-trivial whenever the target
$KO_k$ is non-zero. Since $KO_k=\bZ/2$ whenever it has torsion,
that implies that the map surjects onto the torsion.
Similarly, the part of Theorem~\ref{thm:new-index-diff}$(a)$ about
the rational groups implies rational surjectivity.
\end{proof}
Now we address the case when $M_{\Sigma}$ is not simply-connected. Let
$\theta\co \Gamma_{\beta}\to \Gamma$ be a group homomorphism as in
Section~\ref{sec:main}. We~consider a triple
$(M_{\Sigma},\xi,\xi_{\beta})$ as an object representing an element in
the bordism group $\Omega_n^{\spin,
 (L,G)\fb}\big(B\Gamma_\beta\xrightarrow{\theta}B\Gamma\big)$, i.e.,
$M_{\Sigma}=M\cup_{\p M} N(\beta M)$ comes together with the maps
$\xi\co M\to B\Gamma$, $\xi_{\beta}\co \beta M \to B\Gamma_{\beta}$
satisfying the conditions given in Definition~\ref{def:relbordism-2}.

Thus we obtain the homomorphism
\begin{gather*}
 \mathrm{inddiff}_{g_{M,0}}^{\Gamma}\co\ \pi_q\cR_w^+
 (M_\Sigma)_{g_{\beta M,0}}\xrightarrow{\cong} \pi_q\big(\cR^+(M)_{g_{\p
 M,0}}\big)\xrightarrow{\mathrm{inddiff}_{g_{\p M,0} }}
 KO_{q+n+1}\big(C^*_r(\Gamma)\big),
\end{gather*}
which, together with the index-difference homomorphism
\begin{gather*}
 \mathrm{inddiff}_{g_{\beta M,0} }^{\Gamma_{\beta}}\co\
 \pi_q ( \cR^+ (\beta M))\to KO_{q+n-\ell}\big(C^*_r(\Gamma_{\beta})\big),
\end{gather*}
 determines the homomorphism
\begin{gather*}
\mathrm{inddiff}_{g_0}^{\Gamma,\Gamma_{\beta}}\co\ \pi_q \big(\cR^+_w(M_\Sigma)\big)
\! \xrightarrow{\mathrm{inddiff}_{g_{\p M,0} }\oplus\mathrm{inddiff}_{g_{\beta M,0} }}\!
 KO_{q+n+1}\big(C^*_r(\Gamma)\big)\oplus KO_{q+n-\ell}\big(C^*_r(\Gamma_{\beta})\big).
\end{gather*}
The same argument as above and Theorem~\ref{thm:new-index-diff}$(b)$
prove the following:

\begin{Corollary}\label{cor:homot-groups3}
Let $M_\Sigma$ be a $(L,G)$-fibered compact pseudomanifold with $L$ a
simply connected homogeneous space of a compact semisimple Lie group,
and $n-\ell-1\geq 5$, where $\dim M=n$, $\dim L =\ell$. Let $g_0\in
\mathcal{R}^+_w(M_{\Sigma})\neq \varnothing $ be a base point giving
corresponding base points, the metrics
$g_{\beta M,0}\in \cR^+(\beta M)$, $g_{\p M,0}\in \cR^+ (\p M)$
and $g_{M,0}\in \cR^+ (M)_{g_{\p M,0}}$.
Let $(M_{\Sigma},\xi,\xi_{\beta})$ represent an~ele\-ment in the bordism
group
$\Omega_n^{\spin,(L,G)\fb}\big(B\Gamma_\beta\xrightarrow{\theta}B\Gamma\big)$,
then we have the following commutative diagram:
$$
\begin{diagram}
 \setlength{\dgARROWLENGTH}{0.75em}
\node{0\to \pi_q\mathcal{R}_w^+ (M_\Sigma)_{g_{\beta M,0}}}
 \arrow{e,t}{j_*}
 \arrow{s,l}{\mathrm{inddiff}_{g_{M,0}}^{\Gamma}}
\node{\pi_q\mathcal{R}_w^+ (M_\Sigma)}
 \arrow{e,t}{(\resS)_*}
 \arrow{s,l}{\mathrm{inddiff}_{g_0}^{\Gamma,\Gamma_{\beta}}}
\node{\pi_q\mathcal{R}^+(\beta M)\to 0}
 \arrow{s,l}{\mathrm{inddiff}_{g_{\beta M,0}}^{\Gamma_{\beta}}}
 \\
\node{\!\!\! 0\longrightarrow\!\! KO_{q+n+1}\big(C^*_r(\Gamma)\big)}
 \arrow{e}
\node{KO_{q+n+1}\big(C^*_r(\Gamma)\big)\!\oplus\! KO_{q+n-\ell}\big(C^*_r(\Gamma_{\beta})\big) }
 \arrow{e}
\node{KO_{q+n-\ell}\big(C^*_r(\Gamma_{\beta})\big)\!\to\! 0.}
\end{diagram}
$$
{\samepage
If, in addition, the pair $(\Gamma,\Gamma_{\beta})$ is such that
\begin{itemize}\itemsep=0pt
\item $\Gamma$ and $\Gamma_{\beta}$ both satisfy the rational
 Baum--Connes conjecture,
\item $\Gamma$ and $\Gamma_{\beta}$ are both torsion free and have
 finite rational homological dimension $d$, and
\item $q > d-2n+2\ell+1$,
\end{itemize}}
then the image of the homomorphism
\begin{gather*}
\mathrm{inddiff}_{g_0}^{\Gamma,\Gamma_{\beta}}\co\ \pi_q\mathcal{R}_w^+ (M_\Sigma)\to
KO_{q+n+1}\big(C^*_r(\Gamma)\big)\oplus KO_{q+n-\ell}\big((C^*_r(\Gamma_{\beta})\big)
\end{gather*}
generates the target $KO_{q+n+1}\big(C^*_r(\Gamma)\big)\oplus
KO_{q+n-\ell}\big(\big(C^*_r(\Gamma_{\beta})\big)\big)\otimes \mathbb{Q}$ as a
$\mathbb{Q}$-vector space.
\end{Corollary}
Let $\cR^+(M_{\Sigma})\big[\Bott^{-1}\big]$ be the homotopy colimit
\begin{gather*}
 \cR^+(M_{\Sigma})
 \xrightarrow{\times(\Bott,g_{\flat})}\cR^+(M_{\Sigma}\times \Bott)
 \xrightarrow{\times(\Bott,g_{\flat})} \cR^+(M_{\Sigma}\times
 \Bott\times \Bott) \xrightarrow{\times(\Bott,g_{\flat})}\cdots.
\end{gather*}
We have the following conclusion from Theorems~\ref{thm:stable-homot} and~\ref{thm:pi0ofbdyinjects}:

\begin{Corollary}\label{cor:stable-homot-groups}
Let $M_\Sigma$ be a $(L,G)$-fibered compact pseudomanifold with $L$ a
simply connected homogeneous space of a compact semisimple Lie group,
and $n-\ell-1\geq 5$, where $\dim M=n$, $\dim L =\ell$. Let $g_0\in
\mathcal{R}^+_w(M_{\Sigma})\neq \varnothing $ be a base point giving
corresponding base points, the metrics
$g_{\beta M,0}\in \cR^+(\beta M)$, $g_{\p M,0}\in \cR^+ (\p M)$
and $g_{M,0}\in \cR^+ (M)_{g_{\p M,0}}$.
Let $(M_{\Sigma},\xi,\xi_{\beta})$ represent an~ele\-ment in the bordism
group
$\Omega_n^{\spin,(L,G)\fb}\big(B\Gamma_\beta\xrightarrow{\theta}B\Gamma\big)$.
Assume that the pair $(\Gamma,\Gamma_{\beta})$ is such that
\begin{itemize}\itemsep=0pt
\item $\Gamma$ and $\Gamma_{\beta}$ both satisfy the
 Baum--Connes conjecture,
\item $\Gamma$ and $\Gamma_{\beta}$ are both torsion free.
\end{itemize}
Then the homomorphism
\begin{gather*}
 \mathrm{inddiff}_{g_0}^{\Gamma,\Gamma_{\beta}}\co\
 \pi_q\mathcal{R}_w^+ (M_\Sigma)\big[\Bott^{-1}\big]\to
KO_{q+n+1}\big(C^*_r(\Gamma)\big)\oplus KO_{q+n-\ell}\big((C^*_r(\Gamma_{\beta})\big)
\end{gather*}
is surjective.
\end{Corollary}

\begin{Example}
In~\cite[Section 1.1.3]{ER-W2}, the authors provided
an example of a group $\pi_0$ such that $KO_{7+k}(C^*_r(\pi_0))
\otimes \bQ$ has countably infinite dimension for each
$k\geq 0$,\footnote{Let $\mathrm{Free}_2$ be a free group on two
generators. Then the group $\pi_0$ is a free product
${\rm SB}_3*{\rm SB}_4*{\rm SB}_5*{\rm SB}_6$, where
${\rm SB}_r$ is the $r$-th \emph{Stallings--Bieri group},
the kernel of the homomorphism
$(\mathrm{Free}_2)^r\to \bZ$ sending each generator to
$1$.} as well as a $4$-di\-men\-sional closed spin manifold $X$ with
$\pi_1X =\pi_0$. Then, according to~\cite[Theorem~C]{ER-W2}, the
manifold $Y=X\times S^2$ has the property that the group
$\pi_{7+k}\mathcal{R}^+(Y)\otimes \bQ$ has countably infinite
dimension for each $k\geq 0$. Then it is easy to
construct a pseudomanifold $M_{\Sigma}$ with $L=G=\SU(2)$ and
$\beta M = Y$. Indeed, we let $M= X\times D^3 \times L$, $\p M =
X\times S^2 \times L$, $\beta M = X\times S^2$ and $M_{\Sigma} =
M\cup_{\p M} (Y\times c(L))$. Then, clearly, Theorem~\ref{thm:pi0ofbdyinjects} and Corollary~\ref{cor:homot-groups3} imply
that $\pi_{7+k} \mathcal{R}^+_w(M_{\Sigma})\otimes \bQ$ has
countably infinite dimension for each $k\ge 0$.
\end{Example}

\subsection[A direct approach to the index-difference homomorphism for (L,G)-fibered pseudomanifolds]
{A direct approach to the index-difference homomorphism \\for
 $\boldsymbol{(L,G)}$-fibered pseudomanifolds}

 Let $M_\Sigma$ be as above, thus
a $(L,G)$-fibered pseudomanifold, and let $g_0,g_1\in
\cR^+_w(M_{\Sigma})$ be two well adapted wedge psc metrics. Then the
``difference'' between $g_0$ and $g_1$ could be detected
\emph{directly} by a wedge relative index as follows. We~consider
$M_\Sigma \times [0, 1]$ as a pseudomanifold with boundary, and equip
it with a well adapted wedge metric $\bar g$ of product-type near the
boundary and equal to $g_0$ on $M_\Sigma\times \{0\}$ and to $g_1$ on
$M_\Sigma\times \{1\}$. Now we attach infinite cylinders to $M_\Sigma
\times [0, 1]$ to get a wedge metric on the associated manifold with
cylindrical ends $(M_\Sigma \times [0, 1])_{\infty}$ (which in this
case is nothing but $M_\Sigma \times \bR$), with psc metrics along the
two ends. This is a special case of the situation already encountered
in Section~\ref{subsect:psc-witt-bordism}, see in particular
Theorem~\ref{bordism-invariance}. In~particular, since the resulting
metric on $M_\Sigma \times \bR$ has psc along the cylindrical ends and
is geometric-Witt on~the compact part, there is a well defined class
$\alpha^{\text{rel}}_w(g_0, g_1)$ in $KO_{n+1}$. This class, by
definition, is~equal to the cylindrical $\alpha_w$-class of the
Atiyah--Singer operator on the manifold with cylindrical ends
$(M_\Sigma \times [0, 1])_{\infty}\equiv M_\Sigma \times \bR$ endowed
with the natural extension of the metric~$\bar g$ to a metric~$\bar
g_\infty$.The existence of this class follows immediately, as in
Theorem~\ref{bordism-invariance}, from the $b$-edge calculus developed
in~\cite{a-gr-dirac}.

If the psc metrics $g_0$ and $g_1$ are isotopic, i.e., they lie in the
same component of $\mathcal{R}_w^+ (M_\Sigma)$, then the usual
``isotopy $\Rightarrow$ concordance'' argument shows that we can
construct a~psc metric~$\bar g$ on $M_\Sigma \times [0, 1]$, which
implies that the metric $g_\infty$ is of psc. Hence the element
$\alpha^{\text{rel}}_w(g_0, g_1)$ in~$KO_{n+1}$ gives an obstruction
to an isotopy from $g_0$ to $g_1$ within $\cR^+_w (M_\Sigma)$.
Proceeding in the same way, but using the Atiyah--Singer--Mishchenko
operator instead, we also have a class
\begin{gather*}
\alpha^{\text{rel}, \pi_1
 (M_\Sigma) }_w(g_0, g_1)\in KO_{n+1}\big(C^*_{r,\mathbb{R}}( \pi_1
(M_\Sigma))\big)
\end{gather*}
which is again an obstruction to the existence of an
isotopy between $g_0$ and $g_1$.

More generally, for an $(L,G)$-fibered pseudomanifold
we can adapt the general discussion given just before
Theorem~\ref{thm:new-index-diff} and, using in a crucial way the results
explained in Section~\ref{sec:obstructions-general}, define directly the
\emph{wedge-index-difference homomorphism}
\begin{gather*}
 w\operatorname{-inddiff}_{\bar h_0}^{\Gamma}\co\
 \pi_q\mathcal{R}^+_w (M_\Sigma)\to KO_{q+k+1}\big(C^*_r(\Gamma)\big)
\end{gather*}
either by considering $S^q\times M_\Sigma$ and $D^q\times M_\Sigma$
or by employing the classifying space $\mathbf{KO}_k \big(C^*_r(\Gamma)\big)$.\\
It would be interesting to show, but we shall leave this to future
research, that
$ w\operatorname{-inddiff}_{\bar h_0}^{\Gamma}$ is equal to the index-difference
homomorphism
considered in the previous section, that is, the composition of the
isomorphism
$\pi_q \big(\cR_w^+ (M_\Sigma)_{g_{\beta M,0}}\big)\rightarrow \pi_q \big( \cR^+ (M)_{g_{\p M,0}}\big)$
with the index-difference homomorphism for the manifold with boundary $M$,
$ \mathrm{inddiff}^\Gamma_{g_{\p M,0} }\co \pi_q \big( \cR^+ (M)_{g_{\p M,0}}\big)
\to KO_{q+n+1} (C^*_r\Gamma)$.

\subsection{Rho invariants and torsion fundamental groups}
\label{subsec:rho}
We begin by introducing the relevant bordism groups:
\begin{gather*}
{\rm Pos} ^{\spin,(L,G)\fb}_n (B\Gamma)\qquad\text{and}\qquad
{\rm Pos} ^{\spin,(L,G)\fb}_n \big(B\Gamma_\beta \xrightarrow{B\theta} B\Gamma\big).
\end{gather*}
We only give the definition of the latter of these two in detail,
since ${\rm Pos} ^{\spin,(L,G)\fb}_n (B\Gamma)$ is just a special case
of ${\rm Pos} ^{\spin,(L,G)\fb}_n \big(B\Gamma_\beta \xrightarrow{B\theta} B\Gamma\big)$ when $\Gamma=\Gamma_\beta$ and $\theta$ is the identity map.
Recall that an element of the bordism group $\Omega^{\spin,(L,G)\fb}_n
\big(B\Gamma_\beta \xrightarrow{B\theta}B\Gamma\big)$ is represented by a
\emph{geometric cycle}, i.e., a tuple
\begin{gather*}
\big(M, P\xrightarrow{p}\beta M,\ \xi\co M\to B\Gamma,\ \xi_\beta\co \beta M \to
 B\Gamma_\beta\big)
\end{gather*}
with the same compatibility conditions for $\xi$, $B\theta$,
$\xi_\beta$ as in Definition~\ref{def:relbordism-2} and with
$P\xrightarrow{p}\beta M$ a~principal $G$-bundle such that $\p
M=P\times_G L$.

Similarly, the bordism relation in Definition~\ref{def:relbordism-2}
can be described as follows: two geometric cycles
\begin{gather*}
\big(M,P\xrightarrow{p}\beta M, \xi,\xi_\beta\big)\qquad\text{and}\qquad
 \big(M',P' \xrightarrow{p'}\beta M', \xi',\xi'_\beta\big)
\end{gather*}
are equivalent if there exists a manifold $W$ with corners such that
$\p W =\p^{(0)} W \cup \p^{(1)} W$, where $\p^{(1)} W$ fibers over a
manifold with boundary $\beta W$,
\begin{gather*}
\p^{(0)}W = M\sqcup -M', \qquad
\p \big(\p^{(1)} W\big) = \p M \sqcup - \p M', \qquad
\p(\beta W) = \beta M\sqcup -\beta M',
\end{gather*}
i.e., $\p^{(1)} W \co \p M\rightsquigarrow \p M'$ and $\beta W \co
\beta M\rightsquigarrow \beta M'$ are usual spin bordisms between closed
spin manifolds, and
\begin{gather*}
\p^{(0)} W \cap \p^{(1)} W = \p M \sqcup - \p
M'.
\end{gather*}
We also have a principal $G$-bundle
$\overline{P}\xrightarrow{\overline{p}} \beta W$ restricting to $P$
and $P^\prime$ over $\p (\beta W)=\beta M\sqcup -\beta M'$.

A psc-geometric cycle for
${\rm Pos}^{\spin,(L,G)\fb}_n \big(B\Gamma_\beta \xrightarrow{B\theta} B\Gamma\big)$ is
then given by a tuple
\begin{gather*}
\big(M,P\xrightarrow{p}\beta M, \xi,\xi_\beta,g_M, g_{\beta M},\nabla^p\big),
\end{gather*}
where $g_M$ and $g_{\beta M}$ are Riemannian metrics of psc and
$\nabla^p$ is a bundle-connection on $P$. Then two
 psc-geometric cycles
\begin{gather*}
\big(M,P\xrightarrow{p}\beta M, \xi,\xi_\beta,g_M, g_{\beta M},\nabla^p\big) \qquad\text{and}\qquad
\big(M',P'\xrightarrow{p'}\beta M', \xi',\xi'_\beta, g_{M'}, g_{\beta M'},\nabla^{P'}\big)
\end{gather*}
are equivalent if there exists a bordism $(W,\bar{P})$ as above
together with metrics of psc $\bar g_W$ and $\bar g_{\beta W}$
restricting to $g_M$, $g_{M'}$ and $g_{\beta M}$, $g_{\beta M'}$ on
the boundary, along with a connection $\nabla^{\bar{P}}$ restricting
to the connections $\nabla^p$ and $\nabla^{P'}$. The proof of the
following statement is straightforward and thus omitted.
\begin{Proposition}
There exist natural homomorphisms
\begin{gather}
 R^{{\rm Pos}}\co\ \, {\rm Pos} ^{\spin,(L,G)\fb}_n \big(B\Gamma_\beta \xrightarrow{B\theta} B\Gamma\big)\to {\rm Pos} ^{\spin,(L,G)\fb}_n (B\Gamma)\nonumber
 \\ \hphantom{ R^{{\rm Pos}}\co\ \, }
 \text{``forget $\Gamma_\beta$''}, \label{natural-hom-1}
 \\
i^{\rm Pos}\co\ \ \, {\rm Pos} ^{\spin}_n (B\Gamma)\to {\rm Pos} ^{\spin,(L,G)\fb}_n \big(B\Gamma_\beta \xrightarrow{B\theta} B\Gamma\big)\nonumber
\\ \hphantom{i^{\rm Pos}\co\ \ \, }
\text{``consider a closed manifold as a pseudomanifold with empty singularities''}\label{natural-hom-2},
\\
f_{\Gamma_\beta,\Gamma}\co\ {\rm Pos} ^{\spin,(L,G)\fb}_n \big(B\Gamma_\beta \xrightarrow{B\theta} B\Gamma\big)\to
\Omega^{\spin,(L,G)\fb}_n \big(B\Gamma_\beta \xrightarrow{B\theta} B\Gamma\big)\nonumber
\\ \hphantom{f_{\Gamma_\beta,\Gamma}\co\ }
\text{``forget metric and connection''}.
\label{natural-hom-3}
\end{gather}
We also have a natural map
\begin{gather*}
\pi_0 \big(\mathcal{R}_w^+ (M_\Sigma)\big)\rightarrow {\rm Pos} ^{\spin,(L,G)\fb}_n \big(B\Gamma_\beta \xrightarrow{B\theta} B\Gamma\big).
\end{gather*}
\end{Proposition}

We would like to use the technique of Botvinnik--Gilkey~\cite{Bot-Gil-MA} (and generalizations of it, such as the ones
presented in~\cite{MR2366359}) in order to detect
elements in the groups
\begin{gather*}
 {\rm Pos} ^{\spin,(L,G)\fb}_n \big(B \pi_1 (\beta M)\to B \pi_1(M_\Sigma)\big),
\end{gather*}
with $M_\Sigma$ a pseudomanifold with $(L,G)$-fibered singularities.
A relevant version of the rho invariant on a
depth-one wedge spin stratified pseudomanifold of psc was introduced
in~\cite{Piazza-Vertman}. Then, using a suitable APS index theorem on
spaces with $(L,G)$-fibered singularities, it follows from
\cite{Piazza-Vertman} that the APS rho
invariant and the Cheeger--Gromov rho invariant of such a $M_\Sigma$
are well defined and that they both define maps
$\rho_{w,{\rm APS}}\co \pi_0 \big(\cR_w^+ (M_\Sigma)\big) \to \bR$ and
$\rho_{w,{\rm CG}}\co \pi_0 \big(\cR_w^+ (M_\Sigma)\big)$ $\to \bR$ and group
homomorphisms:
\begin{gather*}
\rho_{w,{\rm APS}}\co\ {\rm Pos} ^{\spin,(L,G)\fb}_n (B\Gamma)\to \mathbb{R},\qquad
\rho_{w,{\rm CG}}\co\ {\rm Pos} ^{\spin,(L,G)\fb}_n (B\Gamma)\to \mathbb{R},
\end{gather*}
with $\Gamma=\pi_1 (M_\Sigma)$. We~shall mainly be interested in $\rho_{w,{\rm CG}}$ and
we shall denote it simply by $\rho_{w}$.
Summarizing, we have a well-defined homomorphism:
\begin{gather}\label{rho-hom-bis}
\rho_{w}\co\ {\rm Pos} ^{\spin,(L,G)\fb}_n (B\Gamma)\to \bR.
\end{gather}
Here is an example where the Cheeger--Gromov rho invariant of a wedge
space can be used in~order to show that the group
${\rm Pos} ^{\spin,(L,G)\fb}_n \big(B \pi_1 (\beta M)\to B \pi_1 (M_\Sigma)\big)$
is infinite.

\begin{Proposition}
Let $L$ be simply connected, $G$ be a compact semisimple Lie group. Let~$M_\Sigma$ be a pseudomanifold with $(L,G)$-fibered singularities and
of dimension $m=4k+3$, with $\Gamma:=\pi_1 (M_\Sigma)$. Assume that
$\Gamma$ has an element of finite order. Then the group
\begin{gather*}
{\rm Pos}^{\spin,(L,G)\fb}_m \big(B \pi_1 (\beta M)\to B \pi_1 (M_\Sigma)\big)
\end{gather*}
is infinite provided $\cR_w^+ (M_\Sigma)\not= \varnothing $.
\end{Proposition}
\begin{proof}
Let $\Gamma_\beta:=\pi_1 (\beta M)$, $\Gamma:= \pi_1 (M_\Sigma)$ and
$\theta\co \Gamma_\beta\to\Gamma$ induced by the inclusion of $\beta M$
into~$M_\Sigma$. Let $\xi_\Sigma\co M_\Sigma\to B\Gamma$ be the
classifying map for the universal cover of $M_\Sigma$; this provides
us with a map $\xi\co M\to B\Gamma$ and we also have a classifying map
$\xi_\beta\co \beta M \to B\Gamma_\beta$ satisfying the compatibility
conditions of Definition~\ref{def:relbordism-2}. As already remarked,
we have a forgetful homomorphism
\begin{gather*}
 f_{\Gamma_\beta,\Gamma}\co {\rm Pos} ^{\spin,(L,G)\fb}_m
 \big(B\Gamma_\beta\xrightarrow{B\theta}B\Gamma\big)\rightarrow
\Omega_m ^{\spin,(L,G)\fb} \big(B\Gamma_\beta\xrightarrow{B\theta}B\Gamma\big),
\end{gather*}
(we forget the metrics). We~shall consider the subset
$C(M_\Sigma,\!\xi,\!\xi_\beta)$ of the group ${\rm Pos} ^{\spin,(L,G)\fb}_m
\!\big(\!B\Gamma_\beta\allowbreak \xrightarrow{B\theta}B\Gamma\big)$ obtained by keeping
$M_\Sigma$ and the classifying maps $\xi\co M\to B\Gamma$ and
$\xi_\beta\co \beta M \to B\Gamma_\beta$ fixed and varying only the psc
metrics $g_M$, $g_{\beta M}$ and the connection $\nabla^p$.

We shall prove that under the present assumptions
\begin{itemize}\itemsep=0pt
\item $\ker f_{\Gamma_\beta,\Gamma}$ has infinite cardinality,
\item there is a free and transitive action of $\ker f_{\Gamma_\beta,\Gamma}$
on $C(M_\Sigma,\xi,\xi_\beta)$.
\end{itemize}
This will imply that $C(M_\Sigma,\xi,\xi_\beta)$ and thus ${\rm Pos}
^{\spin,(L,G)\fb}_m \big(B\Gamma_\beta\xrightarrow{B\theta}B\Gamma\big)$ has
infinite cardinality. Consider then $\ker
(f_{\Gamma_\beta,\Gamma})\subset {\rm Pos} ^{\spin,(L,G)\fb}_m
\big(B\Gamma_\beta \xrightarrow{B\theta} B\Gamma\big)$ and
$C(M_\Sigma,\xi,\xi_\beta)$. Then, using crucially our bordism
theorem, Theorem~\ref{thm:bordism2}, we can prove
as in~\cite[Proposition 2.4]{MR2366359}, that there is well defined action
\begin{gather*}
 \ker (f_{\Gamma_\beta,\Gamma})\times C(M_\Sigma,\xi,\xi_\beta)
 \to C(M_\Sigma,\xi,\xi_\beta)\,,
\end{gather*}
which associates to $x\in \ker (f_{\Gamma_\beta,\Gamma})$ and
$\big[(M_\Sigma, \xi,\xi_\beta, g_M,g_{\beta
 M})\big]$\footnote{Here we skip ``$\nabla^p$''
 from the notations.} the class
\begin{gather*}
 x+ \big[(M_\Sigma, \xi,\xi_\beta, g_M,g_{\beta M})\big].
\end{gather*}
Indeed, as $x$ is null-bordant in $\Omega ^{\spin,(L,G)\fb}_m
\big(B\Gamma_\beta\xrightarrow{B\theta}B\Gamma\big)$, the element $x+
\big[(M_\Sigma, \xi,\xi_\beta)\big]$ is bordant to $\big[(M_\Sigma,
 \xi,\xi_\beta)\big]$ in $\Omega ^{\spin,(L,G)\fb}_m
\big(B\Gamma_\beta\xrightarrow{B\theta}B\Gamma\big)$. Since the element
$x+ \big[(M_\Sigma, \xi,\xi_\beta)\big]$ is represented by a manifold with an
adapted wedge psc metric, we can use our bordism theorem and propagate
this psc metric back to $M_\Sigma$, obtaining a new adapted wedge
metric of psc, i.e., a new element in $C(M_\Sigma,\xi_\beta,\xi)$. One
proves exactly as in Proposition 2.4 in~\cite{MR2366359} that this
action is well defined, free and transitive.

{\sloppy
Consider now
$\Gamma=\mathbb{Z}_n$, the cyclic group of order $n$. It~is explained
in~\cite{MR2366359}, building on~spe\-ci\-fic examples provided in~\cite{Bot-Gil-MA}, that the Cheeger--Gromov rho invariant defines a map
$\rho$: \mbox{${\rm Pos}^{\spin}_m (B\mathbb{Z}_n)\to \mathbb{R}$} which has an
image of infinite cardinality. To prove this, one only needs to show
that $\rho$ is a homomorphism and that it is non-trivial; the property
that it is a homomorphism is a consequence of the definition of rho
invariant whereas the fact that it is non-trivial follows from the
specific examples in~\cite{Bot-Gil-MA} (lens spaces). Remark that we
have a well defined rho-homomorphism $\rho^{{\rm rel}}_w \co {\rm Pos}
^{\spin,(L,G)\fb}_m \big(B\Gamma_\beta \xrightarrow{B\theta} B\Gamma\big)\to
\mathbb{R}$, obtained by composing the homomorphism~\eqref
 {natural-hom-1} with the homomorphism~\eqref{rho-hom-bis}:
\begin{gather*}
\rho^{{\rm rel}}_w:=\rho_w \circ R^{\rm Pos}.
\end{gather*}}\noindent
Recall now that we are under the assumption that
there is an injection $j\co\mathbb{Z}_n \to \Gamma$.
Such an~injec\-tion induces homomorphisms
\begin{gather*}
Bj^{\Omega}_*\co\ \Omega ^{\spin}_m (B\mathbb{Z}_n)\to\Omega ^{\spin}_m (B\Gamma)\qquad\text{and}\qquad
Bj^{{\rm Pos}}_* \co\ {\rm Pos} ^{\spin}_m (B\mathbb{Z}_n)\to {\rm Pos} ^{\spin}_m (B\Gamma).
\end{gather*}
Consider the following diagram, where for typographic reasons
we omit the superscripts $\spin$ and $\mathrm{fb}$,\vspace{-2mm}
\begin{gather*}
\xymatrix{{\rm Pos}_m (B\mathbb{Z}_n) \ar[r]^{Bj^{{\rm Pos}}_*} \ar[d]^{f_{\mathbb{Z}_n}} &{\rm Pos}_m (B\Gamma)\ar[d]^{f_{\Gamma}}\ar[r]^(.4){i^{\rm Pos}}
 & {\rm Pos} ^{(L,G)}_m \big(B\Gamma_\beta \xrightarrow{B\theta} B\Gamma\big) \ar[d]^{f_{\Gamma_\beta,\Gamma}}\\
 \Omega_m (B\mathbb{Z}_n) \ar[r]^{Bj^\Omega_*} & \Omega_m (B\Gamma)\ar[r]^(.4){i^\Omega}
 & \Omega^{(L,G)}_m \big(B\Gamma_\beta \xrightarrow{B\theta} B\Gamma\big).}
\end{gather*}
By naturality this is a commutative diagram. Consider
$K:= \ker f_{\mathbb{Z}_n}$. We~know from
\cite{MR2366359} that $\rho |_{K}$ has an image of infinite cardinality.
For $\kappa\in K$ consider
\begin{gather*}
\Theta_\kappa := i^{\rm Pos} \circ Bj_*^{\rm Pos} (\kappa)\in {\rm Pos} ^{\spin,(L,G)\fb}_m \big(B\Gamma_\beta \xrightarrow{B\theta} B\Gamma\big).
\end{gather*}
Then, by commutativity of the diagram, we have
\begin{gather*}
f_{\Gamma_\beta,\Gamma} (\Theta_\kappa)=0.
\end{gather*}
This means that
\begin{gather*}
\big\{ \Theta_k + [M,\beta M,\xi,\xi_\beta,g], \kappa\in K\big\}\subset
C(M_\Sigma,\xi,\xi_\beta)\subset {\rm Pos}^{\spin,(L,G)\fb}_m \big(B\Gamma_\beta \xrightarrow{B\theta} B\Gamma\big).
\end{gather*}
We also have
\begin{gather*}
\rho_w^{{\rm rel}}\big( \Theta_k + [M,\beta M,\xi,\xi_\beta, g_M,g_{\beta M},\nabla^p]\big)=\rho (k) +
\rho_w [M_\Sigma,\xi_\Sigma,g]
\end{gather*}
with $g$ the wedge metric defined by $g_M$, $g_{\beta M}$ and $\nabla^p$.
Indeed:
\begin{gather*}
\rho_w^{{\rm rel}}\big( [M,\beta M,\xi,\xi_\beta, g_M,g_{\beta M},\nabla^p]\big)
:= \rho_w \big(R^{\rm Pos} \big( [M,\beta M,\xi,\xi_\beta, g_M,g_{\beta M},\nabla^p]\big)\big)
\\ \hphantom{\rho_w^{{\rm rel}}\big( [M,\beta M,\xi,\xi_\beta, g_M,g_{\beta M},\nabla^p]\big){:}}
=
\rho_w [M_\Sigma,\xi_\Sigma,g],
\end{gather*}
whereas, thanks to Lemma 2.2 in~\cite{MR2366359} and naturality, we have
\begin{gather*}
\rho_w^{{\rm rel}}( \Theta_k)= \rho_w^{{\rm rel}}\big( i^{\rm Pos} \circ Bj_*^{\rm Pos} (\kappa)\big)=
\rho_w \big(R^{\rm Pos} \circ i^{\rm Pos} \circ Bj_*^{\rm Pos} (\kappa)\big)=\rho \big(Bj_*^{\rm Pos} (\kappa)\big)=\rho (\kappa).
\end{gather*}
This implies that the set
$\{ \Theta_k + [M,\beta M,\xi,\xi_\beta,g], \kappa\in K\}$ and thus
the group
\begin{gather*}
 {\rm Pos} ^{\spin,(L,G)\fb}_m \big(B\Gamma_\beta
 \xrightarrow{B\theta} B\Gamma\big)
\end{gather*}
that contains it, is of infinite cardinality.
\end{proof}
\begin{Corollary}
Under the previous assumptions we 
have
\begin{gather*}
{}\big| \pi_0 \big(\mathcal{R}_w^+ (M_\Sigma)\big) \big|=\infty.
 \end{gather*}
 \end{Corollary}
\begin{proof}
It suffices to observe that there exists a surjection
$\pi_0 \big(\mathcal{R}_w^+ (M_\Sigma)\big)\rightarrow C(M_\Sigma,\xi,\xi_\beta). $
\end{proof}
The above result is just an example of how rho-invariants can be used
in order to detect elements in bordism groups of wedge psc metrics or in
$\pi_0 \big(\cR_w^+ (M_\Sigma)\big)$. It~should also be possible to use
\emph{higher} rho invariants to get sharper results (under
additional assumptions
on the fundamental groups). We~comment on this in the next section.

\section{Open problems and subjects for future study}
In this section we just sketch a few ideas for additional
projects related to the topics of this paper.

\subsection[The index-difference homomorphism for L-fibered pseudomanifolds]
{The index-difference homomorphism for $\boldsymbol L$-fibered pseudomanifolds}
In this paper, we have not done very much (beyond the obstruction theory
in Section~\ref{sec:obstructions-general}) for~$L$-fi\-be\-red pseudomanifolds
that are not $(L,G)$-fibered. But there is hope of extending the
theory of the index-difference homomorphism to the more general $L$-fibered
case. Assume that $g_0,g_1\in \cR^+_w(M_{\Sigma})$,
with $M_\Sigma$ an $L$-fibered pseudomanifold.
Let us take a path of wedge metrics joining~$g_0$ and~$g_1$ and let us denote
by $g_0^v $ and $g_1^v$ the vertical part of the metric on the link-bundle. It~is not automatic, in this generality, that the path we take joining $g_0$
and $g_1$ will stay within the metrics that are geometric-Witt or psc-Witt,
i.e., with positive scalar curvature along the links.
Under which assumptions can we ensure this? One sufficient condition is that
 $\cR^+ (\p M /\beta M)$, the space of
vertical metrics of psc along the links, is connected, and,
moreover, that we can lift a path joining $g_0^v $ and $g_1^v$
in $\cR^+ (\p M /\beta M)$ to a path joining $g_0$ and $g_1$ in
$\cR^+_w(M_{\Sigma})$. When this is the case, the theory of the
index difference outlined above in Section~\ref{subsect:index-diff}
ought to go through; of course, it will be necessary to work out the details.

\subsection[The theory of R-groups]{The theory of $\boldsymbol R$-groups}
In~\cite{StolzR}, Stolz sketched the theory of concordance groups of
psc metrics (on compact smooth manifolds), which he called $R$-groups.
Assuming for simplicity that we are considering spin manifolds with
fundamental group $\Gamma$, one obtains a long exact sequence
\cite[sequence~(4.4)]{StolzR}
\begin{gather}
\label{eq:Stolzseq}
\cdots
\to R_{n+1}^\spin(B\Gamma) \xrightarrow{\p} \text{Pos}_n^\spin(B\Gamma)
\to \Omega_n^\spin(B\Gamma) \to R_n^\spin(B\Gamma) \xrightarrow{\p} \cdots.
\end{gather}
Here $\text{Pos}_n^\spin(B\Gamma)$ is the bordism group of pairs
$(M,g)$, with $M$ a closed spin manifold with a map to $B\Gamma$
and $g$ a {\psc} metric on $M$. The map from this group to~$\Omega_n^\spin(B\Gamma)$ simply forgets the Riemannian metric.
The group $R_n^\spin(B\Gamma)$~\cite[Definition 4.1]{StolzR}
is a~bordism group of
pairs $(M\to B\Gamma, h)$, with $M$ a compact spin $n$-manifold with
possibly non-empty boundary, and with $h$ a psc metric on $\partial M$.
The map $\p\co R_{n+1}^\spin(B\Gamma) \to \text{Pos}_n^\spin(B\Gamma)$
is just restriction to the boundary.

It should be rather straightforward to replace $\text{Pos}_n^\spin(B\Gamma)$
here with $\text{Pos}_n^{\spin,(L,G)\fb}(B\Gamma)$, and~$\Omega_n^\spin(B\Gamma)$ with $\Omega_n^{\spin,(L,G)\fb}(B\Gamma)$,
to obtain $R$-groups of pseudomanifolds with $(L,G)$-fibered singularities:
$R_n^{\spin,(L,G)\fb} (B\Gamma)$.
\begin{itemize}\itemsep=0pt
\item
 What do these $R$-groups look like and how do they differ from the usual
 $R$-groups?
\item Do they act freely and transitively on the set of concordance classes
 of psc wedge metrics, as in Stolz~\cite{StolzR}?
\item Do we have the analogue of Stolz sequence~\eqref{eq:Stolzseq}
 in the singular setting?
\item If so, can one map the analogue of sequence~\eqref{eq:Stolzseq} to a Higson--Roe type analytic exact sequence
as in~\cite{MR3286895}? This would involve the definition of a higher
rho class associated to a~wedge psc metric. For the signature operator
on stratified spaces that are Witt or~Cheeger, a similar result has been proved
in~\cite{ap-surgery}, mapping the Browder--Quinn surgery sequence to
the Higson--Roe analytic surgery sequence,
and we expect the analysis there to play a role here.
\item Is there an $R$-group
 $R_n^{\spin,(L,G)\fb} \big(B\Gamma_\beta\xrightarrow{B\theta} B\Gamma\big)$
and does it fit into a Stolz sequence
\begin{align*}
\cdots&\to\Omega_n^{\spin,(L,G)\fb} \big(B\Gamma_\beta\xrightarrow{B\theta} B\Gamma\big)\to
R_n^{\spin,(L,G)\fb} \big(B\Gamma_\beta\xrightarrow{B\theta} B\Gamma\big)
\\
&\to{\rm Pos}_{n-1}^{\spin,(L,G)\fb} \big(B\Gamma_\beta\xrightarrow{B\theta} B\Gamma\big)\to \cdots\, ?
\end{align*}
\item Can we use the Stolz $(L,G)$-sequences above, the mapping to
 Higson--Roe, and the techniques of Xie--Yu--Zeidler~\cite{XieYuZeidler}
 in order to give a lower bound on the rank of the group
 $\text{Pos}_n^{\spin,(L,G)\fb}(B\Gamma)$?
\end{itemize}

\subsection{Rho-invariants of wedge psc metrics}
\subsubsection{Torsion-free groups}
In the smooth case one can prove that if $\pi_1 (M)$ is torsion-free and
the Baum--Connes map
$\mu_{{\rm max}}\co K_* (B\Gamma) \to K_* (C^*_{{\rm max}} \Gamma)$
is an isomorphism,
then the Cheeger--Gromov rho invariant of a~{\psc} metric on a spin manifold
$M$ of odd dimension must vanish.
See~\cite{MR2294190} and also~\cite{BRoy} for a new proof based on the
Higson--Roe analytic surgery sequence.
Notice that this is a no-go result: the Cheeger--Gromov rho invariant for
manifolds with fundamental group
satisfying these conditions \emph{cannot} be used in order to
distinguish metrics of psc (indeed, it is equal to $0$).

Can one extend this vanishing result to $(L,G)$-pseudomanifolds and the
Cheeger--Gromov rho invariant of a wedge metric?

A direct approach would build on the proof of Piazza--Schick in~\cite{MR2294190}. A different route would use the results~of the previous Section~\ref{subsec:rho},
the Benameur--Roy map $K_0 (D^*_\Gamma)\to \mathbb{R}$ of~\cite{BRoy}
and the exactness of the Higson--Roe analytic surgery sequence.

\subsubsection{Groups with torsion}
On the other hand, one would like to define and use \emph{higher} rho
invariants
in order to distinguish wedge metrics of positive scalar curvature, especially
for fundamental groups with torsion.
See for example~\cite{azzali-wahl-2cocycle,Le-Pi-psc,Xie-Yu-moduli}
for the case of smooth closed manifolds.

\subsection{Stratifications with higher depth}
For applications to algebraic varieties, moduli spaces, and other
natural examples of stratified spaces, it would be nice to generalize
our theory to pseudomanifolds with higher depth, in other words, with
singular strata of multiple dimensions. While the analytic part
concerning Dirac operators is largely under control, thanks to
\cite{a-gr-dirac} (see also~\cite{PZ}), on the geometric side we
expect several complications:
\begin{itemize}\itemsep=0pt
\item The geometry of the tubular neighborhoods of the singular strata
 gets to be rather complicated, as one would need to consider
 iterated Riemannian submersions and careful curvature estimates.
\item Proving the appropriate bordism theorem would be quite complicated.
\item Perhaps one could say more in the case of Baas--Sullivan singularities,
 as in~\cite{MR1857524}, which did consider some pseudomanifolds with
 higher depth.
\end{itemize}

\subsection{Topological questions}
\begin{itemize}\itemsep=0pt
\item
 The ``mixed fundamental group'' bordism group
 $\Omega_n^{(L,G)} \big(B\Gamma_\beta\xrightarrow{B\theta} B\Gamma\big)$
 seems to be very hard to compute; it doesn't just fit into an exact
 sequence like the one we developed in~\cite{BB-PP-JR}.
 But is there a \emph{spectral sequence} computing
 $\Omega_n^{(L,G)} \big(B\Gamma_\beta\xrightarrow{B\theta} B\Gamma\big)$?
\item What about a Kreck--Stolz $s$-invariant? More precisely, in~\cite[Proposition 2.13]{MR1205446}, Kreck and Stolz defined an
 invariant $s(M,g)$
 of psc metrics $g$ on closed spin manifolds $M$ of~dimen\-sion~$3$
 mod $4$ with vanishing Pontryagin classes. This invariant has the
 properties that it is preserved under spin-structure-preserving
 isometries and that the relative index $i(g, g')$ in~the
 sense of~\cite[Section~4]{MR720933} is given by $s(M,g)-s(M,g')$.
 It would be interesting
 to know if one could do something similar for
 adapted psc wedge metrics on Witt spin pseudomanifolds $M_\Sigma$.
 This would require using the index theorem of~\cite{Piazza-Vertman}
 and replacing vanishing of the Pontryagin classes with vanishing
 of the Goresky--MacPherson homology $L$-class or, more generally,
 of the $K$-homology class of the pseudomanifold signature operator.
\end{itemize}

\appendix
\section[On the spin-bordism invariance of the alpha-invariant]
{On the spin-bordism invariance of the $\boldsymbol{\alpha}$-invariant}\label{sec:appendix}

In this appendix we shall discuss briefly an analytic proof of the spin-bordism invariance of the $\alpha$ class and of
the $\alpha^\Gamma$ class. We~have already observed, building
 on the case of the signature operator
on Witt spaces treated in great detail in~\cite[Section 7]{almp-package}, that given an analytic proof in the smooth case, it extends mutatis mutandis to the
pseudomanifold case.

\subsection[A generator of KO1(R)]{A generator of $\boldsymbol{KO_1 (\bR)}$}
Consider $\mathbb{R}$ with its standard spin structure. Then
we have $\Di_{\mathbb{R}}$ and its class in
$KO_1 (\mathbb{R}):= KKO(C_0 (\mathbb{R}),C\ell_1)$.
\begin{Proposition}
$KO_1 (\mathbb{R})$
 is isomorphic to $\mathbb{Z}$ and
the class
$[\Di_{\mathbb{R}}]\in KO_1 (\mathbb{R})$ is a generator. Equivalently,
$KO_1 ((0,1))$
 is isomorphic to $\mathbb{Z}$ and
the class
$\big[\Di_{(0,1)}\big]\in KO_1 ((0,1))$ is a generator.
\end{Proposition}

\begin{proof}
 This is a very special case of~\cite[Theorems 4.8 and 4.10]{MR918241}
 and of~\cite[Section~5, Lem\-ma~4]{MR1388299}.
\end{proof}

\subsection{Boundary of Atiyah--Singer is Atiyah--Singer}
We now know that $[\Di_{\mathbb{R}}]$ is the generator of $KO_1 (\mathbb{R})$. We~want to prove that prove that
``boundary of Atiyah--Singer is Atiyah--Singer'' in $KO$-homology, i.e.,
\begin{gather*}
\delta [\Di_X]=[\Di_{\partial X}]\qquad\text{in}\quad KO_{\dim \partial X} (\partial X)
\end{gather*}
on a spin manifold $X$ with metrically collared boundary $\partial X$. Here $\delta$ is the connecting homomorphism in KO-theory
associated to the semisplit short exact sequence
\begin{gather*}
0\to C_{\partial X} (X)\to C (X)\to C(\partial X)\to 0
\end{gather*}
with $C_{\partial X} (X)$ the ideal of continuous functions vanishing
on the boundary.

Using excision we can easily reduce to a collar neighbourhood.
Then the result is a consequence of the following:

\begin{Lemma}
Let $X=[0,1)\times N$ with $\p X=\{0\}\times N$. Let $n=\dim N$.
In this case
$C_{\partial X} (X)=C_0 \big((0,1)\times N\big)$. Let us consider the connecting homomorphism
\begin{gather*}
KO_{n+1} \big((0,1)\times N)\big)\xrightarrow{\delta} KO_n (N)
\end{gather*}
associated to the semisplit short exact sequence
\begin{gather*}
0\to C_0 \big((0,1)\times N\big)\to C_0 \big([0,1)\times N\big)\to C(N)\to 0.
\end{gather*}
Then $\delta$ is the inverse of the isomorphism
\begin{gather*}
KO_n (N)\to KO_{n+1} \big((0,1)\times N\big)
\end{gather*}
obtained by taking the Kasparov product with $\big[\Di_{(0,1)}\big]$, the generator of $KO_1 ((0,1))$.
In particular, as $[\Di_{(0,1)\times N}]$ is precisely the Kasparov product of
$\big[\Di_{(0,1)}\big]$ and $[\Di_N]$, we have that $\delta [\Di_{(0,1)\times N}]=[\Di_N]$
which is what we wanted to prove.
\end{Lemma}

\begin{proof} By the naturality of Kasparov product we can reduce to the
 connecting homomorphism
\begin{gather*}
KO_1 ((0,1))\xrightarrow{\delta} KO_0 ({\rm point})
\end{gather*}
associated to the semisplit short exact sequence
\begin{gather*}
0\to C_0 ((0,1))\to C_0 ([0,1))\to \bR\to 0.
\end{gather*}
But $\delta$ is the inverse of the
isomorphism $KO_0 ({\rm point})\to KO_1 ((0,1))$ obtained by taking
the Kasparov product with the generator $\big[\Di_{(0,1)}\big]$ of $KO_1 ((0,1))$.
This again follows from Kasparov's results on Poincar\'e duality
in $KO$-theory,~\cite[Theorems~4.8 and~4.10]{MR918241}
and~\cite[Section~5, Lemma~4]{MR1388299}.
\end{proof}

\subsection{Spin bordism invariance}
 Using $\delta [\Di_X]=[\Di_{\partial X}]$ we can now prove the spin bordism invariance of the $\alpha$-class
 in $KO_n$.

\begin{Theorem} Assume that $Y$ is spin and of dimension $n$ and that $Y=\partial X$, with $X$ spin. Then
$\alpha(Y)=0\in KO_n$.
\end{Theorem}

\begin{proof} We write part of the long exact sequence in $KO$ associated to the
 short exact sequence
\begin{gather*}
0\to C_{\partial X} (X)\to C (X)\to C(\partial X)\to 0.
\end{gather*}
We have the induced exact sequence
\begin{gather*}
KO_{n+1} (X,\partial X)\xrightarrow{\delta} KO_n (\partial X)\xrightarrow{\iota_*}
KO_n (X)
\end{gather*}
with $\iota$ the inclusion of $ \partial X$ into $X$. Notice,
in particular,
that
\begin{gather*}
\iota_* \circ \delta =0.
\end{gather*}
Consider $\pi^Y\co Y\to {\rm point}$; obviously $\pi^Y=\pi^X\circ\iota$. Then we have
\begin{gather*}
\alpha(Y)=\pi^Y_* [\Di_Y]=\pi^X_* \circ \iota_* [\Di_Y]
= \pi^X_* \circ \iota_* \circ \delta [\Di_X]=0.\tag*{\qed}
\end{gather*}\renewcommand{\qed}{}
\end{proof}

\noindent
The same proof applies to the $\alpha^\Gamma$ class if we use
the Atiyah--Singer--Mishchenko operator
$\Di_{{\rm MF}}$ and the exact sequence
\begin{gather*}
KK_{n+1}O\big( C_{\partial X} (X), C^*_r \Gamma\big) \xrightarrow{\delta}
KK_n O\big( C({\partial X}), C^*_r \Gamma\big) \xrightarrow{\iota_*}
KK_n O\big( C_{\partial X} (X), C^*_r \Gamma\big) .
\end{gather*}
See~\cite[Proposition 2.3]{Leichtnam-Piazza-GAFA} for the details.

Note that for the
$\alpha$ class the spin bordism invariance is usually proved by
identifying it with the Atiyah--Milnor--Singer
invariant. See~\cite{lawson89:_spin}.
The generalization of this approach to the $\alpha^\Gamma$ class is
treated in~\cite{MR866507} and (more extensively) in~\cite{Rosenberg-3}.

\subsubsection*{Acknowledgements}

{\sloppy
We thank the Mathematisches Forschungsinstitut Oberwolfach
for hosting Workshop 1732 in~2017 on
Analysis, Geometry and Topology of Positive Scalar
Curvature Metrics, which marked the start of this project. This work
was also supported by U.S.\ NSF grant number
DMS-1607162, by Simons Foundation Collaboration Grant number 708183,
by Sapienza Universit\`a di Roma, and by the {\it
 Ministero Istruzione Universit\`a e Ricerca} through the PRIN {\it
 Spazi di Moduli e Teoria di Lie}. B.B.\ and J.R.\ acknowledge a
very pleasant and productive visit to Rome in May--June 2019, as well
as a visit by J.R.\ to Rome in January 2020. P.P.\ thanks Pierre
Albin and Jesse Gell-Redman for interesting discussions about the
content of Section~\ref{sec:obstructions-general}.

}

We would like to thank the referees of this article for
a careful reading of the previous drafts and for suggestions for
improvements.\vspace{-2mm}

\pdfbookmark[1]{References}{ref}
\LastPageEnding

\end{document}